\documentclass[12pt]{amsart}
\usepackage{amscd,amssymb}
\usepackage{upref}
\usepackage{pstcol,pst-3d}

\psset{unit=0.7cm,linewidth=0.4pt,arrowsize=1.5pt 4}
\def\vertex{\pscircle[fillstyle=solid,fillcolor=black]{0.0714}}
\newpsstyle{fatline}{linewidth=1pt}
\newpsstyle{fyp}{fillstyle=solid,fillcolor=verylight}
\definecolor{verylight}{gray}{0.95}
\definecolor{light}{gray}{0.9}
\definecolor{medium}{gray}{0.85}

\let\oldlabel=\label
\def\prellabel{\marginparsep=1em
    \def\label##1{\oldlabel{##1}\ifmmode\else\ifinner\else
         \marginpar{{\footnotesize\ \\ \tt
                    ##1}}\fi\fi}}

\def\Y{{\textup{\textsf{Y}}}}
\def\bF{{\textbf{\textup{F}}}}
\def\bG{{\textbf{\textup{G}}}}
\def\bR{{\textbf{\textup{R}}}}
\def\bS{{\textbf{\textup{S}}}}
\def\bT{{\textbf{\textup{T}}}}

\let\frak=\mathfrak
\let\cal\mathcal
\let\Bbb\mathbb

\let\la=\langle
\let\ra=\rangle

\def\cd#1{$\textup{CD}_{#1}$}

\def\merc{\lambda}
\def\regcyc{\,\lozenge\,}

\def\ggr.aut{{\operatorname{gr.aut}}}
\def\comp{\operatorname{comp}}
\def\L{\operatorname{L}}

\def\int{{\operatorname{int}}}
\def\inv{{\operatorname{inv}}}
\def\vert{{\operatorname{vert}}}
\def\conv{\operatorname{conv}}
\def\het{\operatorname{ht}}
\def\E{\operatorname{E}}
\def\V{\operatorname{V}}
\def\st{\operatorname{St}}
\def\St{\operatorname{{\Bbb S}t}}

\def\path{\operatorname{path}}
\def\B{\operatorname{B}}
\def\X{\operatorname{X}}
\def\T{\operatorname{\bf T}}

\def\Aff{\operatorname{Aff}}
\def\End{\operatorname{End}}
\def\gp{\operatorname{gp}}
\def\rank{\operatorname{rank}}
\def\Ker{{\operatorname{Ker}}}
\def\Coker{{\operatorname{Coker}}}
\def\Proj{\operatorname{Proj}}
\def\Im{\operatorname{Im}}
\def\Hom{\operatorname{Hom}}
\def\Pol{\operatorname{Pol}}
\def\GL{\operatorname{GL}}
\def\Qu{\operatorname{Q}}
\def\Col{\operatorname{Col}}
\def\Vect{\operatorname{Vect}}

\def\Q{\qedsymbol\kern1pt}
\def\sqq#1#2{{\hbox{\kern 0.5pt\vbox{\hbox{\kern#1ex\vrule width #2pt height #1ex}%
     \hrule height #2pt}\kern 1.5pt}}}
\def\sq{{\mathchoice{\sqq1{0.5}}{\sqq1{0.5}}{\sqq{0.9}{0.4}}
                    {\sqq{0.7}{0.3}}}}

\def\vt{{\kern1pt|}}
\def\RR{{\mathbb R}}

\def\QQ{{\mathbb Q}}
\def\ZZ{{\mathbb Z}}
\def\NN{{\mathbb N}}
\def\TT{{\mathbb T}}
\def\AA{{\mathbb A}}
\def\EE{{\mathbb E}}
\def\FF{{\mathbb F}}
\def\GG{{\mathbb G}}
\def\vv{{\mathbb V}}
\def\fP{{\mathfrak P}}
\def\cG{{\mathcal G}}

\def\cB{{\mathcal B}}

\let\epsilon=\varepsilon
\let\phi=\varphi
\let\theta=\vartheta

\advance\headheight1.15pt

\newtheorem{lemma}{Lemma}[section]
\newtheorem{corollary}[lemma]{Corollary}
\newtheorem{theorem}[lemma]{Theorem}
\newtheorem{proposition}[lemma]{Proposition}

\theoremstyle{definition}
\newtheorem{definition}[lemma]{Definition}
\newtheorem{remark}[lemma]{Remark}
\newtheorem{example}[lemma]{Example}
\newtheorem{question}[lemma]{Question}

\begin{document}

\title[Higher polyhedral $K$-groups]{Higher polyhedral $K$-groups}

\author{Winfried Bruns \and Joseph Gubeladze}
\address{Universit\"at Osnabr\"uck,
FB Mathematik/Informatik, 49069 Osnabr\"uck, Germany}
\email{Winfried.Bruns@mathematik.uni-osnabrueck.de}
\address{A. Razmadze Mathematical Institute, Alexidze St. 1, 380093
Tbilisi, Georgia} \email{gubel@rmi.acnet.ge, soso@math.sfsu.edu}

\thanks{The second author was supported by the Deutsche
Forschungsgemeinschaft, INTAS grant 99-00817 and TMR grant ERB
FMRX CT-97-0107}


\begin{abstract}
We define higher \emph{polyhedral} $K$-groups for commutative
rings, starting from the stable groups of elementary automorphisms
of polyhedral algebras. Both Volodin's theory and Quillen's +
construction are developed. In the special case of algebras
associated with unit simplices one recovers the usual algebraic
$K$-groups, while the general case of lattice polytopes reveals
many new aspects, governed by polyhedral geometry. This paper is a
continuation of \cite{BrG5} which is devoted to the study of
polyhedral aspects of the classical Steinberg relations.  The
present work explores the polyhedral geometry behind Suslin's well
known proof of the coincidence of the classical Volodin's and
Quillen's theories. We also determine all $K$-groups coming from
$2$-dimensional polytopes.
\end{abstract}

\maketitle

\section{Introduction}\label{INTR}

In \cite{BrG5} we have initiated a `polyhedrization' of algebraic
$K$-theory of commutative rings $R$ that is based on the groups of
graded automorphisms of \emph{polyhedral} algebras. The algebra
associated with a unit simplex is just a polynomial ring over $R$,
its automorphism group is the general linear group, and the
resulting theory is nothing else but usual $K$-theory (details in
Subsection \ref{UNIMODULAR}). However, for general lattice
polytopes (and, further, for lattice polyhedral complexes) many
new aspects show up.

The motivation behind such a polyhedrization can be summarized as
follows. The geometry of the affine space $\AA^d_k$ and polyhedral
geometry merge naturally in the concept of toric variety. Linear
algebra constitutes part of the geometry of affine spaces, and it
admits its own polyhedrization resulting in the study of the
category $\Pol(k)$ ($k$ a field) of \emph{polytopal algebras} and
their graded homomorphisms. The objects of $\Pol(k)$ are
essentially the homogeneous coordinate rings of projective toric
varieties, and they are naturally associated with lattice
polytopes. This category contains the category $\Vect(k)$ of
(finite dimensional) vector spaces over $k$ as a full subcategory
and, despite being non-additive, reveals surprising similarities
with $\Vect(k)$ \cite{BrG4}. This leads to \emph{polyhedral}
linear algebra, where the $\Hom$-objects are no longer vector
spaces, but certain $k$-varieties equipped with $k^*$-equivariant
structures. Our general philosophy is that essentially all
standard linear algebra facts should have meaningful geometric
analogues in $\Pol(k)$.

Lower $K$-theory generalizes linear algebra to the study of
projective modules over general rings and the automorphism groups
of free modules. By analogy, in the category $\Pol(k)$ one studies
algebra retractions and the groups of automorphisms. The objects
belong to the area of toric varieties and, pursuing this analogy
(see, \cite{BrG1}--\cite{BrG4}), we have gained new insight into
well known results and proved several new ones, especially for
automorphism groups.

The next natural step in merging $K$-theory and polyhedral
geometry was carried out in \cite{BrG5} where we investigated the
Schur multiplier of the stable group of elementary automorphisms
of polytopal algebras. The elementary automorphisms play the same
r\^ole in $\Pol(k)$ as elementary matrices in $\Vect(k)$. For
detailed comments and explanations see Section \ref{POLYTOPES}.

In this paper we develop higher polyhedral $K$-theory for
commutative rings. From the very beginning it is clear that it
cannot be based on additive structures (say, exact categories).
There are essentially two classical constructions of higher
$K$-theory, which can be adapted to groups of algebra
automorphisms, namely Quillen's + construction \cite{Qu1} and
Volodin's theory \cite{Vo}. We will explore them both, having in
view their equivalence in the classical situation. While
\cite{BrG5} has been devoted to the polyhedral aspects of the
classical Steinberg relations, the present work explores the
polyhedral geometry behind Suslin's well known proof of the
coincidence of the Volodin and Quillen theories \cite{Su1,Su2}.

The analysis of the polyhedral Steinberg relations leads one to an
appropriate class of lattice polytopes which we have called
\emph{balanced}. For them one can proceed similarly to the
development of Milnor's $K_2$. For higher $K$-groups one has to
investigate the extent to which lattice polytopes and their
\emph{column structures} support the standard $K$-theoretical
constructions, and this requirement narrows the class of suitable
polytopes somewhat further.

In Suslin's approach \cite{Su1}, \cite{Su2} the original Volodin
theory is defined in terms of \emph{triangular} linear groups
based on abstract finite posets. In the polyhedral interpretation
these posets correspond to oriented graphs formed by edges of unit
simplices. For a general polytope one considers \emph{column
vectors} which, usually, penetrate the interior of the polytope.
These column vectors are roots of the corresponding linear group,
namely the group of graded automorphisms of the polytopal algebra,
but, as explained in \cite{BrG5}, we are working with essentially
non-reductive linear groups.

The main results can be summarized as:

(a) Generalization of the triangular subgroups of the general
linear groups to subgroups of the stable group of elementary
automorphisms, supported by \emph{rigid systems of column
vectors}. Their structure is studied in Section \ref{TRIANG}.
Thereafter, in Section \ref{HPK}, we construct polyhedral versions
of Volodin's and  Quillen's theories. For the index $2$ they
coincide with the polyhedral version of Milnor's $K_2$.

(b) The homology acyclicity of the space naturally associated to
the Volodin theory. This is shown by a suitable polyhedrization of
Suslin's arguments, where several technical difficulties have to
be overcome.

(c) Introduction of the class of {\it $\Col$-divisible} polytopes.
For polytopes of this class the theories of Quillen and Volodin
agree. In view of (b) this amounts to showing the homotopy
simplicity of the relevant spaces. However, the situation is more
complicated than in the classical situation \cite{Su1}. The class
of $\Col$-divisible polytopes includes all balanced polygons
($2$-dimensional polytopes).

(d) Complete computation of all polygonal theories (i.~e.\ those
corresponding to lattice polygons). This is accomplished with use
of Quillen's \cite{Qu2} and Nesterenko-Suslin's \cite{NSu}
homological computations in linear groups.

Actually two versions of Volodin's theory are defined, one of
which is of an auxiliary character -- it provides the acyclicity
result mentioned in (b) above. Then we prove the coincidence of
all the three theories for $\Col$-divisible polytopes.

Sections \ref{POLYTOPES}--\ref{HPK} are devoted to the development
of the basic notions, including the definition of the polytopal
analogues of triangular subgroups and the description of the
corresponding $K$-theories. The $K$-theoretic calculations are
contained in Sections \ref{ACYCL}--\ref{POLYG}.

The new $K$-groups are actually bifunctors in two covariant
arguments -- a (commutative) ring and a (balanced) polytope. In
other words, \cite{BrG5} and the present paper add a polyhedral
argument to $K$-theory, showing up in $K_i$ for $i\ge 2$. One is
naturally lead to the following two questions:
\begin{enumerate}
\item Are the new groups computable to the same extent as the
usual $K$-groups?

\item Is there a polyhedral $K$-theory of schemes?
\end{enumerate}

In connection with the first question we mention that it is not
yet clear how far the polyhedral $K$-groups are from the usual
ones. Our expectation, supported by the results in Section
\ref{POLYG}, is that the polyhedral theory for $\Col$-divisible
polytopes is very close to the usual one, namely a finite direct
sum of copies. Such a claim should be thought of as a higher
stable analogue of the main result of \cite{BrG1} (see Theorem
\ref{PLg} below): higher syzygies between elementary automorphisms
for $\Col$-divisible polytopes come from unit simplices.
Computations that make such results possible are provided by the +
construction approach.

However, for general balanced polytopes the theories may well
diverge substantially. The simplest hypothetical candidate for
such a divergence is the unit pyramid over the unit square,
showing up several times in \cite{BrG5} and the present paper. But
the corresponding computations remain to be elusive.

In connection with the second question we remark that a more
conceptual approach might lead to the theory of \emph{toric
bundles} over schemes for which the fibers are no longer affine
spaces but general affine toric varieties, equipped with an
algebraic action of the 1-dimensional torus. The morphisms of such
bundles are defined by the requirement that fiberwise they induce
morphisms of algebraic varieties, respecting the torus action. One
could also try to use `Jouanolou's device' \cite{J} based on
appropriate affinization.

Finally, thanks to the results in \cite{BrG5} and the structural
description of triangular groups (Theorems \ref{strutri} and
\ref{strutri1}), we can treat elementary automorphisms in a formal
way, i.~e.\ without referring to their action on polytopal
algebras. Hence the rather combinatorial flavor of the exposition.
However, the paper is really about higher syzygies of elementary
automorphisms. For this reason we are restricted to commutative
rings of coefficients -- for non-commutative rings such
automorphisms simply do not exist.

As remarked, all the rings that appear below are commutative.

\section{Polytopes, their algebras, and their linear groups}\label{POLYTOPES}

Per suggestion of the referee -- to whom we are indeed very
grateful for a number of valuable comments, corrections, and
suggestions to make the paper more readable -- we begin this
section with some basic notions related to convex polytopes.

\subsection{General polytopes}\label{GENPOL}
The interested reader may consult Ziegler \cite{Z} for the
arguments skipped in this section.

By a \emph{polytope} $P\subset\RR^n$, $n\in\NN$, we always mean a
\emph{finite convex} polytope, i.~e.\ $P$ is the convex hull of a
finite subset $\{x_1,\ldots,x_k\}\subset\RR^n$:
$$
P=\conv(x_1,\ldots,x_k):=\{a_1x_1+\dots+a_kx_k\ :\ 0\le
a_1,\dots,a_k\le1,\\ a_1+\cdots+a_k=1\}.
$$

Alternatively, a polytope is a bounded subset of $\RR^n$ that can
be represented as the intersection of a finite system of closed
affine halfspaces,
$$
P=\bigcap_{j=1}^s\mathcal H_j,\quad \mathcal H_j=\{x\in\RR^n\ :\
\mathcal L_j(x)\ge b_j\},
$$
where $\mathcal L_j:\RR^n\to\RR$ is a linear mapping and
$b_j\in\RR$, $j=1,\ldots,s$.

Polytopes of dimension 1 are called \emph{segments} and those of
dimension 2 are called \emph{polygons}.

The \emph{affine hull} $\Aff(X)$ of a subset $X\subset\RR^n$ is
the smallest affine subspace of $\RR^n$ containing $X$, i.~e.
\begin{align*}
\Aff(X)=\{a_1x_1+\cdots+a_kx_k\ :\ k\in\NN,\ x_1,\ldots,x_k\in X,\
a_1,\ldots,a_k\in\RR,\\ a_1+\cdots+a_k=1\}.
\end{align*}
If $\dim \Aff(X)=k$ for a subset $X=\{x_1,\dots,x_k\}$ of
cardinality $k$, then $x_1,\dots,x_k$ are \emph{affinely
independent} and the polytope $P=\conv(x_1,\dots,x_k)$ is called a
\emph{simplex}.

For a halfspace $\mathcal H\subset\RR^n$ containing $P$, the
intersection $P\cap\partial\mathcal H$ of $P$ with the boundary
affine hyperplane $\partial\mathcal H$ of $\mathcal H$ is called a
\emph{face} of $P$. The polytope itself is also considered as a
face.

The faces of $P$ are themselves polytopes. Faces of dimension $0$
are \emph{vertices} and those of codimension $1$ (i.~e.\ of
dimension $\dim P-1$) are called \emph{facets}. A polytope is the
convex hull of the set $\vert(P)$ of its vertices. If $\dim P=n$,
then there is a unique halfspace $\mathcal H$ for each facet
$F\subset P$ such that $\partial\mathcal H\cap P=F$.

\subsection{Lattice polytopes}\label{LATPOL} A polytope $P\subset\RR^n$ is
called a \emph{lattice polytope} if the vertices of $P$ belong to
the integral lattice $\ZZ^n$. More generally, a lattice in $\RR^n$
is a subset $\cG=x_0+\cG_0$ with $x_0\in\RR^n$ and an  additive
subgroup $\cG_0$ generated by $n$ linearly independent vectors. A
polytope $P$ with $\vert(P)\subset \cG$ is called a $\cG$-polytope
if the vertices of $P$ belong $\cG$. However, since all the
properties of $\cG$-polytopes we are interested in remain
invariant under an affine automorphism of $\RR^n$ mapping $\cG$ to
$\ZZ^n$, we can always assume that our polytopes have vertices in
$\ZZ^n$. More generally, lattice polytopes $P$ and $Q$ that are
isomorphic under an integral-affine equivalence of $\Aff(P)$ and
$\Aff(Q)$ are equivalent objects or our theory. We simply speak of
\emph{integral-affinely equivalent polytopes}.

Faces of a lattice polytope are again lattice polytopes.

For a lattice polytope $P\subset\RR^n$ we put $\L_P=P\cap\ZZ^n$. A
simplex $\Delta$  is called \emph{unimodular} if
$\sum_{z\in\vert(\Delta)}\ZZ(z-z_0)$ is a direct summand of
$\ZZ^n$ for some (equivalently, every) vertex $z_0$ of $\Delta$.
All unimodular simplices of dimension $n$ are integral-affinely
equivalent. Such a simplex is denoted by $\Delta_n$ and called a
\emph{unit $n$-simplex}. Standard realizations of $\Delta_n$ are
$\conv(O,e_1,\dots,e_n)\subset\RR^n$ or
$\conv(e_1,\dots,e_{n+1})\subset\RR^{n+1}$. ($e_i$ is the $i$th
unit vector.)

There is no loss in assuming that a given lattice polytope $P$ is
full dimensional (i.~e.\ $\dim P=n$) and that $\ZZ^n$ is the
smallest affine lattice containing $\L_P$.  In fact, we choose
$\Aff(P)$ as the space in which $P$ is embedded and fix a point
$x_0\in\L_P$ as the origin. Then the lattice
$x_0+\sum_{x\in\L_P}\ZZ(x-x_0)$ can be identified with $\ZZ^r$,
$r=\dim P$.

Under this assumption let $F$ be a facet of $P$ and choose a point
$z_0\in F$. Then the subgroup
$$
F_\ZZ:=(-z_0+\Aff(F))\cap\ZZ^n\subset\ZZ^n
$$ is isomorphic to $\ZZ^{n-1}$. Moreover, there is a unique group
homomorphism $\la F,-\ra:\ZZ^n\to\ZZ$, written as $x\mapsto\la F,x\ra$,
such that $\Ker(\la F,-\ra)=F_\ZZ$, $\Coker(\la F,-\ra)=0$, and on the set
$\L_P$, $\la F,-\ra$ attains its minimum $b_F$  at the lattice points of
$F$.

The $\ZZ$-linear form $\la F,-\ra$ can be extended in a unique way
to a linear function on $\RR^n$. The description of $P$ as an
intersection of halfspaces yields that $x\in P$ if and only if
$\la F,x\ra\ge b_F$ for all facets $F$ of $P$.

Our blanket assumption throughout the paper is: \emph{all
polytopes, considered below, are lattice polytopes.}

\subsection{Column structures}\label{COLSTR} 
\begin{figure}[htb]
\begin{center}
\begin{pspicture}(0,0)(5,3)
\pspolygon[style=fyp](0,0)(5,0)(5,2)(4,3)(2,3)(1,2)
\psline[linecolor=gray](2,0)(2,3)
\psline[linecolor=gray](3,0)(3,3)
\psline[linecolor=gray](4,0)(4,3)
\psline[linecolor=gray](1,0)(1,2) \multirput(0,0)(1,0){6}{\vertex}
\multirput(1,1)(1,0){5}{\vertex} \multirput(1,2)(1,0){5}{\vertex}
\multirput(2,3)(1,0){3}{\vertex} \rput(1.3,0.5){$v$}
\psline[style=fatline]{->}(1,1)(1,0.05)
\end{pspicture}
\end{center}
\caption{A column structure}\label{FigCol}
\end{figure}
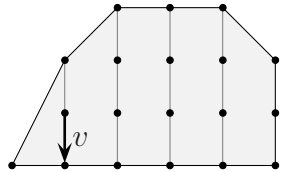

Let $P\subset\RR^n$ be a
polytope. A nonzero element $v\in\ZZ^n$ is called a \emph{a column
vector} for $P$ if there exists a facet $F\subset P$ such that
$x+v\in P$ whenever $x\in\L_P\setminus F$. In this situation $F$
is uniquely determined and called the \emph{base facet} of $v$. We
use the notation $F=P_v$. The set of column vectors of $P$ is
denoted by $\Col(P)$. A \emph{column structure} is a pair of type
$(P,v)$, $v\in\Col(P)$. Figure \ref{FigCol} gives an example of a
column structure.
Familiar examples of column structures are the unit simplices
$\Delta_n$ with their edge vectors (i.~e. the vectors
$z'-z''\in\ZZ^n$ for vertices $z'$ and $z''$ of $\Delta_n$) and
the unit square (the convex hull of the set
$\{(0,0),(1,0)(0,1)(1,1)\}$) with its edge vectors.

Using the description of $P$ in terms of the functions $\la
F,-\ra$ it is not hard to see that $v$ is a column vector of $P$
if and only if there exists exactly one facet $F$ with $\la
F,v\ra=-1$ and $\la G,v\ra\ge 0$ for all facets $G\neq F$.

\subsection{Polytopal semigroups and their rings}\label{POLRINGS}

To a polytope $P\subset\RR^n$ one associates the additive
subsemigroup $S_P\subset\ZZ^{n+1}$, generated by $\{(z,1)\ :\
z\in\L_P\}\subset\ZZ^{n+1}$. Let $C_P\subset\RR^{n+1}$ be the cone
$\{az\ :\ a\in\RR_+,\ z\in P\}$. Then $C_P$ is the convex hull of
$S_P$. It is a \emph{finite rational pointed} cone. In other
words, $C_P$ is the intersection of a finite system of halfspaces
in $\RR^{n+1}$ whose boundaries are rational hyperplanes
containing the origin $O\in\RR^{n+1}$, and there is no affine line
contained in $C_P$.

As in Subsection \ref{LATPOL}, there is no loss of generality in
assuming that $\ZZ^n$ is the lattice spanned affinely by $\L_P$ in
$\RR^n$. This is equivalent to $\gp(S_P)=\ZZ^{n+1}$, and the
condition $S_P=C_P\cap\ZZ^{n+1}$ on the polytope $P$ is known as
the \emph{normality condition} \cite{BrGTr}. This is equivalent to
saying that the affine semigroup ring $k[S_P]$ is normal for some
(equivalently, every) field $k$.

All segments and polygons are normal, but in dimensions $\ge3$
this is no longer the case (the interested reader is referred to
\cite{BrGTr} and \cite{BrG6} for the detailed theory).

It is an easy observation that the normality of $P$ is equivalent
to the normality of the facet $F$ if there exists a column
structure with base facet $F$.

While the points $x\in \L_P$ are identified with
$(x,1)\in\ZZ^{n+1}$, a column vector $v$ is to be identified with
$(v,0)\in\ZZ^{n+1}$.

Let $F$ be a facet of $F$. We use the function $\la F,-\ra$ to
define the \emph{height} of
$x=(x',x'')\in\RR^{n+1}=\RR^n\times\RR$ above the hyperplane
$\mathcal H$ through the facet $C_F$ of $C_P$ by setting
$$
\het_F(x)=\la F,x'\ra-x''b_F.
$$
For lattice points $x$ the function $\het_F$ counts the number of
hyperplanes between $\mathcal H$ and $x$ (in the direction of $P$)
that are parallel to, but different from $\mathcal H$ and pass
through lattice points. If $v$ is a column vector, then $\het_v$
stands for $\het_{P_v}$. Moreover, we are justified in calling
$\het_F(v,0)=\la F,v\ra$ the \emph{height of $v$ with respect to
$F$}, since $v$ is identified with $(v,0)$.

Although the semigroup $S_P$ may miss some integral points in the
cone $C_P$ this cannot happen on the segments parallel to a column
vector $v$. More precisely, the following holds:
\begin{equation}\label{COLUMNS}
z+v\in S_P\ \text{for all}\ z\in S_P\setminus C_{P_v}.
\end{equation}
($C_{P_v}\subset C_P$ is the face subcone, corresponding to
$P_v$.)

Let $R$ be a ring and $P\subset\RR^n$ a lattice polytope. The
semigroup ring $R[P]:=R[S_P]$ -- the \emph{polytopal $R$-algebra
of $P$} -- carries a graded structure $R[P]=R\oplus
R_1\oplus\cdots$ in which $\deg(x)=1$ for all $x\in\L_P$. By
definition of $S_P$ it follows that $R_1$ generates $R[P]$ over
$R$.

We are interested in the group $\ggr.aut_R(P)$ of graded
$R$-algebra automorphisms of $R[P]$. For a field $R=k$ the group
$\ggr.aut_k(P)$ is naturally a $k$-linear group. In fact, it is a
closed subgroup of $\GL_m(k)$, $m=\#\L_P$. We call $\ggr.aut_k(P)$
the \emph{polytopal $k$-linear group} of $P$. Its structure will
be given in Theorem \ref{PLg}.

In the special case when $P$ is a unimodular simplex the ring
$R[P]$ is isomorphic to a polynomial algebra $R[X_1,\ldots,X_m]$,
$m=\#\L_P$. Therefore, the category $\Pol(R)$ of polytopal
$R$-algebras and graded homomorphisms between them contains a full
subcategory that is equivalent to the category of free
$R$-modules. From this perspective $\Pol(R)$ is a `polytopal
extension' of the category of free $R$-modules and one might
wonder to which extent the basic $K$-theoretic facts generalize
from the smaller category to $\Pol(R)$. As mentioned in the
introduction, this direction of investigation is pursued in
\cite{BrG4}.

However, the motivation for the current paper is somewhat
different (and less categorial). This we explain next.

\subsection{Polytopal linear groups}\label{PLG}

Assume $R$ is a ring and $P$ a polytope. Let $(P,v)$ be a column
structure and $\lambda\in R$. As pointed out above, we identify
the vector $v$ with the degree $0$ element $(v,0)\in\ZZ^{n+1}$,
and further with the corresponding monomial in $R[\ZZ^{n+1}]$.
Then we define a mapping from $S_P$ to $R[\ZZ^{n+1}]$ by the
assignment
$$
x\mapsto (1+\lambda v)^{\het_v x}.
$$
Since $\het_v$ is a group homomorphism $\ZZ^{n+1}\to\ZZ$, our
mapping is a homomorphism from $S_P$ to the multiplicative monoid
of $R[\ZZ^{n+1}]$. Now it is immediate from (\ref{COLUMNS}) in
Subsection \ref{POLRINGS} that the (isomorphic) image of $S_P$
lies actually in $R[P]$. Hence this mapping gives rise to a graded
$R$-algebra endomorphism $e_v^\lambda$ of $R[P]$. But then
$e_v^\lambda$ is actually a graded automorphism of $R[P]$ because
$e_v^{-\lambda}$ is its inverse.

Here is an alternative description of $e_v^\lambda$.  By a
suitable integral change of coordinates we may assume that
$v=(0,-1,0,\ldots,0)$ and that $P_v$ lies in the subspace
$\RR^{n-1}$ (thus $P$ is in the upper halfspace). Now consider the
standard unimodular $n$-simplex $\Delta_n$ with vertices at the
origin and standard coordinate vectors:
$\Delta_n=\conv(O,(1,\dots,0),\dots,(0,\dots,1))$. It is clear
that there is a sufficiently large natural number $c$, such that
$P$ is contained in a parallel translate of $c\Delta_n$ by a
vector from $\ZZ^{n-1}$.  Let $\Delta$ denote such a parallel
translate. Then we have a graded $R$-algebra embedding
$R[P]\subset R[\Delta]$. Moreover, $R[\Delta]$ can be identified
with the $c$-th Veronese subring of the polynomial ring
$R[X_0,\ldots,X_n]$ in such a way that $v=X_0/X_1$. Now the
automorphism of $R[X_0,\ldots,X_n]$ mapping $X_1$ to $X_1+\lambda
X_0$ and leaving all the other variables invariant induces an
automorphism $\alpha$ of the subalgebra $R[\Delta]$, and $\alpha$
in turn can be restricted to an automorphism of $R[P]$, which is
nothing else but $e_v^\lambda$.

It is clear from this description of $e_v^\lambda$ that it becomes
an elementary matrix ($e_{01}^\lambda$ in our notation) in the
special case when $P=\Delta_n$, after the identification
$\ggr.aut_R(P)=\GL_{n+1}(R)$. Accordingly, the automorphisms of
type $e_v^\lambda$ are called \emph{elementary}, and the group
they generate in $\ggr.aut_R(P)$ is denoted by $\EE_R(P)$.

In this way we have generalized the basic building blocks of
higher $K$-theory of rings to the polytopal setting: general
linear groups and their elementary subgroups. Actually, the real
motivation for us to pursue the analogy has been the main result
of \cite{BrG1}. It is the polytopal extension of the fact that an
invertible matrix over a field can be diagonalized by elementary
transformations on rows (or columns) -- or, putting it in
different words, the group $SK_1$ is trivial for fields.

\begin{lemma}\label{afemb}
Let $R$ be a ring, $P$ a polytope, and $v_1,\ldots,v_s$ pairwise
different column vectors for $P$ with the same base facet
$F=P_{v_i}$, $i=1,\dots,s$.  Then the mapping
$$
\phi:(R,+)^s\to\ggr.aut_R(P),\qquad
(\lambda_1,\ldots,\lambda_s)\mapsto
e_{v_1}^{\lambda_1}\circ\cdots\circ e_{v_s}^{\lambda_s},
$$
is an embedding of groups. In particular, $e_{v_i}^{\lambda_i}$
and $e_{v_j}^{\lambda_j}$ commute for all $i,j\in\{1,\dots,s\}$,
and the inverse of $e_{v_i}^{\lambda_i}$ is
$e_{v_i}^{-\lambda_i}$.

In the special case, when $R$ is a field the homomorphism $\phi$
is an injective homomorphisms of algebraic groups.
\end{lemma}

This lemma is proved in \cite[Lemma 3.1]{BrG1} for fields $R=k$
(where we use the notation $\Gamma_k(P)$ for $\ggr.aut_k(P)$), but
the general case makes absolutely no difference.

The image of the embedding $\phi$ given by Lemma \ref{afemb} is
denoted by $\AA(F)$. Of course, $\AA(F)$ may consist only of the
identity map of $R[P]$, namely if there is no column vector with
base facet $F$. In the case in which $P$ is the unit simplex and
$R[P]$ is the polynomial ring, $\AA(F)$ is the subgroup of all
matrices in $\GL_{\dim P+1}(R)$ that differ from the identity
matrix only in the non-diagonal entries of a fixed column.

For the rest of this subsection we assume that $k$ is a field and
set $n=\dim P$.

After $\AA(F)$ we introduce some further subgroups of
$\ggr.aut_k(P)$. First, the $(n+1)$-torus $\TT_{n+1}=(k^*)^{n+1}$
acts naturally on $k[P]$ by restriction of its action on
$k[\ZZ^{n+1}]$ that is given by
$$
(\xi_1,\ldots,\xi_{n+1})(e_i)=\xi_ie_i,\ i\in[1,n+1];
$$
here $e_i$ is the $i$-th standard basis vector of $\ZZ^{n+1}$.
This gives rise to an algebraic embedding
$\TT_{n+1}\subset\ggr.aut_k(P)$, and we will identify $\TT_{n+1}$
with its image. It consists precisely of those automorphisms of
$k[P]$ that multiply each monomial by a scalar from $k^*$.

Second, the automorphism group $\Sigma(P)$ of the semigroup $S_P$
is in a natural way a finite subgroup of $\ggr.aut_k(P)$. It is
the group of integral affine transformations mapping $P$ onto
itself.

Third, we have to consider a subgroup of $\Sigma(P)$ defined as
follows. Assume $v$ and $-v$ are both column vectors.  Then for
every point $x\in P\cap\ZZ^n$ there is a unique $y\in P\cap\ZZ^n$
such that $\het_v(x,1)=\het_{-v}(y,1)$ and $x-y$ is parallel to
$v$. The mapping $x\mapsto y$ gives rise to a semigroup
automorphism of $S_P$: it `inverts columns' that are parallel to
$v$. It is easy to see that these automorphisms generate a normal
subgroup of $\Sigma(P)$, which we denote by $\Sigma(P)_{\inv}$.

Finally, $\Col(P)$ is the set of column structures on $P$. Now the
main result of \cite{BrG1} is:

\begin{theorem}\label{PLg}
Let $P$ be an $n$-dimensional polytope and $k$ a field.
\begin{itemize}
\item[(a)]
Every element $\gamma\in\ggr.aut_k(P)$ has a (not uniquely
determined) presentation
$$
\gamma=\alpha_1\circ\alpha_2\circ\cdots\circ\alpha_r\circ\tau\circ\sigma,
$$
where $\sigma\in\Sigma(P)$, $\tau\in\TT_{n+1}$, and
$\alpha_i\in\AA(F_{i})$ such that the facets $F_i$ are pairwise
different and $\#(F_i\cap\ZZ^n)\le \#(F_{i+1}\cap\ZZ^n)$,
$i\in[1,r-1]$.
\item[(b)]
For an infinite field $k$ the connected component of unity
$\ggr.aut_k(P)^0\subset\ggr.aut_k(P)$ is generated by the
subgroups $\AA(F_i)$ and $\TT_{n+1}$. It consists precisely of
those graded automorphisms of $k[P]$ which induce the identity map
on the divisor class group of the normalization of $k[P]$.
\item[(c)] $\dim \ggr.aut_k(P)=\#\Col(P)+n+1$.
\item[(d)]
One has $\ggr.aut_k(P)^0\cap\Sigma(P)=\Sigma(P)_\inv$ and
$$
\ggr.aut_k(P)/\ggr.aut_k(P)^0\approx\Sigma(P)/\Sigma(P)_\inv.
$$
Furthermore, if $k$ is infinite, then $\TT_{n+1}$ is a maximal
torus of $\ggr.aut_k(P)$.
\end{itemize}
\end{theorem}

\section{Stable groups of elementary automorphisms and Polyhedral $K_2$}\label{MILN}

\subsection{Product of column vectors}\label{PRODUCTS}
The notion of product of two column vectors $u,v\in\Col(P)$ has
been introduced in \cite[Definition 3.2]{BrG5}: we say that the
product $uv$ exists if $u+v\neq 0$ and for every point
$x\in\L_P\setminus P_u$ the condition $x+u\notin P_v$ holds. In
this case, we define the product as $uv=u+v$.

\begin{figure}[htb]
\begin{center}
 \psset{unit=1cm}
 \def\vertex{\pscircle[fillstyle=solid,fillcolor=black]{0.05}}
\begin{pspicture}(-0.3,-0.5)(4,2.5)
 \pspolygon[style=fyp,linecolor=darkgray](0,0)(3,0)(3,2)(2,2)
 \footnotesize
 \multirput(0,1)(1,0){4}{\vertex}
 \multirput(0,0)(1,0){4}{\vertex}
 \multirput(0,2)(1,0){4}{\vertex}
 \psline[style=fatline]{->}(1,1)(0,0)
 \psline[style=fatline]{->}(1,0)(0,0)
 \psline[style=fatline]{->}(2,1)(3,1)
 \psline[style=fatline]{->}(1,1)(1,0)
 \psline[style=fatline]{->}(3,2)(2,1)
 \psline[style=fatline]{->}(3,2)(3,1)
 \rput(-0.2,0.6){$w=uv$}
 \rput(1.3,0.5){$u$}
 \rput(0.6,-0.2){$v$}
 \rput(2.5,0.8){$-v$}
 \rput(2.1,1.5){$w$}
 \rput(4.0,1.5){$u=w(-v)$}
\end{pspicture}
\end{center}
\caption{The product of two column vectors} \label{ProdCol}
\end{figure}
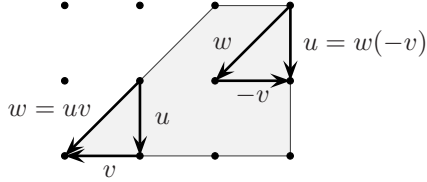

Figure \ref{ProdCol} shows a polytope with all its column vectors
and the two existing products $w=uv$ and $u=w(-v)$.

In the case of a unimodular simplex the product of two oriented
edges, viewed as column vectors, exists if and only if they are
not opposite to each other and the end point of the first edge is
the initial point of the second edge.

The basic properties of such products are given by Proposition 3.3
and Corollary 3.4 in \cite{BrG5}. For easier reference we
summarize these properties in the following

\begin{proposition}\label{maincrit}
Let $P$ be a polytope and $u,v,w\in\Col(P)$. Then:
\begin{itemize}
\item[(a)]
$uv$ exists if and only if $u+v\not=0$ and $\la P_v,u\ra>0$,
\item[(b)]
$uv$ exists if and only if $u+v\in\Col(P)$ and $P_{u+v}=P_u$,
\item[(c)]
if $uv$ exists then $P_{uv}=P_u$ and $v$ is parallel to $P_u$
(i.~e.\ $\la P_u,v\ra=0$),
\item[(d)]
$u+v\in\Col(P)$ if and only if exactly one of the products $uv$
and $vu$ exists,
\item[(e)]
if both $uv$ and $vw$ exist and $u+v+w\neq 0$ then the products
$(uv)w$ and $u(vw)$ also exist and, clearly, $(uv)w=u(vw)$,
\item[(f)]
If $vw$ and $u(vw)$ exist and $u+v\neq 0$ then $uv$ also exists.
On the other hand, the existence of $uv$ and $(uv)w$ does in
general not imply the existence of $vw$, even if $v+w\neq 0$.
\item[(g)]
The following are equivalent for $v\in\Col(P)$:
\begin{enumerate}
\item $-v\in\Col(P)$;
\item there exist facets $F,G$ of $P$ such
that $\la F,v\ra=-1$, $\la G,v\ra=1$, and $\la H,v\ra=0$ for every
facet $H\neq F,G$;
\item there exists $w\in\Col(P)$ and facets
$F,G$ such that $\la F,v\ra=-1$, $\la G,v\ra>0$, $\la F,w\ra>0$,
$\la G,w\ra=-1$.
\end{enumerate}
\item[(h)]
If $w=uv$ and $-w\in\Col(P)$, then $-u,-v\in\Col(P)$ as well.
\end{itemize}
\end{proposition}

The basic observation, already mentioned, for the proof of the
proposition is that $v\in\ZZ^n$ belongs to $\Col(P)$ if and only
if there exists $F\in\FF(P)$ such that
$$
\la F,v\ra=-1\quad\text{and}\quad\la G,v\ra\geq0\ \text{for all
}G\in\FF(P),\ G\neq F.
$$

\subsection{Balanced polytopes}\label{BALANCED}

A polytope $P$ is called \emph{balanced} if $\la P_u,v\ra\leq 1$
for all $u,v\in\Col(P)$. One easily observes that $P$ is balanced
if and only if $|\la P_u,v\ra|\leq 1$ for all $u,v\in\Col(P)$.

The reason we introduce balanced polytopes is that the main result
of \cite{BrG5} is only proved for this class of polytopes.
However, it is not yet excluded that everything generalizes to
arbitrary polytopes.

In order to help the reader visualize the class of balanced
polytopes we recall the classification result in dimension 2. It
uses the notion of projective equivalence: $n$-dimensional
polytopes $P,Q\subset\RR^n$ are called \emph{projectively
equivalent} if and only if $P$ and $Q$ have the same dimension,
the same combinatorial type, and the faces of $P$ are parallel
translates of the corresponding ones of $Q$. An alternative
definition in terms of \emph{normal fans} is given in Subsection
\ref{FUNCT}.

Recall the notation
$\Delta_n=\conv(O,(1,\ldots,0),\ldots,(0,\ldots,1))$ for the unit
$n$-simplex.

\begin{theorem}\label{DIM2}
For a balanced polygon $P$ there are exactly the following
possibilities (up to integral-affine equivalence):
\begin{itemize}
\item[(a)]
$P$ is a multiple of the unimodular triangle $P_a=\Delta_2$. Hence
$\Col(P)= \{\pm u,\pm v,\allowbreak\pm w\}$ and the column vectors
are subject to the obvious relations,
\item[(b)]
$P$ is projectively equivalent to the trapezoid
$P_b=\conv\bigl((0,0),(0,2),(1,1),\allowbreak(0,1)\bigr)$, hence
$\Col(P)=\{u,\pm v,w\}$ and the relations in $\Col(P)$ are $uv=w$
and $w(-v)=u$,
\item[(c)]
$\Col(P)=\{u,v,w\}$ and $uv=w$ is the only relation,
\item[(d)]
$\Col(P)$ has any prescribed number of column vectors, they all
have the same base edge (clearly, there are no relations between
them),
\item[(e)]
$P$ is projectively equivalent to the unit lattice square $P_e$,
hence $\Col(P)=\{\pm u,\pm v\}$ with no relations between the
column vectors,
\item[(f)]
$\Col(P)=\{u,v\}$ so that $P_u\neq P_v$ with no relations in
$\Col(P)$.
\end{itemize}
\end{theorem}

It turns out that polyhedral $K$-groups are invariants of the
projective equivalence classes of polytopes (in arbitrary
dimension); see Lemma \ref{PrEq} below.

\subsection{Doubling along a facet}\label{DOUBLING}

Let $P\subset\RR^n$ be a polytope and $F\subset P$ be a facet. For
simplicity we assume that $0\in F$, a condition that can be
satisfied by a parallel translation of $P$. Denote by
$H\subset\RR^{n+1}$ the $n$-dimensional linear subspace that
contains $F$ and whose normal vector is perpendicular to that of
$\RR^n=\RR^n\oplus0\subset\RR^{n+1}$ (with respect to the standard
scalar product on $\RR^{n+1}$). Then the upper half space
$H\cap\bigl(\RR^n\times\RR_+\bigr)$ contains a congruent copy of
$P$ which differs from $P$ by a $90^\circ$ rotation. Denote the
copy by $P^{\vt_F}$, or just by $P^\vt $ if there is no danger of
confusion.

Note that $P^\vt$ is not always a lattice polytope with respect to
the standard lattice $\ZZ^{n+1}$. However, it is so with respect
to the sublattice $(\ZZ^n)^{\vt_F}$ which is the image of $\ZZ^n$
under the $90^\circ$ rotation.

The operator of doubling along a facet is then defined by
$$
P^{\sq_F}=\conv(P,P^\vt )\subset \RR^{n+1}.
$$

The doubled polytope is a lattice polytope with respect to the
subgroup $(\ZZ^n)^{\sq_F}=\ZZ^n+(\ZZ^n)^{\vt_F}\subset\RR^{n+1}$.
After a change of basis in $\RR^{n+1}$ that does not affect
$\RR^n$ we can replace $(\ZZ^n)^{\sq_F}$ by $\ZZ^{n+1}$, and
consider $P^{\sq_F}$ as an ordinary lattice polytope in
$\RR^{n+1}$. In what follows, whenever we double a lattice
polytope $P\subset\RR^n$ along a facet $F$, the lattice of
reference in $\RR^{n+1}$ is always $\ZZ^n+(\ZZ^n)^{\vt_F}$. For
simplicity of notation this lattice will be denoted by
$\ZZ^{n+1}$.

Sometimes we will refer to $P$ as $P^-$. For an object $z$,
associated to $P$ (say, a lattice point or a column vector),
$z^\vt$ will denote the corresponding object in $P^\vt$, presuming
the facet $F$ is clear from the context.

\begin{figure}[htb]
\begin{center}
\psset{unit=1cm}
\def\vertex{\pscircle[fillstyle=solid,fillcolor=black]{0.07}}
\begin{pspicture}(-2.8,0)(1,2)
\psset{viewpoint=3.5 2 -1.5} \footnotesize \ThreeDput[normal=0 0
-1](0,0,0){
  \pspolygon[style=fyp](0,0)(3,0)(2,1.5)(0,1.5)(0,0)
  \multirput(0,0)(1,0){4}{\vertex}
  \multirput(0,1.5)(1,0){3}{\vertex}
  \psline[linewidth=1.5pt]{->}(0,1.5)(0,0) 
 }
\ThreeDput[normal=0 -1 0](0,0,0){
  \pspolygon[style=fyp](0,0)(3,0)(2,1.5)(0,1.5)(0,0)
  \multirput(0,1.5)(1,0){3}{\vertex}
  \psline[linewidth=2.0pt]{->}(0,1.5)(0,0) 
 }
\ThreeDput[normal=-1 0 0](0,0,0){
  \psline[linewidth=1.0pt]{->}(1.5,0)(0,1.5)
 }
\ThreeDput[normal=-1 0 0](2,0,0){
  \psline[linewidth=1.0pt]{->}(0,1.5)(1.5,0)
 }
 \rput(-2.2,-0.0){$P$}
 \rput(-0.25,2.0){$P^\vt $}
 \rput(-1.3,0.5){$F$}
 \rput(-1.8,1.3){$\delta^-$}
 \rput(-0.5,1.1){$\delta^+$}
 \rput(-0.6,-0.25){$v^-$}
 \rput(0.3,0.5){$v^\vt $}
\end{pspicture}
\caption{Doubling along the facet $F$} \label{Doubling}
\end{center}
\end{figure}
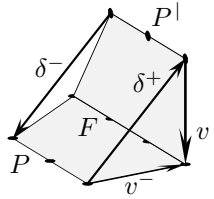

In case $F=P_v$ for some $v\in\Col(P)$ we will use the notation
$P^{\sq_F}=P^{\sq_v}$.

The polytope $P^{\sq_F}$ has two distinguished column vectors,
which are the lattice unit vectors in $\ZZ^{n+1}$ parallel to the
lines connecting the points $x^-\in\L_{P^-}\setminus F$ with the
corresponding points $x^\vt\in\L_{P^\vt}$. The column vector of
these two, which has $P^\vt$ as the base facet, will be denoted by
$\delta^+$, and $\delta^-$ will refer to the other vector. In
particular, $P^{\sq_F}_{\delta^-}=P^-$. Clearly
$\delta^-=-\delta^+$.

Let $\FF(P)$ denote the set of facets of $P$. We have the bijective mapping
$$
\Psi:\FF(P)\cup\{P\}\to\FF(P^{\sq_F})
$$
defined by
$$
\Psi(G)=
\begin{cases}\conv(G^-,G^\vt)&\text{if}\ G\in\FF(P)\setminus\{F\},\\
P^\vt&\text{if}\ G=F,\\
P&\text{if}\ G=P.
\end{cases}
$$

The following equations are easily observed:
\begin{align}
 \la\Psi(G),\delta^+\ra&=\la\Psi(G),\delta^-\ra=0 \
 \text{for all}\ G\in\FF(P)\setminus\{F\},\\
 \la P^-,\delta^+\ra&=\la P^\vt,\delta^-\ra=1,\\
 \la G,z\ra&=\la\Psi(G),z\ra \ \text{for all}\
 z\in\ZZ^n,\ G\in\FF(P).
\end{align}
(In equation (3) the pairings are considered for $P$ and $P^\sq$
respectively and $\ZZ^n$ is thought of as the subgroup
$\ZZ^n\oplus0\subset\ZZ^{n+1}$.)

\begin{lemma}\label{double}
Let $F\subset P$ be a facet and $v\in\Col(P)$. Then:
\begin{itemize}
\item[(a)]
$v\in\Col(P^{\sq_F})$,
\item[(b)]
in $\Col(P^{\sq_v})$ the following equations hold:
$$
v=v^-=\delta^+v^\vt,\qquad v^\vt=\delta^-v^-,
$$
\item[(c)]
if $P$ is balanced then $P^{\sq_F}$ is also balanced and
$$
\Col(P^{\sq_F})=\Col(P)^-\cup\Col(P)^\vt\cup\{\delta^+,\delta^-\}.
$$
\end{itemize}
\end{lemma}
See \cite[Lemmas 4.1, 5.1, Corollary 4.2]{BrG5}.

\subsection{The stable group of elementary automorphisms}\label{ELAUT}
An ascending infinite chain of lattice polytopes ${\mathfrak
P}=(P=P_0\subset P_1\subset\dots)$ is called a \emph{doubling
spectrum} if the following conditions hold:
\begin{itemize}
\item[(i)]
for every $i\in\ZZ_+$ there exists a column vector
$v\subset\Col(P_i)$ such that $P_{i+1}=P_i^{\sq_v}$,
\item[(ii)]
for every $i\in\ZZ_+$ and any $v\in\Col(P_i)$ there is an index
$j\geq i$ such that $P_{j+1}=P_j^{\sq_v}$.
\end{itemize}

Here we use the inclusion $\Col(P_i)\subset\Col(P_{i+1})$, a
consequence of Lemma \ref{double}(a).

One says that $v\in\Col(P_i)$ is \emph{decomposed} at the $j$th
step in ${\frak  P}$ for some $j\geq i$ if $P_{j+1}=P_j^{\sq_v}$.
By \cite[Lemma 7.2]{BrG5} one has

\begin{lemma}\label{infdec}
Every column vector, showing up in a doubling spectrum, gets
decomposed infinitely many times.
\end{lemma}

Associated to a doubling spectrum $\frak P$ is the `infinite
polytopal' algebra
$$
R[\fP]=\lim_{i\to\infty}R[P_i]
$$
and the filtered union
$$
\Col(\fP)=\lim_{i\to\infty}\Col(P_i).
$$
The product of two vectors from $\Col(\fP)$ is defined in the
obvious way, using the definition for a single polytope. Also, we
can speak of systems of elements of $\Col(\fP)$ having the same
base facets, etc.

Elements $v\in\Col(\fP)$ and $\lambda\in R$ give rise to a graded
automorphism of $R[\fP]$ as follows: we choose an index $i$ big
enough so that $v\in\Col(P_i)$. Then the elementary automorphisms
$e_v^{\lambda}\in\EE_R(P_j)$, $j\geq i$, constitute a compatible
system and, therefore, define a graded automorphism of $R[{\frak
P}]$. This automorphism will also be called `elementary' and it
will be denoted by $e_v^{\lambda}$.

The group $\EE(R,\fP)$ is by definition the subgroup of
$\ggr.aut_R(R[\fP])$, generated by all elementary automorphisms.

The next result comprises Propositions 7.3 and 7.4 and Theorem 7.6
from \cite{BrG5}.

\begin{theorem}\label{elaut}
Let $R$ be a ring and $P$ be a polytope (not necessarily balanced)
admitting a column structure. Assume $\fP=(P\subset P_1\subset
P_2\subset\dots)$ is a doubling spectrum. Then:
\begin{itemize}
\item[(a)] $\EE(R,\fP)$ is naturally isomorphic to
$\EE(R,{\mathfrak Q})$ for any other doubling spectrum ${\mathfrak
Q}=(P\subset Q_1\subset Q_2\subset\dots)$.

\item[(b)] $\EE(R,\fP)$ is perfect.

\item[(c)] The center of $\EE(R,\fP)$ is trivial.

\item[(d)] $e_u^{\lambda}\circ e_u^{\mu}=e_u^{\lambda+\mu}$ for
every $u\in\Col(\fP)$ and $\lambda,\mu\in R$. \item[(e)] If $P$ is
balanced, $u,v\in\Col(\fP)$, $u+v\neq 0$ and $\lambda,\mu\in R$
then
$$
[e_u^\lambda,e_v^\mu]=
\begin{cases}e_{uv}^{-\lambda\mu}&\text{if}\ uv\ \text{exists},\\
1&\text{if}\ u+v\notin\Col(\fP).
\end{cases}
$$
\end{itemize}
\end{theorem}

The difficult parts of this theorem are the claims (c) and (e),
which in the special case $P=\Delta_n$ are just standard facts.

Thanks to Theorem \ref{elaut}(a) we can use the notation
$\EE(R,P)$ for $\EE(R,\fP)$.

\begin{remark}\label{steinberg}
Theorem \ref{elaut}(e) is the generalization of Steinberg's
relations between elementary matrices to balanced polytopes. In
order to find the classical Steinberg relation $[e_{ij}^\lambda,
e_{jk}^\mu]=e_{ik}^{\lambda\mu}$ in this equality one must observe
that in our setting the configuration $e_{ij}e_{jk}$ corresponds
to the existence of $vu$ if we associate with $e_{ij}$ the column
vector $\varepsilon_i-\varepsilon_j$ where
$\varepsilon_1,\dots,\varepsilon_n$ are the vectors of the
canonical basis of $\RR^n$, and simultaneously the vertices of
$\Delta_{n-1}$.

That we associate $\varepsilon_i-\varepsilon_j$ with $e_{ij}$ (and
not $\varepsilon_j-\varepsilon_i$) is forced by our notation in
which we add column vectors on the right. Thus the successive
addition of first $u$ and then $v$ corresponds to the product
$uv$.
\end{remark}

\begin{remark}\label{quasidoubling}
One can define the group $\EE(R,P)$ using sequences of polytopes
$\fP'=(P=P'_0\subset P_1'\subset \cdots)$ that are more general
than doubling spectra. In particular, suppose that
$\fP=(P=P_0\subset P_1\subset\cdots)$ is a doubling spectrum for
$P$ and $\fP'=(P_0'\subset P_1'\subset \cdots)$ is a sequence of
polytopes for which there exist isomorphisms $\phi_i:P_i\to P_i'$
of polytopes that commute with the embeddings $P_i\subset P_{i+1}$
and $P_i'\to P_{i+1}'$. Then $\fP'$ need not be a doubling
spectrum in the strict sense since condition (ii) is not invariant
under isomorphisms as just described. However, there evidently
exists a natural isomorphism $\EE(R,\fP)\approx \EE(R,\fP')$.

For instance, if we start from the unimodular simplex $\Delta_n$,
$n\in\NN$, and consider the sequence $\fP'=(\Delta_n=P'_0\subset
P_1'\subset\cdots)$, in which $P'_{i+1}={P'_i}^{\sq_v}$,
$i\in\NN$, for the same column vector $ v\in\Col(\Delta_n)$, then
the resulting sequence of unstable groups is naturally identified
with the familiar sequence of groups of elementary matrices
$$
E_{n+1}(R)\subset E_{n+2}(R)\subset\cdots, \qquad*\mapsto
\begin{pmatrix}
*&0\\
0&1
\end{pmatrix}.
$$
In particular, $\EE(R,\Delta_n)=\E(R)$ for all $n\in\NN$.
\end{remark}

One may ask why the stable group of elementary automorphisms is
not defined as a filtered union of the unstable ones. The reason,
as explained in the Remark \ref{nfiltun}(b) below, is that such a
filtered union is in general \emph{not} possible.

\begin{remark}\label{nfiltun}
(a) In general the groups $\EE_R(P)$ are not perfect. However,
after finitely many steps in the doubling spectrum one arrives at
a polytope $P''$ for which this group is perfect.

In fact, after finitely many steps each base facet in $P$ has been
used for a doubling, and in the polytope $P'$ then constructed
each base facet has an invertible column vector. This property is
preserved under further doublings. In this situation it is easily
observed that all the vectors $\delta^+$ and $\delta^-$ that come
up in doublings of $P'$ are automatically decomposed, and after
finitely many doublings starting from $P'$ one arrives at a
polytope $P''$ in which all the column vectors of $P'$, and thus
all of $P''$, are decomposed. In view of Theorem \ref{elaut}(e)
this is enough for the desired perfectness. (The assertion holds
for all polytopes.)

(b) The group $\EE(R,P)$ is in general not the filtered union of
the unstable groups $\EE_R(P_i)$. Consider the simple example of
the segment $2\Delta_1$. Then the second term in $\fP$ (for an
arbitrary doubling spectrum, starting with $P$) can be identified
with the triangle
$2\Delta_2=\conv\bigl((0,0),(2,0),(0,2)\bigr)\subset\RR^2$ so that
$2\Delta_1$ is the lower edge (see Figure \ref{NonEmb}).
\begin{figure}[htb]
\begin{center}
\psset{unit=1cm}
\def\vertex{\pscircle[fillstyle=solid,fillcolor=black]{0.07}}
\begin{pspicture}(1.0,0)(1,2.2)
 \footnotesize
 \pspolygon[style=fyp](0,0)(2,0)(0,2)
 \multirput(0,0)(0,1){3}{\vertex}
 \multirput(1,0)(0,1){2}{\vertex}
 \multirput(2,0)(0,1){1}{\vertex}
 \psline[linewidth=1pt]{->}(1,0)(0,0)
 \psline[linewidth=1pt]{->}(1,0)(2,0)
 \rput(0.5,-0.3){$-v$}
 \rput(1.5,-0.3){$v$}
\end{pspicture}
\caption{The polytope $2\Delta_2$} \label{NonEmb}
\end{center}
\end{figure}
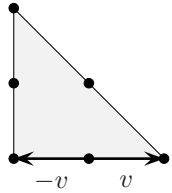
Consider the vectors  $v=(1,0)$ and $-v=(-1,0)$ from
$\Col(2\Delta_1)$. Assume $2\neq 0$ in $R$. Then the element
$\epsilon=(e_v^1\circ e^{-1}_{-v}\circ e_v^1)^2\in\EE(R,\Delta_2)$
is not the identity automorphism of $R[\fP]$ (it switches signs on
the second horizontal layer of $2\Delta_2$) whereas the
restriction of $\epsilon$ to $R[2\Delta_1]$ is the identity
automorphism. In particular there is no natural group homomorphism
$\EE_R(2\Delta_1)\to\EE(R,2\Delta_2)$.

(c) On the other hand, every element $e\in\EE_R(P_i)$, $i\in\ZZ_+$
is a restriction to $R[P_i]$ of some element of $\EE(R,P)$ and
every element $\epsilon\in\EE(R,P)$ restricts to an element of
$E_R(P_i)$ whenever $i$ is big enough. Clearly, we have the
following approximation principle: if two elements
$\epsilon,\epsilon'\in\EE(R,P)$ restrict to the same elements of
$\EE_R(P_i)$ for all sufficiently large $i$ then
$\epsilon=\epsilon'$.

(d) Unlike the group $\EE(R,P)$, the Steinberg group $\St(R,P)$,
to be introduced in Subsection \ref{SCHUR}, is the direct limit of
the corresponding unstable groups.
\end{remark}

\subsection{The Schur multiplier}\label{SCHUR}

Let $P$ be a balanced polytope and $\fP=(P\subset P_1 \subset
P_2\subset\dots)$ be a doubling spectrum. Then for a ring $R$ we
define the \emph{stable polytopal Steinberg group} $\St(R,P)$ as
the group generated by symbols $x_v^\lambda$, $v\in\Col({\mathfrak
P})$, $\lambda\in R$, which are subject to the relations
$$
x_v^\lambda x_v^\mu=x_v^{\lambda+\mu}
$$
and
$$
[x_u^\lambda,x_v^\mu]=
\begin{cases}
x_{uv}^{-\lambda\mu}&\text{if}\ uv\ \text{exists,}\\
\\
1&\text{if}\ u+v\notin\Col(\fP)\cup\{0\}.
\end{cases}
$$
The use of the notation $\St(R,P)$ is justified by the fact that,
like in Theorem \ref{elaut}(a), the stable Steinberg groups are
determined by the underlying doubling spectra (with the same
initial polytope) up to canonical isomorphism. Clearly, $\St(R,P)$
is a perfect group.

It follows from Remark \ref{steinberg} that for every $n\in\NN$ we
have $\St(R,\Delta_n)=\st(R)$ -- the usual Steinberg group.

\begin{lemma}\label{stemb}
Let $v_1,\dots,v_s$ be pairwise different elements of $\Col(\fP)$
with the same base facet. Then the mapping
$$
R^s\to\St(R,P),\qquad (\mu_1,\cdots,\mu_s)\mapsto
x_{v_1}^{\mu_1}\cdots x_{v_s}^{\mu_s}
$$
is a group isomorphism.
\end{lemma}

This follows from the Claim in the proof of \cite[Proposition
8.2]{BrG5}.

\begin{remark}
One can introduce (as we do in \cite{BrG5}) the notion of
`unstable' Steinberg polytopal group. The sequence of such groups,
associated to the members of the spectrum $\fP$, then forms an
inductive system of groups whose limit is $\St(R,P)$.
\end{remark}

The central result of \cite{BrG5} is the following

\begin{theorem}\label{schur}
For a ring $R$ and a balanced polytope $P$ the natural surjective
group homomorphism $\St(R,P)\to\EE(R,P)$ is a universal central
extension whose kernel coincides with the center of $\St(R,P)$.
\end{theorem}

See Proposition 8.2 and Theorem 8.4 in \cite{BrG5}.

The group $\Ker\bigl(\St(R,P)\to\EE(R,P)\bigr)$ is called the
\emph{polyhedral Milnor group}. We denote it by $K_2(R,P)$.

In \cite{BrG5} we have developed the notion of polytopal Steinberg
group for arbitrary lattice polytopes. However, one should note
that the proof of Theorem \ref{schur}, as presented in
\cite{BrG5}, uses the fact that $P$ is balanced in a crucial way.

\section{Rigid systems of column vectors}\label{RIGID}

The triangular subgroups of $\GL(n,R)$ on which Volodin's
$K$-theory is based are defined in terms of partial orders on
$\{1,\dots,n\}$ (see Subsection \ref{UNIMODULAR}). In this section
we develop a polyhedral analogue in terms of column vectors,
called rigid and $\Y$-rigid systems. We start with basic
properties of long products of column vectors.

It follows from Proposition \ref{maincrit}(e) that we can speak of
the product $\prod_{i=1}^mv_i$ of elements $v_i\in\Col(P)$
whenever the following two conditions are satisfied:
\begin{itemize}
\item[(i)] the products $v_iv_{i+1}$ exist for all $i\in[1,m-1]$,

\item[(ii)] $\sum_{i=r}^sv_i\neq 0$ for all $1\leq r<s\leq m$.
\end{itemize}
In this case every bracket structure on the sequence $v_1v_2\dots
v_m$ yields pairs of column vectors whose product exist.

Proposition \ref{maincrit}(c) implies that $v_j\|P_{v_i}$ whenever
$i<j$, and so $\rank_{\QQ}(v_1,\dots,v_m)=m$:

\begin{lemma}\label{diff}
If $\prod_{i=1}^mv_i$ exists, then $v_1,\dots,v_m$ are linearly
independent. In particular, $\sum_Iv_i\neq 0$ for all subsets
$I\subset[1,m]$, and $v_1,\dots,v_m$ are pairwise different column
vectors.
\end{lemma}

Obviously, if $\prod_{i=1}^mv_i$ exists, then $\prod_{i=r}^sv_i$
exists as well for all $r,s$, $1\leq r< s\leq m$.

For a system of column vectors $V\subset\Col(P)$ we set
$$
[V]=\{v\in \Col(P): \text{there exist } v_1,\dots,v_m \text{ with
} v=v_1\cdots v_m\}.
$$

It is useful to have another, weaker notion of product. We say
that $\prod_{i=1}^mv_i$ exists \emph{weakly} if there is a bracket
structure on the sequence
$$
v_1v_2\cdots v_m
$$
such that all the recursively defined products of pairs of column
vectors exist. Since $v_1\cdots v_n=v_1+\dots + v_n$ in the case
of weak existence, the value of the product does not depend on the
bracket structure.

It follows from Proposition \ref{maincrit}(f) that
$\prod_{i=1}^mv_i$ exists if and only if
$$
(v_1(v_2(v_3\dots(v_{m-3}(v_{m-2}(v_{m-1}v_m)))\dots))).
$$
exists and $\sum_{i=r}^sv_i\not=0$ for all $r,s$, $1\leq r<s\leq
m$.

By $\la V\ra$ we denote the hull of $V$ in $\Col(P)$ under
products (of two column vectors). One has $v\in\la V\ra$ if and
only if there exist $v_1,\dots,v_m\in V$ such that $v=v_1\cdots
v_m$ is their weak product.

Clearly $\la\la V\ra\ra=\la V\ra$, but in general $[[V]]\neq [V]$.
In fact, $[[V]]=[V]$ if and only $[V]=\la V\ra$. (A simple example
for $[[V]]\neq [V]$ will be discussed in Remark \ref{rigsyst}(b).)

Both $\la V\ra$ and $[V]$ carry an \emph{associative} partial
product structure. However, the partial product structure on $[V]$
is not always the restriction of that on $\Col(P)$. For
$w_1,w_2\in[V]$ the product may exist in $\Col(P)$, but it need
not belong to $[V]$ if $[V]\neq \la V\ra$.

For simplicity we introduce the following convention: $v_1\cdots
v_m\in [V]$ means that the product of $v_1,\dots, v_m$ exists (in
the strong sense), whereas $v_1\cdots v_m\in \la V\ra$ means that
the product of $v_1,\dots, v_m$ exists in the weak sense.

We will represent certain partial product structures on sets of
column vectors by equivalence classes of directed paths in graphs.
The \emph{graphs} considered by us will always be finite directed
graphs $\bG$ satisfying the following conditions:
\begin{itemize}
\item[(i)] $\bG$ has no isolated vertices;
\item[(ii)] $\bG$ has no multiple edges and no edges from a vertex
to itself;
\item[(iii)] if vertices $a$ and $b$ are connected by an edge, then
there is no other directed path connecting $a$ and $b$.
\end{itemize}
Condition (iii) implies that there are no directed cycles in $\bG$
(but the existence of non-directed cycles is not excluded). A
\emph{path} is always assumed to be oriented.

By definition, a \Y-\emph{graph} is a graph $\bF$ that at each
vertex $a$ satisfies the following condition:
$$
(\Y)\qquad \text{$a$ is the end point of at most one edge of
$\bF$.}
$$
In other words, if we direct all edges upwards, then branching is
only allowed in the form of a \Y\ (with any numbers of `arms').

The set of nonempty paths in a graph $\bF$ carries a natural
partial product structure -- $ll'$ exists if the end point of the
path $l$ is the initial point for $l'$. The set of all paths in
$\bF$ is denoted by $\path({\bF})$. There is an equivalence
relation on $\path{\bf F}$: two paths are considered to be
equivalent if they have the same initial and the same end point.
We let $\overline{\path{\bF}}$ denote the corresponding quotient
set. Thus for \Y-graphs (or more generally, for graphs without
non-oriented cycles) we have $\overline{\path{\bF}}=\path{\bF}$.
The aforementioned partial product operation on $\path{\bF}$
induces a partial product operation on $\overline{\path{\bF}}$. We
write $\overline{\path \bF}=\path\bF$ in order to indicate that
every equivalence class contains exactly one path.

In the following a vertex $a$ of $\bF$ is called \emph{terminal}
if there is no edge with initial vertex $a$.

\begin{definition}\label{rigid}
A system of column vectors $V\subset\Col(P)$ is called
\emph{rigid} if the following conditions are satisfied:
\begin{itemize}
\item[(a)] $[V]$ does not contain a subset of type $\{v,-v\}$,
$v\in\Col(P)$;

\item[(b)] $[V]=\la V\ra$;

\item[(c)] there exist a graph $\bF$ and an isomorphism
$[V]\approx\overline{\path{\bF}}$ of partial product structures.
\end{itemize}
Furthermore, $V$ is called a \Y-\emph{rigid system} if $\bF$ is a
\Y-graph.
\end{definition}

The graph $\bF$ and the isomorphism $[V]\approx\overline{\path{\bF}}$ in
Definition \ref{rigid}(c) are \emph{not} part of the data defining a rigid
system. We only require their existence. In general, $V$ does not uniquely
determine the graph $\bF$ (see Remark \ref{rigsyst}(a)). A graph
satisfying condition (c) will be referred to as a graph that
\emph{supports} the rigid system $V$, or a graph \emph{associated} to $V$.
Moreover, whenever a graph $\bF$ is associated to a rigid system $V$ it is
implicitly assumed that we have also fixed an isomorphism
$[V]\approx\overline{\path{\bF}}$ as above.

\begin{remark}\label{rigsyst}
(a) It is easy to find examples of rigid systems for which the
associated graph $\bF$ is not uniquely determined by $V$, not even
if $V$ happens to be \Y-rigid. For instance, when $V$ consists of
two vectors $u$ and $v$ such that $u+v\neq0$ and
$u+v\notin\Col(P)$, then $V$ is a \Y-rigid system, and the graphs
$\bG$, $\bF$ and $\bF'$ in Figure \ref{VarGraph1} support it. Two
of them, $\bF$ and $\bF'$, are non-isomorphic \Y-graphs.
\begin{figure}[htb]
$$
\psset{unit=1cm}
\def\vertex{\pscircle[fillstyle=solid,fillcolor=black]{0.05}}
\begin{pspicture}(-0.3,-0.3)(2.3,1.5)
 \rput(0,0){\vertex}
 \rput(1,1){\vertex}
 \rput(2,0){\vertex}
 \psline{->}(0,0)(1,1)
 \psline{->}(2,0)(1,1)
 \rput(0,0.5){\footnotesize$\bG$}
\end{pspicture}
\quad
\begin{pspicture}(-0.7,-0.3)(2.3,1.5)
 \rput(0,0){\vertex}
 \rput(0,1){\vertex}
 \rput(1,0){\vertex}
 \rput(1,1){\vertex}
 \psline{->}(0,0)(0,1)
 \psline{->}(1,0)(1,1)
 \rput(-0.3,0.5){\footnotesize$\bF$}
\end{pspicture}
\quad
\begin{pspicture}(-0.3,-0.3)(2.3,1.5)
 \rput(1,0){\vertex}
 \rput(0,1){\vertex}
 \rput(2,1){\vertex}
 \psline{->}(1,0)(0,1)
 \psline{->}(1,0)(2,1)
 \rput(0.1,0.5){\footnotesize$\bF'$}
\end{pspicture}
$$
\caption{}\label{VarGraph1}
\end{figure}
On the other hand, if $V$ is \Y-rigid, then the corresponding
\Y-graph is unique if we additionally require that there is only
one edge leaving each of its roots (vertices without an entering
edge). However, we will not make such a requirement.

(b) Even if a set of vectors can be arranged geometrically as a
graph, this does not imply the  rigidity of the system. Consider
the balanced polytope
$$
P=\conv\bigl((0,0,0),(1,0,0),(0,1,0),(1,1,0),(0,0,1)\bigr)\subset\RR^3
$$
and the system $V$ of its column vectors $u=(0,0,-1)$, $v=(1,0,0)$
and $w=(0,1,0)$; see Figure \ref{ThreeCol}.
\begin{figure}[htb]
\begin{center}
\psset{unit=2cm}
\def\vertex{\pscircle[fillstyle=solid,fillcolor=black]{0.05}}
\begin{pspicture}(-0.5,-0.7)(1,0.7)
\psset{viewpoint=-1 -2 -1.5} \ThreeDput[normal=0 0 1](0,0,0){
  \pspolygon[linewidth=0pt, style=fyp, fillcolor=medium](1,0)(0,0)(0,1)(1,1)
} \ThreeDput[normal=0 -1 0]{
  \pspolygon[linewidth=0pt, style=fyp](1,0)(0,0)(0,1)
} \ThreeDput[normal=1 0 0]{
  \pspolygon[linewidth=0pt, style=fyp](1,0)(0,0)(0,1)
} \ThreeDput[normal=0 0 1](0,0,0){
  \rput(0,0){\vertex}
  \rput(0,1){\vertex}
  \rput(1,0){\vertex}
  \rput(1,1){\vertex}
  \psline[linestyle=dashed](0,1)(0,0)
  \psline[linewidth=1.5pt]{->}(0,0)(1,0)
  \psline[linewidth=1.5pt]{->}(0,0)(0,1)
  \psline(0,1)(1,1)
} \ThreeDput[normal=0 -1 0]{
  \psline(0,1)(1,0)
  \psline[linewidth=1.5pt]{->}(0,1)(0,0)
}
 \ThreeDput[normal=1 0 0]{\psline(0,1)(1,0)}
 \ThreeDput[normal=1 -1 0]{\psline(0,1)(1.41,0)}
 \ThreeDput(0,0,1){\rput(0,0){\vertex}}
 \rput(-0.12,0.1){\footnotesize$u$}
 \rput(0.6,-0.05){\footnotesize$w$}
 \rput(-0.5,-0.3){\footnotesize$v$}
\end{pspicture}
\end{center}
\caption{The pyramid over the unit square} \label{ThreeCol}
\end{figure}
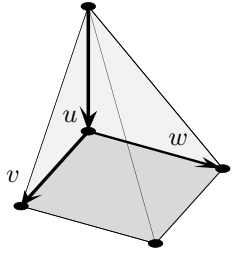
Let ${\bF}$ be the \Y-graph determined by the three vectors $u$,
$v$ and $w$ as shown in the figure. Then $[V]\approx\path{\bF}$
but $(uv)w\in\la V\ra\setminus[V]$.
\end{remark}

Let $V\subset \Col(P)$. The elements $w$ of $V$ that have no
decomposition $w=uv$ with $u,v\in \la V\ra$ are called the
\emph{irreducible} or \emph{indecomposable} elements of $V$.

In the next proposition we collect some properties of rigid
systems, among them the decomposition into irreducible elements.

\begin{proposition}\label{basrig}
Let $V\subset \Col(P)$ be a rigid system and $\bF$ a graph
associated to $V$. Then the following hold:
\begin{itemize}
\item[(a)] $[[V]]=[V]$.
\item[(b)] The product $v_1\cdots v_n$ of elements $v_i\in [V]$
exists if (and only if) it exists weakly.
\item[(c)] $V$ is \Y-rigid if and only if $\bF$ satisfies the
condition $(\Y)$ at all its nonterminal vertices $c$, and $\path
\bF=\overline{\path \bF}$.
\item[(d)] Every $w\in [V]$ has a decomposition into irreducible
elements. The irreducible elements are those represented by edges
of $\bF$. The decomposition is unique if $V$ is \Y-rigid.
\item[(e)] Let $v_1,\dots,v_n$, $n\ge 1$, be arbitrary elements of $[V]$.
Then $\sum_{i=1}^n v_i\neq 0$.
\end{itemize}
\end{proposition}

\begin{proof}
(a) follows immediately from $[V]=\la V\ra$.

(b) For the partial product structures $\path \bF$ and
$\overline{\path \bF}$ we can define strongly and weakly existing
products as for column vectors. In contrast to $\Col(P)$, strong
and weak existence are evidently equivalent in these partial
product structures. Therefore they must be equivalent in a partial
product structure isomorphic to $\overline{\path\bF}$. (Note that (b) does
not necessarily hold for an arbitrary subset $V$ of $\Col(P)$ that just
satisfies the equation $\la V\ra=[V]$: for example, it is in general false
for $V=\Col(P)$ itself.)

(c) Suppose that $\bF$ satisfies (\Y) at all nonterminal vertices.
Two equivalent paths with the same initial points and endpoints
would necessarily end in a terminal vertex $c$ of $\bF$.

If such paths do not exist, we can replace the edges ending in $c$
by edges with separate endpoints without changing the path
structure. The resulting graph is a \Y-graph.

Conversely, if such paths exist, then it is impossible to find a
\Y-graph $\bF'$ with $\overline{\path \bF'}\approx \overline{\path
\bF}$.

(d) The existence of the decomposition is clear since every
element in $[V]$ has a longest representation as a product
$v_1\cdots v_n$. Its factors must be irreducible. The irreducible
elements correspond to the edges of $\bF$ since (in our class of
graphs) no edge is equivalent to a path of length $\ge2$. The
uniqueness of the decomposition in \Y-rigid systems follows from
the fact that a path in a \Y-graph is uniquely determined by its
initial and endpoint.

(e) Suppose that $\sum_{i=1}^n v_i=0$. We can assume that
$v_1\cdots v_m$ is a longest possible product that can be formed
from the vectors $v_1,\dots,v_n$. Then $w=v_1\cdots v_m$ is a
column vector. Let $F$ be its base facet. Among the remaining
vectors $v_{m+1},\dots,v_n$ there must be one, say $v_{m+1}$, with
$\la F,v_{m+1}\ra>0$: note that $\sum_{i=1}^n\la F, v_i\ra=0$, but
$\sum_{i=1}^m\la F, v_i\ra=\la F,v_1\cdots v_m\ra <0$.

It is impossible that $v_{m+1}=-w$. So $v_{m+1}(v_1\cdots v_m)$ exists
weakly. By (b) it exists also in the strong sense, and we obtain a
contradiction to the choice of $v_1,\dots,v_m$. \end{proof}

The basic examples of rigid systems are provided by the unit
simplices and their edge vectors:

\begin{example}\label{basic}
The set of column vectors $\Col(\Delta_n)$ of the unit $n$-simplex
$\Delta_n$ coincides with the set $\{a-b: a \text{ and } b \text{
are different vertices of } \Delta_n\}$. We can think of these
column vectors as oriented edges of $\Delta_n$ from $b$ to $a$.

(a) Then $uv$ exists for column vectors $u$ and $v$ if and only if
they form a broken line of length 2.  More generally, for every
system $\{v_1,\dots,v_m\}\subset\Col(\Delta_n)$ the following
conditions are equivalent:
\begin{itemize}
\item[(i)]
$\prod_{i=1}^mv_i$ exists,
\item[(ii)]
$\{v_1,\dots,v_m\}$ is a \Y-rigid system and its underlying
\Y-graph is the linear directed graph of length $m$:
$$
\psset{unit=1cm}
\def\vertex{\pscircle[fillstyle=solid,fillcolor=black]{0.05}}
\begin{pspicture}(-0.3,-0.3)(4.3,0.3)
 \multirput(0,0)(1,0)5{\vertex}
 \psline{->}(0,0)(1,0)
 \psline{->}(1,0)(2,0)
 \psline{->}(3,0)(4,0)
 \multirput(2.2,0)(0.2,0){4}{.}
 \rput(0.5,0.2){\footnotesize 1}
 \rput(1.5,0.2){\footnotesize 2}
 \rput(3.5,0.2){\footnotesize $m$}
\end{pspicture}
$$
\item[(iii)]
$[v_1\,\dots,v_m]\subset\partial\Delta_n$ is a broken line without
self-intersections ($\partial$ denotes the boundary).
\end{itemize}
Observe that the weak existence of $\prod_{i=1}^mv_i$ together
with $\sum_Iv_i\neq 0$ for every subset $I\subset[1,m]$ is
equivalent to the existence of $\prod_{i=1}^mv_i$ in the strong
sense.

(b) The rigid systems in $\Col(\Delta_n)$ are exactly the
non-empty subsets $V$ for which $\la V\ra$ does not contain a pair
$\{v,-v\}$. By induction on $m$ it follows easily from (a) that
all products $v_1\cdots v_m\in\la V\ra$ exist also in the strong
sense, and if we take the longest possible representation
$w=v_1\cdots v_m$ of $w\in\la V\ra$ with $v_i\in V$, $i\in[1,m]$,
then the $v_i$ must be irreducible. The graph $\bF$ formed by the
irreducible elements in the boundary of $\Delta_n$ then satisfies
the condition $\overline{\path F}\approx [V]$. (If $V$ is
\Y-rigid, then the graph just produced need not be \Y-rigid; see
Remark \ref{rigsyst}(a).)

(c) In preparation of Subsection \ref{UNIMODULAR} we note that the
rigid systems $V$ in $\Col(\Delta_n)$ can be identified with those
partial orders on subsets $X$ of the vertex set $\vert(\Delta_n)$
for which no $x\in X$ is simultaneously maximal and minimal. The
set $X$ corresponding to $V$ consists of all vertices $x$ of
$\Delta_n$ such that there exists $y\in\vert(\Delta_n)$  for which
$x-y\in[V]$ or $y-x\in [V]$, and it is partially ordered by the
condition $x\le y \iff y-x\in[V]$.

Among these partially ordered sets the \Y-rigid systems in
$\Col(\Delta_n)$ are characterized by the following two
conditions:
$$
a<c \text{ and }b<c \text{ for some }c\notin\max(X) \implies a\leq
b \text{ or } b\leq a.
$$
and
$$
a<b,c \text{ and } b,c<d\in\max(X)\implies b\le c \text{ or } c\le
b.
$$
These conditions reflect Proposition \ref{basrig}(c).

(d) Since any finite poset can be augmented to a linear order, we
conclude that for any rigid system $V\subset\Col(\Delta_n)$ there
is a `linear' \Y-rigid system $W\subset\Col(\Delta_n)$ as in (ii)
above such that $[V]\subset[W]$. Namely, we augment the partial
order on $X=\{x_1,\dots,x_m\}$ (as defined in (b)) to a linear
one. we may assume the vertices are labelled such that $x_i<x_j$
if $i<j$, and choose $W=\{x_{i+1}-x_i: i=1,\dots,m-1\}$. See
Figure \ref{VarGraph} where $\bR$ is the graph of a rigid system
$V$ and $\bS$ is the graph of a `linear' \Y-rigid system $W$
containing it.
\begin{figure}[htb]
$$
\psset{unit=0.7cm}
\def\vertex{\pscircle[fillstyle=solid,fillcolor=black]{0.05}}
\begin{pspicture}(-1,0)(1,2)
 \rput(0,0){\vertex}
 \rput(-1,1){\vertex}
 \rput(1,1){\vertex}
 \rput(0,2){\vertex}
 \psline{->}(0,0)(-1,1)
 \psline{->}(0,0)(1,1)
 \psline{->}(-1,1)(0,2)
 \psline{->}(1,1)(0,2)
 \rput(0.0,0.7){\footnotesize$\bR$}
\end{pspicture}
\qquad
\begin{pspicture}(-1,0)(1,2)
 \rput(0,0){\vertex}
 \rput(-1,1){\vertex}
 \rput(1,1){\vertex}
 \rput(0,2){\vertex}
 \psline{->}(0,0)(-1,1)
 \psline{->}(-1,1)(1,1)
 \psline{->}(1,1)(0,2)
 \rput(0.0,0.7){\footnotesize$\bS$}
\end{pspicture}
\qquad
\begin{pspicture}(-1,0)(1,2)
 \rput(-1,0){\vertex}
 \rput(1,0){\vertex}
 \rput(0,1){\vertex}
 \rput(0,2){\vertex}
 \psline{->}(-1,0)(0,1)
 \psline{->}(1,0)(0,1)
 \psline{->}(0,1)(0,2)
 \rput(0.0,0.2){\footnotesize$\bT$}
 \rput(-0.8,0.6){\footnotesize$u$}
 \rput(0.85,0.7){\footnotesize$u'$}
 \rput(0.3,1.5){\footnotesize$v$}
\end{pspicture}
$$
\caption{}\label{VarGraph}
\end{figure}
\end{example}

\begin{remark}\label{whybal}
In general it is not possible to embed a rigid system in $\Col(P)$
into a \Y-rigid system. The pyramid $P$ over the unit square (see
Figure \ref{ThreeCol}) serves again as an example. The column
vectors $u,u'=u+w,v$ form a rigid system $V$ with associated graph
$\bT$ as shown in Figure \ref{VarGraph}. A \Y-rigid system
$V'\supset V$ would have to contain a column vector $t$ with
$u'=tu$ or $u=tu'$, in other words, $t=w=u'-u\in V'$ or $t=-w\in
V'$. However, $u'=uw$ and $u=u'(-w)$ -- the products exist in the
wrong order.
\end{remark}

As far as the partial product structure is concerned, every rigid
system can be realized in a unit simplex:

\begin{lemma}\label{colemb}
Let $P$ be a polytope and $V\subset \Col(P)$ a rigid system with
associated graph $\bF$. Set $n=\#\vert(\bF)-1$.
\begin{itemize}
\item[(a)] Then $[V]$ is isomorphic to a rigid system $W$ in
$\Col(\Delta_n)$ as a partial product structure.
\item[(b)] Every subset $U\subset[V]$ is a rigid system.
\end{itemize}
\end{lemma}

\begin{proof}
(a) We identify the vertices of $\bF$ with those of $\Delta_n$,
and let $W'$ be the set of edge vectors of $\Delta_n$
corresponding to the edges of $\bF$. Then $W=\la W'\ra$ is the set
of edge vectors in $\Delta_n$ that as a path are equivalent to
some broken line formed by the elements of $W'$. None such broken
line intersects itself: $\bF$ has no directed cycles. Moreover it
is impossible that $\{w,-w\}\subset W$ for some
$w\in\Col(\delta_n)$ -- again we would have a directed cycle in
$\bF$. Clearly $W$ is a rigid system with associated graph $\bF$.

(b) Since $\{v,-v\}\not\subset [U]$ for any $v\in\Col(P)$, the
rest is only a condition on the partial product structure of
$[U]$. By (a) we can therefore assume $V\subset \Delta_n$, and
then the claim follows from the observation in Example
\ref{basic}(b).
\end{proof}

Note that part (b) of Lemma \ref{colemb} has no \Y-version: in
general a subset of a \Y-rigid system need not be \Y-rigid, as is
clear from Example \ref{basic}(d).

The construction of stable groups uses the doubling of polytopes.
Therefore we must analyze how rigid systems can be extended to the
doubled polytope. The result will be described in Lemma
\ref{decomp}. First we formulate an auxiliary lemma that sheds
more light on the structure of rigid (and somewhat weaker) systems
$V$. It shows that the table of values $\la F,v\ra$ for such a
system looks like that of a system of column vectors in a unit
simplex, at least if one pays attention only to those facets $F$
that appear as base facets of elements of $V$.

\begin{lemma}\label{linear}
Let $P$ be a polytope, $V\subset \Col(P)$, and $\cB(V)$ be the set
of the base facets of the vectors $v\in V$.
\begin{itemize}
\item[(a)] Then the following are equivalent:
\begin{itemize}
\item[(i)]
Whenever $v_1\cdots v_n\in\la V\ra$ with $n\in\NN$ and
$v_1,\dots,v_n\in V$, then $v_iv_{i+1}$ exists or $v_{i+1}=-v_i$
for each $i\in[1,n-1]$.
\item[(ii)] For each $v\in V$ one has $\la F,v\ra=0$, $F\neq P_v$, for all
$F\in\cB(V)$ with at most one exception $G$ for which $\la
G,v\ra=1$.
\end{itemize}
\item[(b)] Suppose one of the conditions in \textup{(a)} holds.
For every product $w=v_1\cdots v_n\in\la V\ra$ with
$v_1,\dots,v_n\in V$ there exists at most one facet $G\in\cB(V)$
with $\la G,w\ra>0$. Moreover, in the case of existence one has
$\la G,w\ra=\la G, v_n\ra=1$.
\end{itemize}
\end{lemma}

The face $G\in\cB(V)$ with $\la G,w\ra=1$ (if it exists) is in
some sense opposite to the base facet of $w$. Therefore we denote
it by $P_V^w$.

\begin{proof}
In this proof and that of Lemma \ref{decomp} we use Proposition
\ref{maincrit} and Lemma \ref{diff} freely, without mentioning
each application explicitly.

Note that the property in (ii) is automatically (and `globally')
satisfied for all column vectors $v$ with $-v\in\Col(P)$; see
Proposition \ref{maincrit}(g).

We first prove (i)$\implies$(ii). More generally than necessary,
let $v\in \la V\ra$, $v=v_1\cdots v_n$ with $v_1,\dots,v_n\in V$,
and suppose first that $\la F,v\ra\ge 2$ for some $F\in\cB(V)$.
Then $-v\notin\Col(P)$, and so $vw=(v_1\cdots v_n)w$ exists for
every $w\in V$ with $F=P_w$. Moreover, $(vw)w$ also exists (since
$vw=-w$ is evidently impossible), but $ww$ does not. This shows
$\la F,v\ra\le 1$ for all $F\in\cB(V)$.

Next suppose that there exist different $F_1,F_2\in\cB(V)$ with
$\la F_1,v\ra=\la F_2,v\ra=1$ and choose $w_i\in V$, $i=1,2$, with
$F_i=P_{w_i}$. (Again $-v\notin \Col(P)$.) Both $vw_1$ and $vw_2$
exist. One has $\la F_2,vw_1\ra=\la F_2, v\ra+\la F_2,w_1\ra\ge 1$
since $\la F_2,v\ra=1$ and $\la F_2,w_1\ra\ge 0$. We have already
seen that $\la F_2,vw_1\ra\le 1$, and so $\la F_2,w_1\ra=0$. Thus
$w_1\neq -w_2$ and $w_1w_2$ does not exist, but $(vw_1)w_2$
exists: it is impossible that $vw_1=-w_2$, since $-v\in\Col(P)$
otherwise.

For (ii) $\implies$ (i) and (ii) $\implies$ (b) we use induction
on $n$. For $n=1$ there is nothing to show in (i), and (b) is
identical to (ii). Let $n>1$. Then
$$
w=v_1\cdots v_n=w_1w_2,\qquad w_1=v_1\cdots v_m,\
w_2=v_{m+1}\cdots v_n,
$$
with $1\le m<n$ and $w_1,w_2\in \la V\ra$. We can apply induction
to each of the shorter products $w_1$ and $w_2$.

Clearly $F=P_{w_2}=P_{v_{m+1}}\in\cB(V)$, and $\la F,w_1\ra>0$. By
the induction hypothesis for (b) we have $F=P_V^{w_1}=P_V^{v_m}$
and $\la F,w_1\ra=\la F,v_m\ra=1$. By an application of the
induction hypothesis to $w_2$ it follows immediately that among
all the facets $G\in\cB(V)$ there is at most one with $\la
G,w_2\ra >0$, namely $P_V^{v_n}$ if it exists. This completes (b).

If $v_{m+1}\neq -v_m$, then the product $v_mv_{m+1}$ exists since
$\la F, v_m\ra>0$, and (i) is also complete.
\end{proof}

Now we can describe how rigid systems can be extended under
doublings of balanced polytopes along facets.

\begin{lemma}\label{decomp}
Suppose $V\subset\Col(P)$ is a rigid system and $v\in V$ is one of
its irreducible elements. Then the system
$V'=V\cup\{\delta^+,v^\vt\} \subset\Col(P^{\sq_v})$ is also rigid.
It is \Y-rigid if $V$ is so.
\end{lemma}

\begin{proof}
Let ${\bF}$ be $V$'s underlying graph, and $E\in E({\bF})$ be the
edge corresponding to the column vector $v$. According to the
equation $v=\delta^+v^\vt$ we replace $E$ by a path
\begin{pspicture}(-0.25,-0.2)(2.2,0.4)
 \rput(0,0){\vertex}
 \rput(1,0){\vertex}
 \rput(2,0){\vertex}
 \psline{->}(0,0)(1,0)
 \psline{->}(1,0)(2,0)
\end{pspicture}
where the left edge represents $\delta^+$, and the right $v^\vt$.
Let $\bF'$ be the extended graph.

The rigid system $V$ has property (i) of Lemma \ref{linear}(a) (in
which the possibility $v_{i+1}=-v_i$ is excluded for $V$ rigid).
Let $W$ consist of all the irreducible elements $w\neq v$ of $V$
and the `new' column vectors $\delta^+$ and $v^\vt$. We have
$\cB(W)=\cB(V)\cup\{P^-\}$, identifying each facet of $P$ with its
extension to $P^{\sq_v}$ as given by the mapping $\Psi$ in
Subsection \ref{DOUBLING}. Since $V$ satisfies (ii) in Lemma
\ref{linear}(a), one sees immediately that $W$ also satisfies it.
In fact, each $u\in V$ is parallel to $P^-$, $\delta^+$ is an
invertible column vector, and if $\la F,v^\vt\ra>0$, then $\la F,
v\ra>0$.

Thus Lemma \ref{linear} allows us to control the pairs
$w_iw_{i+1}$ in weakly existing products $w_1\cdots w_n\in \la
W\ra=\la V'\ra$. By inspection of the base facets and the
`positive' facets with respect to $W$ one sees that only the
following types can occur for $w_i,w_{i+1}$:
\begin{align*}
&w_i,w_{i+1}\in V\text{ and } w_iw_{i+1} \text{ exists in $V$},\\
&w_i\in V,\ w_{i+1}=\delta^+ \text{ and } w_iv \text{ exists in $V$},\\
&w_i=v^\vt, w_{i+1}\in V \text{ and } vw_{i+1} \text{ exists in $V$},\\
&w_i=\delta^+,w_{i+1}=v^\vt.
\end{align*}
Therefore each product $w_1\cdots w_n\in\la W\ra$ represents a
path in $\bF'$.

Now we choose a product $w_1\cdots w_n$ `along' a path in $\bF'$
and show that it exists strongly. It is enough to consider a
maximal path since strong existence is inherited by segments. We
must verify the existence of
$$
w_i(w_{i+1}\cdots w_n),\qquad i=1,\dots,n-1
$$
and show that $\sum_{i=r}^s w_i\neq 0$ for all $r,s$, $1\le r<s\le
m$.

If $w_1,\cdots,w_n\in V$, the strong existence follows from that
in $V$. Otherwise $w_{i-1}=\delta^+$ and $w_i=v^\vt$ for exactly
one $i$, whereas $w_j\in V$ for $j\neq i-1,i$.

Let us first take care of the conditions on sums over segments. If
the segment contains none or both of $\delta^+$ and $v^\vt$, then
we are summing column vectors in $V$ along a path of $\bF$, and
such a sum is necessarily non-zero. If the segment contains
$\delta^+$, but not $v^\vt$, then $\la P^-,\sum_{i=r}^s w_i\ra=1$,
and if it contains $v^\vt$, but not $\delta^+$, then $\la
P^-,\sum_{i=r}^s w_i\ra=-1$. In any case $\sum_{i=r}^s w_i\neq 0$.
If we have to check the existence of $w_i(w_{i+1}\cdots w_n)$ in
the following, then we can use that $w_i\neq-w_{i+1}\cdots w_n$,
as just shown.

Clearly $w_{i+1}\cdots w_n\in[W]$ since it represents a path in
$\bF$ (unless it is empty).

If $i=n$, there is nothing to show for the existence of $w_i\cdots
w_n$. If $i<n$, then $vw_{i+1}$ exists, since $vw_{i+1}$ is a path
in $\bF$. The base facet of $w_{i+1}\cdots w_n$ is
$P_W^v=P_W^{v^\vt}$ and so the product $w_i(w_{i+1}\cdots w_n)$
exists.

Next we must attach $w_{i-1}=\delta^+$ at the left side of
$w_i\cdots w_n$. One has $\la P^-,\delta^+\ra=1$ for the base
facet $P^-$ of $w_i=v^\vt$ which is also the base facet of
$w_i\cdots w_n$. Again we are done.

After having shown the existence of $w_{i-1}(w_i\cdots w_n)$, we
can replace $w_{i-1}w_i=\delta^+v^\vt$ by $v$, and from now on the
product $vw_{i+1}\cdots w_n$ and the succeeding ones are prefixed
by elements from $V$, a harmless operation.

To sum up: we have shown that all the weakly existing products
represent paths in $\bF'$, and that each such path yields a
strongly existing product. That the equivalence classes of the
paths in $\bF'$ represent the elements of $[W]$ follows
immediately from the corresponding property of $\bF$ for $V$.

It remains to show that $[W]=[V']$ does not contain a column
vector $u$ and its inverse $-u$. The only critical pairs of
products are those in which one element contains $v^\vt$, but not
$\delta^+$, and the other contains $\delta^+$, but not $v^\vt$.
The first product must end in $\delta^+$ and the second must start
with $v^\vt$. But then we can concatenate them to a path in
$\bF'$, and the sum over such paths is nonzero, as shown already.
\end{proof}

The last lemma of this section will be used in Section \ref{ACYCL}
in the context of Mayer-Vietoris sequences.

\begin{lemma}\label{balint}
Assume $U,V\subset\Col(P)$ are rigid systems. Then the
intersection $[U]\cap[V]$ is also a rigid system. Moreover, if $U$
and $V$ are \Y-rigid, then $[U]\cap[V]$ is \Y-rigid, too.
\end{lemma}

\begin{proof}
Set $W=[U]\cap[V]$. It follows from Lemma \ref{colemb}(b) that $W$
is a rigid system.

Only the claim on \Y-rigidity has yet to be proved. Suppose $U$
and $V$ are \Y-rigid, and let ${\bF}_U$ and ${\bF}_V$ be the
corresponding \Y-graphs. Then we construct a graph $\bF$
associated with $W$ as in the proof of Lemma \ref{colemb}(b) (say
from $\bF_U$). Let $W'$ be the set of irreducible elements in $W$.

First we have to make sure that there is no triple of distinct
elements $w_1,w_2,w_3\in W'$ such that the products $w_1w_3$ and
$w_2w_3$ exist -- this will show that all nonterminal vertices $c$
of ${\bF}$ satisfy the condition (\Y). But if both $w_1w_3$ and
$w_2w_3$ existed, then the corresponding paths
$l(w_1),l(w_2)\in\path{\bF}_U$ would have the same terminal point
and either $l(w_1)\subset l(w_2)$ or $l(w_2)\subset l(w_1)$. We
may assume $l(w_1)\subset l(w_2)$; then $w_2=w'w_1$ for some
$w'\in[U]$. Since the vector $w'\in\Col(P)$ is uniquely
determined, the same arguments applied to $V$ show that
$w'\in[V]$, that is, $w_2$ is decomposable within $W$ -- a
contradiction.

By Proposition \ref{basrig}(c) the only remaining obstruction to
the \Y-rigidity of $W$ is the existence of a (non-directed) cycle
in $\bF$ that is the union of two oriented paths with the same
initial point and the same terminal point.

Assume such a cycle exists. Then there are elements
$w_1,\dots,w_r,w'_1,\dots,w'_s\in W'$, $r,s\in\NN$ such that
$$
w_1\cdots w_r=w'_1\cdots w'_s.
$$
We have to show that $r=s$ and $w_i=w_j$ for $i\in[1,r]$. Consider
the corresponding paths $l(w_i),l(w'_j)\in\path{\bF}_U$. Since
$\overline{\path{\bF}_U}={\path{\bF}_U}$ we get
$$
l(w_1)\cdots l(w_r)=l(w'_1)\cdots l(w'_s)\in\path{\bF}_U,
$$
where the multiplication is the concatenation of paths. We may
assume that $l(w_r)\neq l(w'_s)$. Since ${\bF}_U$ is a \Y-graph
either $l(w_r)\subset l(w'_s)$ or $l(w'_s)\subset l(w_r)$. Then
the same arguments as in checking the condition (\Y)  for
nonterminal vertices show that $w'_s$ or $w_r$ is decomposable in
$W$ -- a contradiction.
\end{proof}

We conclude this section by a further examination of products
$v_1\cdots v_n$. It has been observed in Example \ref{basic}(a)
that the existence of $v_1\cdots v_m$ in $\Col(\Delta_n)$ implies
the rigidity of $\{v_1,\cdots,v_m\}$. This is not true for all
polytopes:

\begin{example}\label{nonrig}
Let $P$ be the 3-simplex $\conv((0,0,0),(2,0,0),(0,2,0),(0,0,1))$
and consider its column vectors $u=(1,0,-1),v=(-1,0,0),w=(0,1,0)$
(see Figure \ref{ExNonRig}).
\begin{figure}[htb]
\begin{center}
\psset{unit=1.5cm}
\def\vertex{\pscircle[fillstyle=solid,fillcolor=black]{0.05}}
\begin{pspicture}(-0.2,-1.0)(1,0.9)
 \psset{viewpoint=-1 -2 -1.5}
\ThreeDput[normal=0 0 1](0,0,0){
  \pspolygon[linewidth=0pt, style=fyp, fillcolor=medium](2,0)(0,0)(0,2)
}
  \ThreeDput[normal=0 -1 0]{
  \pspolygon[linewidth=0pt, style=fyp](2,0)(0,0)(0,1)
}
  \ThreeDput[normal=1 0 0]{
  \pspolygon[linewidth=0pt, style=fyp](2,0)(0,0)(0,1)
}
  \ThreeDput[normal=0 0 1](0,0,0){
  \rput(0,0){\vertex}
  \rput(0,1){\vertex}
  \rput(0,2){\vertex}
  \rput(1,0){\vertex}
  \rput(1,1){\vertex}
  \rput(2,0){\vertex}
  \psline[linewidth=1.5pt]{->}(1,0)(0,0)
  \psline[linewidth=1.5pt]{->}(0,0)(0,1)
}
  \ThreeDput[normal=0 -1 0]{
  \rput(0,1){\vertex}
  \psline[linewidth=1.5pt]{->}(0,1)(1,0)
}
 \rput(0.64,0.25){\footnotesize$u$}
 \rput(0.45,0){\footnotesize$v$}
 \rput(-0.35,-0.2){\footnotesize$w$}
\end{pspicture}
\end{center}
\caption{} \label{ExNonRig}
\end{figure}
Then $uvw$ and $uw$ exist. This excludes the rigidity of
$\{u,v,w\}$.
\end{example}

This example and many other observations in \cite{BrG5} and in
this paper, have naturally lead us to the class of balanced
polytopes. Their combinatorial properties allow one to develop
polyhedral $K$-theory. The next proposition shows that the
phenomenon just observed is indeed impossible in balanced
polytopes.

\begin{proposition}\label{broklin}
Let $P$ be a balanced polytope. Assume  $\prod_{i=1}^mv_i$ exists
for $V=\{v_1,\dots,v_m\}\subset\Col(P)$. Then $V$ is a \Y-rigid
system whose associated \Y-graph is the directed linear graph
\begin{pspicture}(-0.25,-0.2)(3.2,0.4)
 \rput(0,0){\vertex}
 \rput(1,0){\vertex}
 \rput(2,0){\vertex}
 \rput(3,0){\vertex}
 \psline{->}(0,0)(1,0)
 \psline{->}(2,0)(3,0)
 \psline[linestyle=dotted, linewidth=0.6pt](1,0)(2,0)
\end{pspicture}
of length $m$.
\end{proposition}

\begin{proof}
In this proof we use Proposition \ref{maincrit} heavily. Set
$F_i=P_{v_i}$. We have to show that
$$
\bigl(\la F_j,v_i\ra\bigr)=
\begin{pmatrix}
 -1 & 1 & 0 &\cdots & 0\\
 0 & -1 & 1 &\ddots &\vdots\\
 \vdots & \ddots &\ddots & \ddots &0\\
 \vdots & & \ddots & -1 & 1\\
 0& \cdots &\cdots &0 &-1
\end{pmatrix}
$$
where row $i$ corresponds to $v_i$ and column $j$ to $F_j$. The
case $j=i$ is clear by definition. Clearly $\la F_j,v_i\ra=0$ for
all $i>j$ by the existence of $v_j\cdots v_n$ (independently of
the fact that $P$ is balanced). But also $w=v_1\cdots v_{j-1}$
exists, and since $P$ is balanced, we must have $\la
F_j,w\ra=\sum_{i=1}^{j-1}\la F_j,v_i\ra\le 1$. Since $\la
F_j,v_{j-1}\ra>0$ and $\la F_j,v_i\ra\ge 0$ for $i<j$, this
implies $\la F_j,v_{j-1}\ra=1$ and $\la F_j,v_i\ra=0$ for $i<j-1$.
\end{proof}

\begin{remark}\label{weakbl}
If the product $\prod_{i=1}^m v_i$ exists only in the weak sense
and $\sum_{i\in I}v_i\neq 0$ for every subset $I\subset[1,m]$,
then $\{v_1,\dots,v_n\}$ need not be a \Y-rigid system. Consider
the balanced polytope $P$ of Remark \ref{rigsyst}(b) and the same
column vectors $u=(0,0,-1)$, $v=(1,0,0)$ and $w=(0,1,0)$.
\end{remark}

\section{Triangular subgroups in $\EE(R,P)$ and $\St(R,P)$}\label{TRIANG}

In this section we generalize the notion of a triangular group of
matrices to the polyhedral setting. These groups play a crucial
r\^ole in the definition of Volodin simplicial sets (Section
\ref{HPK}).

Let $R$ be a ring and $P$ a balanced polytope admitting a column
structure. We fix a doubling spectrum $\fP=(P\subset P_1
\subset\cdots)$. Thanks to Theorem \ref{elaut}(a) (and its
straightforward analogue for polyhedral Steinberg groups) all the
objects defined below are independent of the fixed spectrum.

We say that $V\subset\Col(\fP)$ is a rigid (\Y-rigid) system if
there exists an index $j\in\NN$ such that $V\subset\Col(P_j)$ and
is rigid (\Y-rigid) in the sense of Definition
\ref{rigid}.\pagebreak[3]

\begin{definition}\label{triang}\leavevmode\par
\begin{itemize}
\item[(a)] A subgroup $G\subset\EE(R,P)$ is called
\emph{triangular} if there exists a rigid system
$V\subset\Col(\fP)$ such that $G$ is generated by the elementary
automorphisms $e_v^\lambda$, where $\lambda$ runs through $R$ and
$v$ through $V$. The triangular subgroup corresponding to a rigid
system $V$ is denoted by $G(R,V)$, and $\T(R,P)$ is the family of
all triangular subgroups of $\EE(R,P)$.

\item[(b)] The triangular subgroups of $\St(R,P)$ are defined
similarly, and $G'(R,V)$ and $\T'(R,P)$ denote the corresponding
objects.

\item[(c)] The \Y-triangular subgroups in $\EE(R,P)$ and
$\St(R,P)$ are those supported by \Y-rigid systems. Their families
are $\T(R,P)^\Y$ and $\T'(R,P)^\Y$.
\end{itemize}
\end{definition}

Let $\bF$ be a \Y-graph underlying a \Y-rigid system
$V\subset\Col(\fP)$ and let $E({\bF})$ be the set of edges of
$\bF$. There is a natural partial order on $E({\bF})$ defined as
follows: for $f,g\in E({\bF})$ we put $f\leq g$ if there is
$l\in\path{\bF}$ with first edge $f$ and last edge $g$. We have
the disjoint partition
$$
E({\bF})=E_1\cup E_2\cup\cdots\cup E_t
$$
where each of the $E_r$ consists of those elements $f$ of $E({\bf
F})$ that admit sequences of type $f_1<f_2<\cdots<f_r=f$ and do
not admit sequences $f_0<f_1<f_2<\cdots<f_r=f$. ($E_t$ is the set
of maximal elements.) Edges in $E_r$ have \emph{degree} $r$.

We get the partition
$$
\path{\bF}=\bigcup_{r,s}\path_{rs}{\bF}
$$
into disjoint sets where $1\leq r\leq s\leq t$ and
$$
\path_{rs}{\bF}=\{[f_r,\cdots,f_s]\ |\ f_r\in E_r,\cdots,f_s\in
E_s\}.
$$

Let $v_l\in[V]$ denote the column vector corresponding to a path
$l\in\path({\bF})$. Then we have the analogous disjoint partition
$$
[V]=\bigcup_{rs}[V]_{rs}
$$
where $1\leq r\leq s\leq t$ and
$$
[V]_{rs}=\{v_l\ |\ l\in\path_{rs}{\bF}\}.
$$
We introduce the following notation:

\begin{itemize}
\item
$[V]^r=\{v_{r1},\dots,v_{rN_r}\}=\bigcup_{s=r}^t[V]_{rs}$, for
$r\in[1,t]$ ($N_r=\#\bigcup_{s=r}^t[V]_{rs}$). That is, $[V]^r$
consists of the column vectors which correspond to paths with
initial edges of degree $r$.
\item
$[V]_1=[V]$ and $[V]_r=[V]_{r-1}\setminus[V]^{r-1}$ for
$r\in[2,t]$. That is, $[V]_r$ consists of the column vectors
corresponding to the paths with initial edges of degree $\geq r$.
\item
$G_r(R,V)$ (resp.\ $G'_r(R,V)$) is the subgroup of $G(R,V)$
(resp.\ $G'_r(R,V)$) generated by $e_v^{\lambda}$ (resp.\
$x_v^{\lambda}$) with $\lambda\in R$ and $v\in[V]_r$, where
$r\in[1,t]$.
\end{itemize}
Observe that $[V]_r$ is a rigid system for each $r\in[1,t]$. We
have the following ascending sequence of triangular subgroups
$$
G_t(R,V)\subset G_{t-1}(R,V)\subset\cdots\subset G_1(R,V)=G(R,V).
$$
Consider the mappings
$$
\tau_r:R^{N_r}\to G_r(R,P),\quad r\in[1,t],
$$
given by
$$
(\lambda_1,\dots,\lambda_{N_r})\mapsto
e^{\lambda_1}_{v_{r1}}\circ\cdots \circ
e^{\lambda_{N_r}}_{v_{rN_r}},\qquad
\{v_{r1},\dots,v_{rN_r}\}=[V]^r.
$$
and set $A_r(R,V)=\Im(\tau_r)\subset G(R,V)$, $r\in[1,t]$. Then
$A_r(R,V)\subset G_r(R,V)$ for $1\leq r\leq t$ and
$A_t(R,V)=G_t(R,V)$.

\begin{remark}\label{standard}
The name `triangular' is explained by the following observation.
Let $T_m(R)\subset\E(R)$ denote the usual triangular subgroup of
upper $m\times m$-matrices over $R$ with diagonal entries $1$.
Then it is a \Y-triangular subgroup of $\E(R)$ (viewed as
$\EE(R,\Delta_n)$ for some $n\in\NN$), supported by the simplest
graph
\begin{pspicture}(-0.25,-0.2)(3.2,0.4)
 \rput(0,0){\vertex}
 \rput(1,0){\vertex}
 \rput(2,0){\vertex}
 \rput(3,0){\vertex}
 \psline{->}(0,0)(1,0)
 \psline{->}(2,0)(3,0)
 \psline[linestyle=dotted, linewidth=0.6pt](1,0)(2,0)
\end{pspicture}.
In this situation $t=m-1$ and $G_r(R,V)$ becomes the subgroup of
$T_m(R)$ consisting of those matrices which admit non-diagonal
entries only in the rows of index $i\geq r$. The subset $A_r(R,V)$
is just the abelian subgroup of the matrices
$\prod_{j=r+1}^te_{rj}^{\lambda_j}$, $\lambda_j\in R$.

The next two theorems generalize the properties of these objects
to the polyhedral situation.
\end{remark}

\begin{theorem}\label{strutri}
Let $V$ be a rigid system.
\begin{itemize}
\item[(a)]
The mappings $\tau_r$, $r\in[1,t]$, are injective group
homomorphisms \textup{(}defined on $(R^{N_r},+)$\textup{)}.
\item[(b)]
Every element $\epsilon\in G_r(R,V)$ admits a unique
representation of type
$$
\epsilon=\epsilon_r\circ\cdots\circ\epsilon_t,\quad \epsilon_s\in
A_s(R,V),\quad s\in[r,t],
$$
which we will call the \emph{canonical representation}.
\end{itemize}
\end{theorem}

\begin{proof}
(a) We know that $[V]^r$ consists of those column vectors whose
corresponding paths in $\bF$ have an initial edge of degree $r$.
No two such paths can be multiplied in $\path\bF$. Therefore, by
Definition \ref{rigid} the product $ww'$ does not exist for any
elements $w,w'\in[V]^r$. The definition of a rigid system also
excludes that $w+w'=0$ for $w,w'\in[V]^r$. So by Theorem
\ref{elaut}(e) the mappings $\tau_r$ are in fact group
homomorphisms. The injectivity follows from the second part of the
proof of (b).

(b) Choose $\epsilon\in G_r(R,V)$ and fix a representation
$\epsilon=e_{v_1}^{\lambda_1}\circ e_{v_2}^{\lambda_2}\circ\dots$
where $v_1,v_2,\dots\in[V]_r$. Clearly, there is no loss of
generality in assuming that $v_i\in[V]^r$ for some $i$. We pick
the minimal such index. Assume $i\neq 1$. Then, using Theorem
\ref{elaut}(e) (and Definition \ref{rigid}), we see that, by
commuting $e_{v_i}^{\lambda_i}$ successively with the elementary
automorphisms $e_{v_{r-1}}^{\lambda_{r-1}}$,
$e_{v_{r-2}}^{\lambda_{r-2}}$, we can draw the factor
$e_{v_i}^{\lambda_i}$ to the left end of the new representation of
$\epsilon$. It is of course essential that $e_{v_i}^{\lambda_i}$
commutes with all the commutators, produced along the way it moves
towards the initial position in the representation. In fact, all
of these commutators belong to $A_r(R,V)$, an abelian group. In
particular, they all commute with $e_{v_i}^{\lambda_i}$ and with
each other. Next we apply the same procedure to that element of
$A_r(R,V)$ showing up first from the left in the obtained
representation etc. The crucial observation is that the iteration
of the process terminates after finitely many steps since each
nontrivial commutator corresponds to a column vector with a path
longer than those of each of its factors. But the length of paths
is bounded.

The uniqueness is shown by induction on $\#V$. For $\#V=0$ there
is nothing to prove. Assume we have shown the uniqueness for every
rigid system in $\Col(\fP)$ with $<m=\#V$ elements. Replacing the
original polytope $P$ by some element of the doubling spectrum, we
may assume $V\subset\Col(P)$.

Consider an element $\epsilon\in G_r(R,V)$ and fix a
representation
$$
\epsilon=e_{v_1}^{\lambda_1}\circ\cdots\circ
e_{v_N}^{\lambda_N},\quad N=\sum_{s=r}^tN_s
$$
in which the column vectors in $[V]^r$ appear first, then those of
$[V]^{r+1}$ etc. For simplicity we consider a fully expanded
representation in the sense that all the factors
$e_{v_i}^{\lambda_i}$, $i\in[1,N]$ are present, what we can
achieve by choosing $\lambda_i=0$ if necessary.

For an element $v\in[V]^r$ either $v_1+v\notin\Col(P)$ or the
product $v_1v$ exists. By Proposition \ref{maincrit} (and since
$\la P_{v_1},v\ra<0$ is obviously equivalent to $P_{v_1}=P_v$)
there are only two possibilities: either $\la P_{v_1},v\ra=0$ or
$P_{v_1}=P_v$. Let $1=i_1<i_2<\cdots$ and $j_1<j_2<\cdots$ be the
indices determined correspondingly by the conditions:
$P_{v_1}=P_{v_{i_2}}=P_{v_{i_3}}=\cdots$ and $\la
P_{v_1},v_{j_1}\ra=\la P_{v_1},v_{j_2}\ra=\cdots=0$. (The sequence
of the $j$ may be empty.) Let us show that
$$
\epsilon=\epsilon'\circ\epsilon'',\quad
\epsilon'=e_{v_1}^{\lambda_1}\circ e_{v_{i_2}}^{\lambda_{i_2}}
\circ e_{v_{i_3}}^{\lambda_{i_3}}\circ\cdots,\quad
\epsilon''=e_{v_{j_1}}^{\lambda_{j_1}} \circ
e_{v_{j_2}}^{\lambda_{j_2}}\circ\cdots.
$$
By Theorem \ref{elaut}(e) it is enough to show that
$v_{i_k}+v_{j_l}\notin\Col(\fP)$ for every pair of indices
$j_l<i_k$. Assume to the contrary $v_{i_k}+v_{j_l}\in\Col (\fP)$
for some $j_l<i_k$. Then by the definition of the sets $[V]^s$ and
by Proposition \ref{maincrit}(d) a product of type $v_jv_i$,
$j<i$, exists. But this is a contradiction to Proposition
\ref{maincrit}(a) because $\la P_{v_i},v_j\ra=0$.

The proper subset $V'=\{v\in[V]_r\ |\ \la
P_{v_1},v\ra=0\}\subset[V]_r$ is a rigid system. This follows from
the fact that $V'=[U]$ for a certain subset $U\subset V$ of
irreducible elements of $V$. In fact, assume $\la
P_{v_1},u_1\cdots u_k\ra=0$ for some $u_1,\dots,u_k\in[V]$
corresponding to the edges of our graph. As observed above, each
of the $u_1,\dots,u_k$ is either parallel to $P_{v_1}$ or has
$P_{v_1}$ as the base facet. What we claim is that the latter case
is impossible. Assume to the contrary that one of the
$u_1,\dots,u_k$ has the base facet $P_{v_1}$. By Proposition
\ref{maincrit}(a) this can only be $u_1$. But then Proposition
\ref{maincrit}(c) implies $P_{u_1\cdots u_k}=P_{v_1}$ -- a
contradiction with the assumption.

The restrictions of $\epsilon$ and $\epsilon''$ to the polytopal
ring $R[P_{v_1}]$ coincide and $\#V'<\#V$  (maybe $\#V'=0$).
Therefore, by the induction hypothesis the elements
$\lambda_{j_1},\lambda_{j_2},\dots\in R$ are uniquely determined
by $\epsilon$. Observe that we can apply the induction hypothesis
because the factors of $\epsilon''$ are already ordered in the
right way:
$$
\{v_{j_1},v_{j_2},\dots\}=\{v_{j_1},\dots,v_{j_a}\}\cup
\{v_{j_{a+1}},\dots,v_{j_{a+b}}\}\cup\{v_{j_{a+b+1}},\dots,v_{j_{a+b+c}}\}
\cup\dots
$$
where successive subsets belong to $A_{s'}(R,V')$ and
$A_{s''}(R,V')$ for indices $s'<s''$.

By Lemma \ref{afemb}(a) $\epsilon\circ(\epsilon'')^{-1}$ -- and,
therefore, $\epsilon$ itself -- determine uniquely the elements
$\lambda_1,\lambda_{i_2},\lambda_{i_3},\dots\in R$ as well.
\end{proof}

After the hard work has been done, we draw some consequences.

\begin{theorem}\label{strutri1}
Let $V$ be a rigid system supported by the graph $\bF$.
\begin{itemize}
\item[(a)]
For each $r\in[1,t-1]$ we have the exact sequence
$$
0\to A_r(R,V)\to G_r(R,V)\to G_{r+1}(R,V)\to0
$$
where the mapping $G_r(R,V)\to G_{r+1}(R,V)$ is determined by
$\epsilon\mapsto\epsilon_r^{-1}\epsilon$ \textup{(}notation as in
Theorem \ref{strutri}\textup{(b))}. This surjective homomorphism
is split by the identity embedding $G_{r+1}(R,V)\to G_r(R,V)$.
\item[(b)]
If $Q$ is another polytope and $W\subset \Col(Q)$ is a rigid
system supported by the same graph $\bF$, then the assignment
$e_w^\lambda\mapsto e_v^\lambda$, where $w$ and $v$ correspond to
the same path in $\bF$, gives rise to a group isomorphism
$G_r(R,W)\to G_r(R,V)$.
\item[(c)]
Let $U\subset\Col(\fP)$ be another rigid system. Then
$$
G(R,U)\cap G(R,V)=G(R,[U]\cap[V]).
$$
\item[(d)]
The natural surjective mappings $G'_r(R,V)\to G_r(R,V)$,
$r\in[1,t]$ are isomorphisms, and the assertions of Theorem
\ref{strutri} and \textup{(a)--(c)} hold analogously.
\end{itemize}
\end{theorem}

\begin{proof}
(a) follows from Theorem \ref{strutri}(b) once one observes that
the same arguments as in the first half of its proof imply the
following: for arbitrary elements $\epsilon,\epsilon'\in A_r(R,V)$
and $\rho,\rho'\in G_{r+1}(R,V)$ there exists $\epsilon''\in
A_r(R,V)$ such that
$$
\epsilon\circ\rho\circ\epsilon'\circ\rho'=\epsilon''\circ\rho\circ\rho'.
$$

(b) We have the bijective mapping $G_r(R,W)\to G_r(R,V)$ defined
in a natural way via the canonical representations. It restricts
to the assignment $e_w^\lambda\mapsto e_v^\lambda$ as above. (That
this mapping is a bijection follows from Theorem \ref{strutri}.)
In order to see that it is a group homomorphism one notices that
only the structure of the underlying graph and the commuting rules
of Theorem \ref{elaut}(e) are used in deriving the canonical
representation of an element of $G(V,r)$ or $G(W,r)$ from an
arbitrary representation (see the proof of Theorem
\ref{strutri}(b)).

(c) Without loss of generality we can assume $U,V\subset\Col(P)$.
We use induction on $\#U+\#V$ starting from the $\#U+\#V=0$ for
which there is nothing to prove.

Assume we have shown the claim when $\#U+\#V<m$ and consider the
case $\#U+\#V=m$. Pick an element $\epsilon\in G(R,U)\cap G(R,V)$
and consider canonical representations
$\epsilon=e^{\lambda_1}_{u_1}\circ e^{\lambda_2}_{u_2}\circ\cdots$
and $\epsilon=e^{\mu_1}_{v_1}\circ e^{\mu_2}_{v_2}\circ\cdots$
with respect to $U$ and $V$. Let $U'\subset U$ and $V'\subset V$
be the subsets determined correspondingly by the conditions $\la
P_{v_1},-\ra=0$ and $\la P_{u_1},-\ra=0$. Then $U'$ and $V'$ are
proper rigid subsystems.

We claim that $P_{u_1}=P_{v_1}$. As we know (from the proof of
Theorem \ref{strutri}(b)) every vector from $U$ is either parallel
to $P_{u_1}$ or has the base facet $P_{u_1}$. For $x\in\L_P$ this
implies that $\epsilon(x)$ is an $R$-linear combination of the
points $\{y\in\L_P\ |\ \la P_{u_1},y\ra\leq\la P_{u_1},x\ra\}$.
Since the same is true with respect to the facet $P_{v_1}$ the
claim follows.

Now consider the representations
$\epsilon=\epsilon'_U\circ\epsilon''_U$ and
$\epsilon=\epsilon'_V\circ\epsilon''_V$ with respect to $U$ and
$V$, similar to those in the proof of Theorem \ref{strutri}(b). We
see that
$$
\epsilon|_{R[P_{u_1}]}=\epsilon''_U|_{R[P_{u_1}]}=
\epsilon''_V|_{R[P_{v_1}]}\in G_r(R,U')\cap G_r(R,V').
$$
By the induction hypothesis one concludes
$\epsilon|_{R[P_{u_1}]}\in G(R,[U']\cap[V'])\subset
G(R,[U]\cap[V])$. (The same arguments, as at the end of the proof
of Theorem \ref{strutri}(b), show that we can use the induction
hypothesis.)

The equation $\epsilon'_U=\epsilon'_V$ and Lemma \ref{afemb}(b)
imply that $\epsilon'_U=G(R,[U]\cap[V])$. Therefore, $\epsilon\in
G(R,[U]\cap[V])$.

(d) That the mappings are isomorphisms follows from Theorem
\ref{strutri}(a) and part (a) of this theorem, with the use of
Lemma \ref{stemb} (and the induction on $\#V$). The rest is clear.
\end{proof}

\begin{remark}\label{trisimp}
It follows from Lemma \ref{colemb}(a) and Theorem
\ref{strutri1}(b) that the triangular groups associated with rigid
systems are just usual triangular matrix groups used in the
construction of the Volodin theory (see Section \ref{UNIMODULAR}).
The essential point is that in general it is not possible to
realize the whole set $\Col(P)$ with its partial product structure
in $\Col(\Delta_n)$ for any $n$. Even if this is the case (as for
the class of Col-divisible polytopes discussed in Section
\ref{Q=V}) we do not have \emph{all} triangular matrix groups used
in the classical theory.
\end{remark}

\section{Higher polyhedral $K$-groups}\label{HPK}

We now present the polytopal versions of the standard
$K$-theoretical constructions. Despite the fact that this paper
exclusively treats the case of single polytopes we use the
attribute `polyhedral' in order to indicate the possibility of a
further generalization to \emph{polyhedral algebras}, defined in
terms of \emph{lattice polyhedral complexes} \cite{BrG2}. (See
also Remark \ref{s.compl}.)

\subsection{Volodin's theory}\label{VOLOD}

Let $R$ be a ring and $P$ a balanced polytope, admitting a column
structure.

\begin{definition}\label{volod}
\leavevmode\par
\begin{itemize}
\item[(a)]
The $d$-simplices of the Volodin simplicial set $\vv(\EE(R,P))$
are those sequences
$(\epsilon_0,\dots,\epsilon_d)\in(\EE(R,P))^{d+1}$ for which there
exists a triangular group $G\in\T(R,P)$  such that
$\epsilon_k\epsilon_l^{-1}\in G$, $k,l\in[0,d]$. The $i$th face
(resp.\ degeneracy) of $\vv(\EE(R,P))$ is obtained by omitting
(resp.\ repeating) $\epsilon_i$.
\item[(b)]
The Volodin \Y-simplicial set $\vv(\EE(R,P))^\Y$ is defined
similarly using the \Y-triangular subgroups of $\EE(R,P)$. (In
particular, $\vv(\EE(R,P))^\Y\subset\vv(\EE(R,P))$ as simplicial
sets.)
\item[(c)]
The simplicial sets $\vv(\St(R,P))$ and $\vv(\St(R,P))^\Y$ are
defined analogously.
\item[(d)]
The higher Volodin polyhedral $K$-groups of $R$ are defined by
$$
K_i^{\V}(R,P)=\pi_{i-1}\bigl(|\vv(\EE(R,P))|,({\text{\bf
Id}})\bigr),\quad i\geq2,
$$
and
$$
K_i^{\V}(R,P)^\Y=\pi_{i-1}\bigl(|\vv(\EE(R,P))^\Y|,({\text{\bf
Id}})\bigr),\quad i\geq2,
$$
\end{itemize}
where $|-|$ refers to the geometric realization of a simplicial
set.
\end{definition}

It is clear that the Volodin complexes are connected. Also, the
definition of the Volodin simplicial set is independent of the
choice of $\frak P$ (see Theorem \ref{elaut}(a) and the
corresponding remarks in Subsection \ref{SCHUR}).

Below, for Volodin simplicial sets, we will usually omit $|-|$,
i.~e. we will use the same notation for a simplicial set and its
geometric realization. Also, the base points will be omitted since
they are always assumed to be the unit elements.

\subsection{The case of a unimodular simplex}\label{UNIMODULAR}

We work out the case of a unimodular simplex in detail. Volodin's
$K$-theory for a ring $R$ is defined as follows (see \cite{Su1},
\cite{Su2}, \cite{Vo}):

Suppose $\sigma$ is a partial order on $\{1,\ldots,n\}$. Define
$\T^\sigma_n(R)$ to be the subgroup of $\GL_n(R)$ consisting of
those matrices $\alpha=(a_{ij})$ for which $a_{ii}=1$ and
$a_{ij}=0$ if $i\nleq^\alpha j$. In particular, for the natural
order $\sigma:=(1<\cdots<n)$ the corresponding group
$\T^\sigma_n(R)$ is just the group of upper triangular matrices
with 1s on the diagonal.

The simplicial complex
$V_n(R):=\V(\GL_n(R),\{\T^\alpha_n(R)\}_\alpha)$ is defined in the
same way as in Definition \ref{volod} where we now take $\GL_n(R)$
instead of $\EE(R,P)$ and the $\T^\alpha_n(R)$ instead of the
triangular subgroups of Section \ref{TRIANG}. The embeddings
$\GL_n(R)\to\GL_{n+1}(R)$,
$*\mapsto\begin{pmatrix}*&0\\0&1\end{pmatrix}$, together with the
induced embeddings $\T^\sigma_n(R)\to\T^{\sigma'}_{n+1}(R)$ where
$\sigma'$ is the extension of $\sigma$ to the partial order of
$\{1,\ldots,n,n+1\}$ under which $n+1$ is the biggest element,
define embeddings of simplicial complexes $\V_n(R)\to\V_{n+1}(R)$.

Finally, for $i\ge1$ we put $K^{\V}_{i,n}(R)=\pi_{i-1}(\V_n(R))$
and $K^{\V}_i(R)=\pi_{i-1}(\V(R))=\varinjlim(K^{\V}_{i,n}(R))$
where $\V(R):=\varinjlim\V_n(R)$. The base point for the homotopy
groups is the identity matrix.

Using the same construction we could define another simplicial
complex, based on the group $\E_n(R)$ instead of $\GL_n(R)$.
Denote this complex by $\V^{\E}(R)$. Then it is clear that
$\V^{\E}(R)$ is the connected component of $\V(R)$ containing the
identity element of $\GL(R)$. Therefore,
$K^{\V}_i(R)=\pi_{i-1}(\V^{\E}(R))$ for $i\ge2$ and, of course, we
have the natural identification of $K_1(R)=\GL(R)/\E(R)$ with the
set $\pi_0(\V(R))$.

In order to see that $\V(R)$ is the same as $\vv(R,\Delta_m)$ (for
every natural number $m$) we have to recognize the triangular
subgroups of Section \ref{TRIANG} in the groups $\T^\sigma_n$
after the natural identification $\E(R)=\EE(R,\Delta_m)$. We use
the sequences $\fP'=(\Delta_m=P'_0\subset P'_1\subset\cdots)$ in
Remark \ref{quasidoubling} for the identification. Then the
vertices of a polytope in $\fP$ can be identified with the indices
used to enumerate them. Moreover, we can even assume that
$\{1,\ldots,m+1\}=\vert(\Delta_m)$ and that the new vertex of
$P_{l+1}$ is larger (as a natural number) than the vertices of
$P_l$ for $l\in\NN$.

Pick a natural number $n$ and a partial order $\sigma$ on
$\{1,\ldots,n\}$. Then $\{1,\ldots,n\}\subset\vert(P_l)$ for a
sufficiently large index $l$. Next we consider the \emph{reverse}
partial order $\sigma^{\text{op}}$ on $\{1,\ldots,n\}$. As
observed in Example \ref{basic}(c), an arbitrary partial order on
a subset of the vertices of a unimodular simplex gives rise to a
rigid system of column vectors. In particular, so does the order
$\sigma^{\text{op}}$. For the resulting rigid system $V$ we have
the equality $\T^\sigma_n(R)=G(R,V)$. It is also clear that this
process of assigning the triangular groups $G(R,V)$ to the groups
$\T^\sigma_n(R)$ can be reversed. Thus in the special case of a
unimodular simplex we recover Volodin's usual $K$-groups.

The reason that we need to pass to the reverse partial orders in
the assignment between the triangular groups of automorphisms and
the triangular groups of matrices is explained by Remark
\ref{steinberg} in Subsection \ref{ELAUT}.

Finally, we remark that the analogous claim for Quillen's theory,
introduced in Subsection \ref{QUIL} below, is just obvious and
needs no detailed explanation.

\subsection{Connection with the Milnor group}\label{CONNECTION}

The group $K_2(R,P)$ acts freely both on $\vv(\St(R,P))$  and
$\vv(\St(R,P))^\Y$ by multiplication on the right and we have
$$
\vv(\St(R,P))/K_2(R,P)=\vv(\EE(R,P)),
$$
$$
\vv(\St(R,P))^\Y/K_2(R,P)=\vv(\EE(R,P))^\Y.
$$
In fact, the equality on vertex sets follows from the very
definition of $K_2(R,P)$ and, hence, the equality for higher
dimensional simplices follows from Theorem \ref{strutri1}(d).

As in the classical case we have

\begin{lemma}\label{simpcon}
Both $\vv(\St(R,P))$ and $\vv(\St(R,P))^\Y$ are simply connected.
\end{lemma}

\begin{proof} We consider $\vv(\St(R,P))$. The arguments are completely
similar for the \Y-theory.

Clearly, we have only to show that any loop $l$ through
$1\in\vv(\St(R,P))$, consisting of {\em edges} of the simplicial
complex $\vv(\St(R,P))$, is contractible. Thus we can assume that
there are a natural number $k$, vectors $v_i\in\Col(\fP)$ and
elements $\lambda_i\in R$, $i\in[1,k]$ such that
$\prod_{i=1}^kx_{v_i}^{\lambda_i}=1$ and
$$
l=\bigl[[1,s_1],[s_1,s_2],\cdots[s_k,1]\bigr],
$$
where $s_1=x_{v_1}^{\lambda_1}$ and
$s_i=x_{v_i}^{\lambda_i}s_{i-1}$ for $i\in[1,k-1]$. For simplicity
denote this loop by $(x_1,\dots,x_k)$ where
$x_i=x_{v_i}^{\lambda_i}$.

Let ${\mathfrak F}$ denote the free (non-commutative) monoid
generated by
$$
\bigl\{|x|\ :\ x\ \text{a standard generator of}\ \St(R,P)\bigr\}.
$$
Moreover, we define the group $\mathfrak G$ as the quotient of the
free group generated by these elements modulo the relations
$$
|x^{-1}|=|x|^{-1}.
$$

For a word $w=|y_1|\cdots|y_t|\in{\frak F}$ put
$w^*=|y_t^{-1}|\cdots|y_1^{-1}|$. For words $w',w''\in\mathfrak F$
we say that $w'$ is obtained from $w''$ by {\em elementary
cancellation} if $w'=w_1w_2$ and $w''=w_1ww^*w_2$ for some (maybe
empty) words $w,w_1,w_2\in\mathfrak F$. Further, $w'$ is obtained
from $w''$ by {\em cancellation} if there is a finite sequence of
elements $w_1,\ldots,w_t$ such that $w'$ is obtained from $w_t$ by
elementary cancellation, $w_{t}$ is obtained from $w_{t-1}$ by
elementary cancellation, $\dots$, $w_1$ is obtained from $w''$ by
elementary cancellation.

The natural monoid homomorphism $\mathfrak F\to\mathfrak G$
satisfies the following condition:
\begin{itemize}
\item
$w_1,w_2\in{\mathfrak F}$ map to the same element in ${\mathfrak
G}$ if and only if there is $w_3\in\mathfrak F$ such that both
$w_1$ and $w_2$ are obtained from $w_3$ by cancellation.
\end{itemize}
Let $\mathfrak W$ denote the smallest submonoid of $\mathfrak F$
determined by the following conditions:
\begin{itemize}
\item[(i)]
the words of types
\begin{gather*}
|x_u^{-\lambda-\mu}||x_u^{\lambda}||x_u^{\mu}|,\quad
|x_{uv}^{\lambda\mu}||x_u^{\lambda}||x_v^{\mu}||x_u^{-\lambda}|
|x_v^{-\mu}|\text{ if $uv$ exists,}\\
|x_u^{\lambda}||x_v^{\mu}||x_u^{-\lambda}||x_v^{-\mu}|,
\quad\text{if $u+v\notin\Col(\fP)\cup\{0\}$},
\end{gather*}
and their $-^*$ versions are in ${\mathfrak W}$,
\item[(ii)]
if $w\in{\mathfrak W}$ then $AwA^*\in{\mathfrak W}$ for arbitrary
$A\in{\mathfrak F}$.
\end{itemize}
By the observation above the equation $x_1\cdots x_k=1$ in
$\St(R,P)$ is equivalent to the existence of words $w'\in\mathfrak
W$ and $w''\in\mathfrak F$ such that $|x_1|\cdots|x_k|$ and $w'$
are obtained from $w''$ by cancellation.

Let $w=|y_1|\cdots|y_t|$ be a word in $\mathfrak F$ such that
$y_1\cdots y_t=1$ in $\St(R,P)$. In a natural way it defines a
loop in $\vv(\St(R,P))$ consisting of edges of $\vv(\St(R,P))$. We
denote this loop by $l(w)$. Furthermore, if a word $w\in\mathfrak
F$ is a obtained from another word $w'\in\mathfrak F$ by
cancellation then they define the same element in $\St(R,P)$.

Summing up all these observations, we see that it is enough to
show the following two claims:
\begin{itemize}
\item[(i)]
if a word $w_1\in\mathfrak F$ is obtained from a word
$w_2\in\mathfrak W$ by cancellation, then $l(w_1)$ is homotopic to
$l(w_2)$;
\item[(ii)]
for every word $w\in\mathfrak W$ the loop $l(w)$ is contractible.
\end{itemize}
Now Claim (i) follows from the fact that the loop of $w_2$ only
differs from that of $w_1$ by finitely many attached ``tails'',
consisting of edges of $\vv(\St(R,P))$ -- for each of these tails
we perform forward and backward movements when we go along
$l(w_2)$. By a similar argument for Claim (ii) we only need that
the loops of types
\begin{gather*}
l(|x_u^{-\lambda-\mu}||x_u^{\lambda}||x_u^{\mu}|),\quad
l(|x_{uv}^{\lambda\mu}||x_u^{\lambda}||x_v^{\mu}||x_u^{-\lambda}|
|x_v^{-\mu}|)\text{ if $uv$ exists}\\
l(|x_u^{\lambda}||x_v^{\mu}||x_u^{-\lambda}||x_v^{-\mu}|),
\quad\text{if $u+v\notin\Col(\fP)\cup\{0\}$}
\end{gather*}
and those corresponding to the $-^*$ versions are contractible.
But due to Proposition \ref{maincrit}(d), systems in $\Col(\fP)$
of the types
$$
\quad \{u,v,uv\}\quad\text{and}\quad\{u,v\ |\ u+v\notin\Col(\fP),
\ u+v\neq 0\}
$$
are \Y-rigid. Therefore, in view of Definition \ref{volod}(c)
these loops are boundaries of simplices of $\vv(\St(R,P))$.
\end{proof}

By Lemma \ref{simpcon} $\vv(\St(R,P))$ (resp.\ $\vv(\St(R,P))^Y$)
is a universal cover of $\vv(\EE(R,P))$ (resp.\
$\vv(\EE(R,P))^Y$). Therefore we have the following

\begin{proposition}\label{coinc2}\leavevmode\par
\begin{itemize}
\item[(a)]
$K_2(R,P)=K_2^{\V}(R,P)=K_2^{\V}(R,P)^\Y$,
\item[(b)]
$K_i^{\V}(R,P)=\pi_{i-1}(\vv(\St(R,P)))$ and
$K_i^{\V}(R,P)^\Y=\pi_{i-1}(\vv(\St(R,P))^\Y)$ for all $i\geq3$.
\end{itemize}
\end{proposition}

As usual, $\B G$ denotes the classifying space of a group $G$. By
Theorem \ref{strutri1}(c) we have the following formula for rigid
systems $U,V\subset\Col(\fP)$:
\begin{equation}
\B G(R,U)\cap\B G(R,V)=\B G(R,[U]\cap[V]).\label{VolEqu1}
\end{equation}
The group $\EE(R,P)$ acts freely both on $\vv(\EE(R,P))$ and
$\vv(\EE(R,P))^\Y$ by multiplication on the right, and the
corresponding quotient spaces admit the following representations:
\begin{equation}
\X(R,P)=\bigcup_{G\in\T(R,P)}\B G\qquad \text{and}\qquad
\X(R,P)^\Y= \bigcup_{G\in\T(R,P)^\Y}\B G\label{VolEqu2}
\end{equation}
-- a general observation valid for abstract Volodin simplicial
sets associated to an arbitrary group $H$ and a system of
subgroups $\{H_\alpha\}$. Similarly, the quotient spaces of the
action of $\St(R,P)$ by multiplication on the right on
$\vv(\St(R,P))$ and $\vv(\St(R,P))^\Y$ admit the representations
\begin{equation}
\X'(R,P)=\bigcup_{G'\in\T'(R,P)}\B G'\qquad \text{and}\qquad
\X'(R,P)^\Y=\bigcup_{G'\in\T'(R,P)^\Y}\B G'.\label{VolEqu3}
\end{equation}

\begin{proposition}\label{X}
We have
\begin{itemize}
\item[(a)]
$\X(R,P)=\X'(R,P)$ and $\X(R,P)^\Y=\X'(R,P)^\Y$,
\item[(b)]
$\pi_1(\X(R,P))=\pi_1(\X(R,P)^\Y)=\St(R,P)$,
\item[(c)]
$\pi_{i-1}(\X(R,P))=K_i^{\V}(R,P)$ and
$\pi_{i-1}(\X(R,P)^\Y)=K_i^{\V}(R,P)^\Y$ for $i\geq3$.
\end{itemize}
\end{proposition}

\begin{proof}
(a) follows from Theorem \ref{strutri1}(d). As shown above, the
spaces $\vv(\St(R,P))$ and  $\vv(\St(R,P))^\Y$ are simply
connected. This implies the rest of the proposition.
\end{proof}

\subsection{Quillen's theory}\label{QUIL}

We define \emph{Quillen's higher polyhedral $K$-groups} by
$$
K_i^{\Qu}(R,P)=\pi_i(\B\EE(R,P)^+),\qquad i\geq2,
$$
where $\B\EE(R,P)^+$ refers to Quillen's $+$ construction applied
to $\B\EE(R,P)$ with respect to the whole group
$\EE(R,P)=[\EE(R,P),\EE(R,P)]$ (Theorem \ref{elaut}(b)). As
remarked in Subsection \ref{UNIMODULAR}, Quillen's polyhedral
$K$-groups coincide with the ordinary $K$-groups \cite{Qu1} when
$P$ is a unimodular simplex.

By a well known argument (see \cite{Ge}) we have the equations
\begin{equation}
K_i^{\Qu}(R,P)=\pi_i(\B\St(R,P)^+),\qquad i\geq3,\label{QuilEq}
\end{equation}
where the $+$ construction is considered with respect to the whole
group $\St(R,P)$.

We need the following general fact (see \cite{Su1},\cite{Su2}).

\begin{proposition}\label{gener}
For a group $G$ and a perfect subgroup $H\subset G$ the homotopy
fiber $Y$ of Quillen's $+$ construction $\B G\to\B G^+$ has the
following properties:
\begin{itemize}
\item[(a)]
$Y$ has the homotopy type of a CW-complex,
\item[(b)]
$Y$ is simple in dimension $\geq2$ (i.~e.\ the fundamental group
acts trivially on the higher homotopy groups),
\item[(c)]
the (reduced) singular integral homology $\tilde H_*(Y)=0$,
\item[(d)]
$Y$ is connected and $\pi_1(Y)$ is a universal central extension
of $H$,
\item[(e)]
$\pi_i(Y)=\pi_{i+1}(\B G^+)$ for $i\geq2$.
\end{itemize}
The properties \textup{(a)--(d)} characterize $Y$ up to homotopy
equivalence.
\end{proposition}

By Theorem \ref{schur}, Proposition \ref{coinc2} and Proposition
\ref{gener}(d),(e) we obtain

\begin{proposition}\label{v=q2}
$K_2^{\Qu}(R,P)=K_2(R,P)=K^{\vv}_2(R,P)^\Y=K^{\vv}_2(R,P)$.
\end{proposition}

In the next sections we will establish the coincidence of
Quillen's and Volodin's theories for all higher groups for a
certain class of balanced lattice polytopes. However we do not
know whether these theories coincide for all balanced polytopes.
As mentioned in the introduction, the strategy is as follows:
Volodin's \Y-theory has an auxiliary function -- we are able to
obtain certain acyclicity results for the spaces associated to
this theory, while the expectation is that the right theory is the
one based on arbitrary rigid systems. The acyclicity for this
theory remains an open question. Fortunately there are many
polytopes for which Volodin's complexes coincide with their
\Y-subcomplexes. As it will become clear in Section \ref{Q=V},
this problem is related to certain subtle properties of column
vectors in lattice polytopes.

\begin{lemma}\label{inorder}
 Quillen's and Volodin's theories coincide if the space
$\X(R,P)$ is acyclic (i.~e.\ the reduced integral homologies are
trivial) and simple in dimension $\geq2$.
\end{lemma}

\begin{proof}
In view of Theorem \ref{schur} and the homotopy uniqueness of $Y$
in Proposition \ref{gener}, the acyclicity and simplicity of
$\X(R,P)$ in dimension $\geq2$ identify the groups
$K^{\V}_i(R,P)=\pi_{i-1}(\X(R,P))$ (Proposition \ref{X}(c)) and
$K_i^{\Qu}(R,P)=\pi_{i-1}(Y)$ (Proposition \ref{gener}(e)) for
$i\geq3$, the case $i=2$ being settled by Proposition \ref{v=q2}.
\end{proof}

\subsection{Functorial properties}\label{FUNCT}

Let $P$ and $Q$ be balanced polytopes and $R$ a ring. If there
exists a mapping $\mu:\Col(P)\to\Col(Q)$, such that the conditions
$$
\text{(i)}\quad\la P_w,v\ra=\la
Q_{\mu(w)},\mu(v)\ra\qquad\text{and}\qquad\text{(ii)}\quad
\mu(vw)=\mu(v)\mu(w)\text{ if $vw$ exists,}
$$
hold for all $v,w\in\Col(P)$, then the assignment
$x_v^\lambda\mapsto x_{\mu(v)}^\lambda$ induces a homomorphism
$$
\St(R,\mu):\St(R,P)\to\St(R,Q).
$$
This has been proved in \cite[Proposition 9.1]{BrG5}. Moreover, if
$\mu$ is bijective, then
$$
\St(R,P)\approx\St(R,Q),\quad\EE(R,P)\approx\EE(R,Q),\quad
K_2(R,P)\approx K_2(R,P).
$$
This observation allows one to study polyhedral $K$-theory as a
functor also in the polytopal argument. The map $\mu$ is called a
\emph{$K$-theoretic morphism} from $P$ to $Q$. Though we cannot
prove $K_2$-functoriality for all maps $\mu$ (see, however,
\cite[Proposition 9.1]{BrG5} for partial results) it is useful to
note the $\St$-functoriality, since it implies bifunctoriality of
the higher polyhedral $K$-groups with covariant arguments:
\begin{multline*}
K_i^{\Qu}(-,-),K_i^{\V}(-,-):\underline{\text{\it Commutative
Rings}}
\ \times\ \underline{\text{\it Balanced Polytopes}}\to\\
\to\underline{\text{\it Abelian Groups}},\quad i\geq3.
\end{multline*}
For Quillen's theory this follows from equation (\ref{QuilEq}) in
Subsection \ref{QUIL}. For Volodin's theory one observes that the
mapping $\mu$ as above extends naturally to the column vectors in
doubling spectra $\Col(\fP)\to\Col({\mathfrak Q})$ so that the
analogous conditions are satisfied. But then the extended mapping
sends rigid systems to rigid systems. In fact, thanks to
Proposition \ref{maincrit}, for every balanced polytope $P$ one
can decide from the matrix
$$
\bigl(\la P_u,v\ra\bigr)_{u,v\in\Col(P)}
$$
when a subset $\{v_1,\ldots,v_n\}\subset\Col(P)$ defines the
product $v_1\cdots v_n$ in the strong or in the weak sense. It is,
of course, also important that $v=-w$ if and only if
$\iota(v)=-\iota(w)$. For details we refer the reader to
\cite[Section 9]{BrG5}. Now the functoriality of the groups
$K_i^{\V}(-,-)$, $i\geq3$ follows from Proposition \ref{X}(a,c).

In particular,
$$
K_i(R,P\times Q)=K_i(R,P)\oplus K_i(R,Q),\quad i\geq2,
$$
for both theories because the analogous equations hold for $\St$
and $\EE$ (\cite[Section 9]{BrG5}).

Finally, we want to point out that the $K$-theoretic groups only
depend on the projective toric variety associated with a polytope
$P$.

The \emph{normal fan} ${\cal N}(P)$  of a finite convex (not
necessarily lattice) polytope $P\subset\RR^n$ is defined as the
complete fan in the dual space $(\RR^n)^*=\Hom(\RR^n,\RR)$ given
by the system of cones
$$
\bigl(\{\phi\in(\RR^n)^*\mid\max_P(\phi)=F\},\ F\ \text{a face
of}\ P\bigr).
$$
Two polytopes $P,Q\subset\RR^n$ are called \emph{projectively
equivalent} if ${\cal N}(P)={\cal N}(Q)$.

Next we recall the relationship with projective toric varieties.
Let $P$ and $Q$ be {\it very ample} polytopes in the sense of
\cite[\S5]{BrG1}. This means that for every vertex $v\in P$ the
affine semigroup $-v+(C_v\cap\ZZ^n)\subset\ZZ^n$ is generated by
$-v+\L_P$, where $C_v$ is the cone in $\RR^n$ spanned by $P$ at
$v$ (we assume that $P\subset\RR^n$ and $\gp(S_P)=\ZZ^{n+1}$), and
similarly for $Q$. Then ${\mathcal N}(P)={\mathcal N}(Q)$ if and
only if the projective toric varieties $\Proj(k[P])$ and
$\Proj(k[Q])$ are naturally isomorphic for some field $k$. These
varieties are normal, but not necessarily projectively normal
\cite[Example 5.5]{BrG1}.

Projectively equivalent polytopes $P$ and $Q$ have the same set of
column vectors: $\Col(P)=\Col(Q)$ (see \cite{BrG1}), and the
identity map on this set is a $K$-theoretic morphism $P\to Q$.
Therefore, we have

\begin{proposition}\label{PrEq}
If $P$ and $Q$ are projectively equivalent balanced polytopes,
then $K_i^{\Qu}(R,P)\approx K_i^{\Qu}(R,Q)$ and
$K_i^{\V}(R,P)\approx K_i^{\V}(R,Q)$ for $i\geq2$.
\end{proposition}

\section{Acyclicity of $\X(R,P)^\Y$}\label{ACYCL}

In this section we follow Suslin \cite{Su1}. However, a number of
changes in Suslin's arguments \cite{Su1} are necessary. Actually,
the polyhedral constructions below do \emph{not} specialize to
those from \cite{Su1} in the classical situation of unit
simplices. In fact, a direct analogue of \cite{Su1} seems to be
impossible for general balanced polytopes.

As usual, $P$ will denote a balanced polytope, admitting a column
structure, and $\fP$ denote a doubling spectrum, starting with
$P$.

\begin{definition}\label{kdecomp}
Let $U,V\subset\Col(\fP)$ be rigid systems and $k\in\NN$. We say
that $U$ is \emph{$k$-decomposable} in $V$ if every irreducible
vector $u\in U$ admits a representation $u=v_1\cdots v_k$ with
$v_1,\dots,v_k\in V$ and, moreover, the sets of irreducible
elements of $V$ that appear in the $V$-decomposition of two
different irreducible elements of $U$ are disjoint.
\end{definition}

Clearly, if $U$ is $k$-decomposable in $V$, then $[U]$ is
$k$-decomposable in $V$.

\begin{lemma}\label{inters}
Let $k\geq2$ be a natural number and
$U_1,\ldots,U_m\subset\Col(P)$ be rigid systems for some
$m\in\NN$. Then there exist rigid systems
$V_1,\ldots,V_m\subset\Col(\fP)$ such that $U_i\subset V_i$ for
$i\in[1,m]$ and $\bigcap_{i=1}^m[U_i]$ is $k$-decomposable in
$\bigcap_{i=1}^m[V_i]$. If the $U_i$ are \Y-rigid, then also the
$V_i$ can be chosen to be \Y-rigid.
\end{lemma}

\begin{proof} Obviously, the lemma follows by an iterated use of the
following
\smallskip

\noindent\emph{Claim.}\enspace For an irreducible element
$u\in\bigcap_{i=1}^m[U_i]$ there exist rigid systems
$V_i\subset\Col(\fP)$, $i\in[1,m]$, satisfying the conditions:
\begin{itemize}
\item
$U_i\subset V_i$ for all $i$,
\item
the irreducible elements in $\bigcap_{i=1}^m[U_i]$, except $u$,
remain irreducible in $\bigcap_{i=1}^m[V_i]$,
\item
there are exactly $k$ new irreducibles, say $v_1,\ldots,v_k$, in
$\bigcap_{i=1}^m[V_i]$, belonging neither to
$\bigcap_{i=1}^m[U_i]$ nor to the affine hull of $P$, such that
$u=v_1\cdots v_k$.
\end{itemize}
(The condition that $v_1,\ldots,v_k$ do not belong to the affine
hull of $P$ yields the separation property of irreducible
elements, required in the second half of Definition
\ref{kdecomp}.)

Consider the (uniquely determined) factorizations
$$
u=u_{i1}u_{i2}\cdots u_{ir_i}\qquad i\in[1,m]
$$
where the $u_{i1},\dots,u_{ir_i}$ are irreducible elements in
$U_i$. By Proposition \ref{maincrit}(c) we have
$$
P_{u_{11}}=P_{u_{21}}=\cdots=P_{u_{m1}}.
$$
Consider the polytope
$$
Q=P^{\sq_{u_{11}}}\quad
(=P^{\sq_{u_{21}}}=\cdots=P^{\sq_{u_{m1}}}).
$$
We have $\Col(Q)\subset\Col(\fP)$ (Lemma \ref{double}(a)) and
$$
u_{i1}=\delta^+u_{i1}^\vt,\qquad i\in[1,m].
$$
By Lemma \ref{decomp} the systems
$$
W_i=U_i\cup\{\delta^+,u_{i1}^\vt\}\subset\Col(\fP),\qquad
i\in[1,m]
$$
are rigid. In particular, the products $u_{i1}^\vt u_{i2}\cdots
u_{ir_j}$, $i\in[1,m]$, exist and, clearly, they are equal. Let
$w$ denote this product.

Since neither $\delta^+$ nor $u_{i1}^\vt$ is in the affine hull of
$P$, $W_i$ has the same irreducibles as $U_i$, except that
$u_{i1}$ is replaced by the pair of new irreducibles $\delta^+$
and $u_{i1}^\vt$. By an iterative application of Lemmas
\ref{infdec} and \ref{decomp} to $\delta^+$ we can produce vectors
$\delta_1,\ldots,\delta_{k-1}$ such that
\begin{itemize}
\item
$\delta_j\notin\bigl(\text{the affine hull of}\
Q\cup\{\delta_1,\ldots, \delta_{j-1}\}\bigr)$ for $1\leq j\leq
k-1$,
\item
the sets
$V_i=U_j\cup\{\delta_1,\dots,\delta_{k-1},u_{i1}^\vt\}\subset
\Col(\fP),\ i\in[1,m]$ are rigid systems,
\item
$\delta^+=\delta_1\cdots\delta_{k-1}$.
\end{itemize}
We have
$$
\delta_1,\dots,\delta_{k-1},w\in\bigcap_{j=1}^lV_j\quad
\text{and}\quad u=\delta_1\cdots\delta_{k-1}\cdot w.
$$
By definition of the $V_i$, what remains to show is the
irreducibility of the elements $\delta_1,\ldots,\delta_{k-1}$ and
$w$ in $\bigcap_{i=1}^m[V_i]$. But for the vectors of type
$\delta$ this is obvious, and the irreducibility of $w$ is an easy
consequence of the irreducibility of $u$ in $\bigcap_{i=1}^m[U_i]$
-- one argues in terms of supporting graphs, the irreducibility
being interpreted as the condition that only the initial and
terminal points belong simultaneously to all the corresponding
paths coming from different rigid systems.
\end{proof}

The  only place in the paper where we use \Y-rigid systems
essentially is the proof of the next lemma. It is a polyhedral
translation of Suslin's arguments \cite{Su1}. The possibility of
such a translation depends heavily on Theorems \ref{strutri} and
\ref{strutri1}. We do not know how to apply these arguments to
rigid systems in general.

\begin{lemma}\label{OldStep1}
Let $U,V\subset\Col(\fP)$ be \Y-rigid systems such that $U$ is
$(k+1)$-decompo\-sable in $V$ for some $k\in\NN$. Then the
homomorphisms
$$
H_i\bigl(G(R,U),\ZZ\bigr)\to H_i\bigl(G(R,V),\ZZ\bigr),\quad
i\in[1,k]
$$
of integral homologies, induced by the embedding $[U]\subset[V]$,
are zero-maps.
\end{lemma}

\begin{proof} It is sufficient to prove that the homomorphisms
$$
H_i\bigl(G(R,U),F\bigr)\to H_i\bigl(G(R,V),F\bigr),\quad i\in[1,k]
$$
are zero-maps for arbitrary field $F$ (acted trivially by the
triangular groups). We use induction on the pairs $(k,\#[U])$ with
respect to the lexicographical order (implicitly used in
\cite{Su1}.)

Because of Theorem \ref{elaut}(b), $\Im\bigl(G(R,U)\to
G(R,V)\bigr)$ dies in the abelianization of $G(R,V)_{\text{ab}}$
whenever $U$ is 2-decomposable in $V$. In other words, our
statement is true for $k=1$ and arbitrary $\#[U]$. Therefore, we
can assume $k\geq2$. Observe that the claim is vacuously true for
arbitrary $k$ when $\#[U]=0$ -- we assume $G(R,\emptyset)=0$ by
convention.

By the induction hypothesis we only need to show that the $k$th
homology homomorphism is zero.

Let ${\bF}_U$ and ${\bF}_V$ denote the corresponding underlying
\Y-graphs. Since the sets of irreducible elements of $V$ that
appear in the $V$-decomposition of two different irreducible
elements of $U$ are disjoint (by Definition \ref{kdecomp}), one
easily observes that there is no loss of generality in assuming
that ${\bF}_V$ arises from $\bF_U$ by subdivision of each edge
into $k+1$ edges: first one considers the $\Y$ rigid subsystem of
$V$ consisting of these irreducible elements and then one changes
certain subsets of irreducibles (corresponding to suitable paths
in $\bF_V$) by the corresponding products. It is clear that the
new system will be again $\Y$-rigid. For an edge $f\in E({\bF}_U)$
the corresponding path in ${\bF}_V$ will be denoted by
$[\phi^f_1,\dots,\phi^f_{k+1}]$. We will use the notation
introduced in Section \ref{TRIANG}.

The graph ${\bF}'_V$ is defined as follows: For every $f\in
E({\bF}_U)$ we replace the path $[\phi^f_k,\phi^f_{k+1}]
\in\path{\bf F}_V$ of length 2 by a new edge $\phi^f$ from the
initial to the end point of $[\phi^f_k,\phi^f_{k+1}]$ and omit the
endpoint of $[\phi^f_k]$. These new edges together with the
remaining original edges of ${\bF}_V$ define the graph ${\bF}'_V$.

We now construct a further graph ${\bF}_V''$ which contains the
omitted vertices. Every path $l=[f_1,\dots,f_s]\in
\path_{1s}{\bF}_U$ gives rise to a system of paths
$$
\Phi^l=\{\Phi^l_1,\dots,\Phi^l_s\}\subset\path{\bF}_V
$$
where
\begin{multline*}
\qquad\Phi^l_1=[\phi^{f_1}_1,\dots,\phi^{f_1}_k],\quad
\Phi^l_2=[\phi^{f_1}_{k+1},\phi^{f_2}_1,\dots,\phi^{f_2}_k],\\
\Phi^l_3=[\phi^{f_2}_{k+1},\phi^{f_3}_1,\dots,\phi^{f_3}_k],\dots,
\Phi^l_s=[\phi^{f_{s-1}}_{k+1},\phi^{f_s}_1,\dots,\phi^{f_s}_k].
\end{multline*}
Observe that each path ends in a vertex of $\bF_V$ that was
omitted in the construction of $\bF'_V$. We take these omitted
vertices and the initial points of the paths $\Phi^l_1$ as the
vertices of $\bF''_V$ and insert edges for each of the paths
$\Phi^l_i$. Clearly, both ${\bF}_V'$ and ${\bF}_V''$ are
\Y-graphs. (The reader should observe that the \Y-property of the
involved graphs is used crucially in the definition of $\bF_V''$.)
We illustrate the construction in Figure \ref{GraphConst}.
\begin{figure}[htb]
\begin{center}
\psset{unit=0.45cm}
\def\vertex{\pscircle[fillstyle=solid,fillcolor=black]{0.15}}
\def\overtex{\pscircle{0.15}}
\begin{pspicture}(-3,0)(3,6)
 \psset{linewidth=0.6pt}
 \psline{->}(0,0)(0,3)
 \psline{->}(0,3)(-3,6)
 \psline{->}(0,3)(3,6)
 \rput(0,0){\vertex}
 \rput(0,3){\vertex}
 \rput(-3,6){\vertex}
 \rput(3,6){\vertex}
 \rput(-1.5,2){$\bF_U$}
\end{pspicture}
\quad
\begin{pspicture}(-3,0)(3,6)
 \psset{linewidth=0.6pt}
 \psline{->}(0,0)(0,1)\psline{->}(0,1)(0,2)\psline{->}(0,2)(0,3)
 \psline{->}(0,3)(-1,4)\psline{->}(-1,4)(-2,5)\psline{->}(-2,5)(-3,6)
 \psline{->}(0,3)(1,4)\psline{->}(1,4)(2,5)\psline{->}(2,5)(3,6)
 \rput(0,0){\vertex}\rput(0,1){\vertex}\rput(0,2){\vertex}
 \rput(0,3){\vertex}
 \rput(-1,4){\vertex}\rput(-2,5){\vertex}\rput(-3,6){\vertex}
 \rput(1,4){\vertex}\rput(2,5){\vertex}\rput(3,6){\vertex}
 \rput(-1.5,2){$\bF_V$}
\end{pspicture}
\quad
\begin{pspicture}(-3,0)(3,6)
 \psline[linewidth=0.6pt]{->}(0,0)(0,1)\psline[linestyle=dashed]{->}(0,1)(0,2)
      \psline[linestyle=dashed]{->}(0,2)(0,3)
 \psline[linewidth=0.6pt]{->}(0,3)(-1,4)\psline[linestyle=dashed]{->}(-1,4)(-2,5)
      \psline[linestyle=dashed]{->}(-2,5)(-3,6)
 \psline[linewidth=0.6pt]{->}(0,3)(1,4)\psline[linestyle=dashed]{->}(1,4)(2,5)
      \psline[linestyle=dashed]{->}(2,5)(3,6)
 \rput(0,0){\vertex}\rput(0,1){\vertex}\rput(0,2){\overtex}
 \rput(0,3){\vertex}
 \rput(-1,4){\vertex}\rput(-2,5){\overtex}\rput(-3,6){\vertex}
 \rput(1,4){\vertex}\rput(2,5){\overtex}\rput(3,6){\vertex}
 \psset{linestyle=solid}
 \psset{linewidth=0.6pt}
 \psarc{->}(-2,5){1.414}{-45}{135}
 \psarc{->}(0,2){1}{-90}{90}
 \psarc{->}(2,5){1.414}{215}{45}
 \rput(-1.5,2){$\bF_V'$}
\end{pspicture}
\quad
\begin{pspicture}(-3,0)(3,6)
 \psset{linestyle=dashed}
 \psline{->}(0,0)(0,1)\psline{->}(0,1)(0,2)
      \psline{->}(0,2)(0,3)
 \psline{->}(0,3)(-1,4)\psline{->}(-1,4)(-2,5)
      \psline{->}(-2,5)(-3,6)
 \psline{->}(0,3)(1,4)\psline{->}(1,4)(2,5)
      \psline{->}(2,5)(3,6)
 \rput(0,0){\vertex}\rput(0,1){\overtex}\rput(0,2){\vertex}
 \rput(0,3){\overtex}
 \rput(-1,4){\overtex}\rput(-2,5){\vertex}\rput(-3,6){\overtex}
 \rput(1,4){\overtex}\rput(2,5){\vertex}\rput(3,6){\overtex}
 \psset{linestyle=solid}
 \psset{linewidth=0.6pt}
 \psarc{->}(0,1){1}{-90}{90}
 \psline{->}(0,2)(-2,5)
 \psline{->}(0,2)(2,5)
 \rput(-1.5,2){$\bF_V''$}
\end{pspicture}
\end{center}
\caption{The construction of $\bF_V'$ and $\bF_V''$}
\label{GraphConst}
\end{figure}

There are natural embeddings $\path{\bF}_V',\path{\bf
F}_V''\subset\path{\bF_V}$. Let $V'$ and $V''$ be the \Y-rigid
subsystems of $[V]$ supported by ${\bF}'$ and ${\bF}''$, and
denote the natural embedding $G(R,U)\subset G(R,V')$ by $\iota$.

By construction every edge $\bar\phi\in E({\bF}_V'')$ corresponds
to a unique element $f\in E({\bF}_U)$ and, clearly, this
correspondence is even an isomorphism of the graphs ${\bf
F}_V''\approx{\bF}_U$.

By Theorem \ref{strutri1}(b) we have the natural group isomorphism
$G(R,U)\approx G(R,V'')$, denoted by $\psi$.

It follows from Theorem \ref{elaut}(e) that the elements of
$G(R,U)\subset G(R,V)$ and $\Im\psi=G(R,V'')\subset G(R,V)$
commute with each other. Hence we have the group homomorphism
$$
\iota\cdot\psi: G(R,U)\to G(R,V),\qquad
\bigl(\iota\cdot\psi\bigr)(g)=g\cdot\psi(g).
$$
Since the homologies are taken with coefficients in a field, by
the K\"unneth formula the induced homomorphism of the $k$th
homology groups can be decomposed as follows:
\begin{multline*}
\begin{CD}
H_k\bigl(G(R,U),F\bigr)@>{\Delta_*}>>H_k\bigl(G(R,U)\times
G(R,U),F\bigr) =
\end{CD}
\\
\begin{CD}
\bigoplus_{i+j=k}H_i\bigl(G(R,U),F\bigr)\otimes
H_j\bigl(G(R,U),F\bigr) @>{\oplus_{i+j=k}H_i(\iota,F)\otimes
H_j(\psi,F)}>>
\end{CD}
\\
\begin{CD}
\bigoplus_{i+j=k}H_i\bigl(G(R,V'),F\bigr)\otimes
H_j\bigl(G(R,V''),F\bigr)=H_k\bigl(G(R,V')\times G(R,V''),F\bigr)
\end{CD}
\\
\begin{CD}
@>{m_*}>>H_k\bigl(G(R,V),F\bigr).
\end{CD}
\end{multline*}
where $\Delta_*$ is the homomorphism induced by the diagonal
mapping and $m_*$ is induced by the multiplication in $G(R,V)$.

Now $U$ is $k$-decomposable in $V'$. Therefore, by the induction
hypothesis the homomorphisms $H_i\bigl(G(R,U),F\bigr)\to
H_i\bigl(G(R,V'),F\bigr)$ are zero for $1\leq i\leq k-1$ and the
decomposition above implies
$H_k(\iota\cdot\psi)=H_k(\iota)+H_k(\psi)$.

Consider the element
$$
e=\prod_{f\in E({\bF}_U)}e_{v(f)}^1\in G_2(R,V)
$$
where $v(f)$ is the element of $[V]$ that corresponds to the edge
$\phi^f_{k+1}\in E({\bF}_V)$ (and the product is understood as
composition).

The  following equation is proved by a routine verification on
generators (using Theorem \ref{elaut}(e)):
$$
\iota\cdot\psi=\bigl(\psi^{e}\bigr)\cdot\bigl((\iota')^e\circ\pi\bigr),
$$
where
\begin{itemize}
\item
$\iota'$ is the restriction of $\iota$ to $G_2(R,U)$,
\item
for any element $\epsilon\in G(R,V)$ we put
$\psi^e(\epsilon)=e\circ\psi(\epsilon)\circ e^{-1}$ and similarly
for $(\iota')^e$,
\item
$\pi:G(R,U)\to G_2(R,U)$ is the surjection from Theorem
\ref{strutri1}(a) ($r=1$),
\item
the dot between the two homomorphisms on the right has the same
meaning as that in $\iota\cdot\psi$. (The  images again commute,
since they the $e$-conjugates of commuting sets.)
\end{itemize}
We have
$$
H_k(\iota)+H_k(\psi)=H_k(\iota\cdot\psi)=H_k(\psi)+H_k(\iota')\circ
H_k(\pi)
$$
where the second equation is proved similarly to the first. Since
$\#[U]_2<\#[U]$ the induction hypothesis implies $H_k(\iota')=0$.
Therefore, $H_k(\iota)=0$.
\end{proof}

By Lemmas \ref{balint}, \ref{inters} and \ref{OldStep1} we get

\begin{corollary}\label{filter}
Let $U_1,\ldots,U_m\subset\Col(\fP)$ be \Y-rigid systems and $k$
be a natural number. Then there are \Y-rigid systems
$V_1,\ldots,V_m\subset \Col(\fP)$ such that $U_i\subset V_i$,
$i\in[1,m]$ and the homomorphisms
$$
H_i\biggl(G(R,\bigcap_{i=1}^m[U_i]),\ZZ\biggr)\to
H_i\biggl(G(R,\bigcap_{i=1}^m[V_i]),\ZZ\biggr),\quad i\in[1,k]
$$
are zero-maps.
\end{corollary}

It is clear that if $V_1,\ldots,V_m$ satisfy the condition of the
corollary then arbitrary rigid systems
$W_1,\ldots,W_m\subset\Col(\fP)$ with $V_i\subset W_i$,
$i\in[1,m]$ do also.

\begin{theorem}\label{acycl}
$\X(R,P)^\Y$ is acyclic.
\end{theorem}

\begin{proof}
In view of the equality (\ref{VolEqu2}) in Subsection \ref{VOLOD}
we have to show that for arbitrary \Y-rigid systems
$U_i\subset\Col(\fP)$, $i\in[1,m]$ the natural homomorphism
$$
H_k\biggl(\bigcup_{i=1}^m\B G(R,U_i),\ZZ\biggr)\to
H_k\bigl(X(R,P)^{\Y},\ZZ\bigr)
$$
is zero for every $k\in\NN$.

We show the following stronger claim in which we use the follow
notation: for a family of sets $Z_i$, $i\in[1,m]$ and a subset
$I\subset[1,m]$ we let $Z_I$ denote the intersection
$\bigcap_{i\in I}Z_i$.
\smallskip

\noindent\emph{Claim.}\enspace Let $s$ and $k$ be natural numbers
and $I_1,\ldots,I_s\subset[1,m]$ be nonempty subsets. Then there
exist \Y-rigid systems $W_1,\ldots,W_m\subset\Col(\fP)$ such that
$U_i\subset W_i$, $i\in[1,m]$ and the homomorphisms
$$
H_i\biggl(\bigcup_{j=1}^s\B G(R,[U]_{I_j}),\ZZ\biggr)\to
H_i\biggl(\bigcup_{j=1}^s\B G(R,[W]_{I_j}),\ZZ\biggr),\quad
i\in[1,k]
$$
are zero.
\smallskip

The claim gives the acyclicity as in the theorem when $s=m$ and
$I_1=\{1\}$,$\ldots$,$I_m=\{m\}$. On the other extreme, Corollary
\ref{filter} implies the claim when $s=1$.

We will use induction on $s$, the case $s=1$ being already proved.

Assume $s\geq2$. By the induction hypothesis there exist \Y-rigid
systems $U_i\subset V_i\subset W_i$, $i\in[1,m]$ such that the
upper-right vertical and lower left vertical homomorphisms in the
commutative diagram below are zero:
 {\footnotesize
$$
\begin{CD}
H_k\bigl(\bigcup\limits_{j=1}^{s-1}\B G([U]_{I_j})\bigr)\oplus
H_k\bigl(\B
G([U]_{I_s})\bigr)&\to&H_k\bigl(\bigcup\limits_{j=1}^s\B
G([U]_{I_j})\bigr)&\to&\tilde
H_{k-1}\bigl(\bigcup\limits_{j=1}^{s-1}\B G([U]_{I_j}\cap[U]_{I_s}\bigr)\\
@VVV@VVV@VVV\\
H_k\bigl(\bigcup\limits_{j=1}^{s-1}\B G([V]_{I_j})\bigr)\oplus
H_k\bigl(\B
G([V]_{I_s})\bigr)&\to&H_k\bigl(\bigcup\limits_{j=1}^s\B
G([V]_{I_j})\bigr)&\to&\tilde
H_{k-1}\bigl(\bigcup\limits_{j=1}^{s-1}\B G([V]_{I_j}\cap[V]_{I_s}\bigr)\\
@VVV@VVV@VVV\\
H_k\bigl(\bigcup\limits_{j=1}^{s-1}\B G([W]_{I_j})\bigr)\oplus
H_k\bigl(\B
G([W]_{I_s})\bigr)&\to&H_k\bigl(\bigcup\limits_{j=1}^s\B
G([W]_{I_j})\bigr)&\to&\tilde
H_{k-1}\bigl(\bigcup\limits_{j=1}^{s-1}\B G([W]_{I_j}\cap[W]_{I_s}\bigr)\\
\end{CD}
$$
 }
Here the rows represent Mayer-Vietoris sequences in which the
$\tilde H_{k-1}$ terms are being identified according to equation
(\ref{VolEqu1}). (For typographical reasons we have omitted the
ring $R$ and the coefficient group $\ZZ$.)

It follows that $W_1,\ldots,W_m$ are the desired \Y-rigid systems.
\end{proof}

\begin{question}\label{xacicl}
Is $\X(R,P)$ acyclic for every balanced polytope $P$?
\end{question}

\section{On the coincidence of Quillen's and Volodin's theories}\label{Q=V}

All polytopes are assumed to be balanced and to admit a column
vector, unless specified otherwise.

In this section we single out the class of {\it $\Col$-divisible}
polytopes and prove that Quillen's and Volodin's $K$-theories
coincide for them. The reason for the introduction of this class
of polytopes is that a rigid system of column vectors should be
embeddable into a \Y-rigid system. Fortunately,
$\Col$-divisibility persists in doubling spectra. Moreover, it is
closely related with the desired homotopy properties (simplicity)
of the relevant spaces, thus yielding the coincidence of all the
three theories $K^{\V}(-,-)$, $K^{\V}(-,-)^Y$ and $K^{\Qu}(-,-)$.
These polytopes are not so rare. For instance, all balanced
polygons are such (see Section \ref{POLYG}).

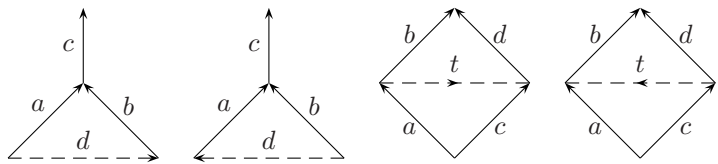
\begin{figure}[htb]
\begin{center}
\footnotesize
 \psset{unit=1cm}
\begin{pspicture}(-1,0)(1,2)
 \psline{->}(-1,0)(0,1)
 \psline{->}(1,0)(0,1)
 \psline{->}(0,1)(0,2)
 \psline[linestyle=dashed]{->}(-1,0)(1,0)
 \rput(-0.6,0.7){$a$}
 \rput(0.6,0.7){$b$}
 \rput(-0.2,1.5){$c$}
 \rput(0,0.25){$d$}
\end{pspicture}
\quad
\begin{pspicture}(-1,0)(1,2)
 \psline{->}(-1,0)(0,1)
 \psline{->}(1,0)(0,1)
 \psline{->}(0,1)(0,2)
 \psline[linestyle=dashed]{->}(1,0)(-1,0)
 \rput(-0.6,0.7){$a$}
 \rput(0.6,0.7){$b$}
 \rput(-0.2,1.5){$c$}
 \rput(0,0.25){$d$}
\end{pspicture}
\quad
\begin{pspicture}(-1,0)(1,2)
 \psline{->}(0,0)(-1,1)
 \psline{->}(0,0)(1,1)
 \psline{->}(-1,1)(0,2)
 \psline{->}(1,1)(0,2)
 \psline[linestyle=dashed](-1,1)(1,1)
 \psline[linestyle=dashed]{->}(-0.05,1)(0.05,1)
 \rput(-0.6,0.4){$a$}
 \rput(0.6,0.4){$c$}
 \rput(-0.6,1.65){$b$}
 \rput(0.6,1.65){$d$}
 \rput(0,1.25){$t$}
\end{pspicture}
\quad
\begin{pspicture}(-1,0)(1,2)
 \psline{->}(0,0)(-1,1)
 \psline{->}(0,0)(1,1)
 \psline{->}(-1,1)(0,2)
 \psline{->}(1,1)(0,2)
 \psline[linestyle=dashed](1,1)(-1,1)
 \psline[linestyle=dashed]{->}(0.05,1)(-0.05,1)
 \rput(-0.6,0.4){$a$}
\rput(0.6,0.4){$c$}
 \rput(-0.6,1.65){$b$}
 \rput(0.6,1.65){$d$}
 \rput(0,1.25){$t$}
\end{pspicture}
\end{center}
\caption{$\Col$-divisibility} \label{FigColDiv}
\end{figure}

\begin{definition}\label{divis}
A (balanced) polytope $P$ is \emph{$\Col$-divisible} if its column
vectors satisfy the conditions following:
\begin{itemize}
\item[(\cd1)]
if $ac$ and $bc$ exist and $a\neq b$, then $a=db$ or $b=da$ for
some $d$;
\item[(\cd2)]
if $ab=cd$ and $a\not=c$, then there exists $t$ such that $at=c$,
$td=b$, or $ct=a$, $tb=d$.
\end{itemize}
(See Figure \ref{FigColDiv}.)
\end{definition}

\begin{remark}\label{CDrem}
(a) It is enough in (\cd2) that $a=ct$. Then the product $tb$
exists and is necessarily equal to $d$. For the existence we note
that $\la F,t\ra >0$ for all facets $F$ with $\la F,a\ra>0$, and
so $\la P_b,t\ra>0$ since $\la P_b,a\ra>0$. (Compare Proposition
\ref{maincrit}.) Moreover, $t=-b$ is evidently impossible.

Similarly one sees that $d=tb$ is sufficient.

(b) The vector $d$ required for (\cd1) exists as soon as $a$ (or
$b$) is invertible. Then $\la P_c, -a\ra=-\la P_c, a\ra <0$, and
$P_{-a}=P_c$. But $\la P_c,b\ra >0$ by the existence of $bc$. So
$b(-a)$ exists ($b=a$ is excluded by hypothesis), and also
$(b(-a))a=b$ exists. (Again one uses Proposition \ref{maincrit}.)

The same argument (in conjunction with (a)) shows that the
invertibility of $b$ or $d$ is enough for the existence of $t$ in
(\cd2). But also the invertibility of $a$ or $c$ implies the
existence of $t$, as the reader may check. (One obtains
$a=c((-c)a)$ if $c$ is invertible.)
\end{remark}

\begin{remark}\label{affine}\label{s.compl}
(a) Not all polytopes are $\Col$-divisible. Consider the polytope
$P$ from Remark \ref{rigsyst}(b) -- the unit pyramid over the unit
square. It violates (\cd1) and (\cd2) simultaneously as one can
easily check by listing the column vectors and their products.

(b) The conditions (\cd1) and (\cd2) are independent of each
other. This is illustrated by the following examples. As already
remarked, the whole theory generalizes to lattice polyhedral
complexes (in the sense of \cite{BrG2}) whose opposite extreme
cases are single polytopes, treated here, and simplicial
complexes, viewed as lattice polyhedral complexes. (The
corresponding algebras are \emph{Stanley-Reisner} rings of
simplicial complexes.)

For simplicial complexes, for instance, a column vector is just an
oriented edge such that the facets of the complex that contain the
terminal point of the edge contain also the initial point. It is
not difficult to see that the complex $\Pi_1$ of Figure
\ref{CDIndep} satisfies (\cd2), but not (\cd1), whereas $\Pi_2$
satisfies (\cd1) and violates (\cd2).
\begin{figure}[htb]
\begin{center}
\psset{unit=1cm}
\def\vertex{\pscircle[fillstyle=solid,fillcolor=black]{0.05}}
\begin{pspicture}(-2,0)(2,2)
 \pspolygon[style=fyp](0,0)(-1,0.8)(0,2)
 \pspolygon[style=fyp, fillcolor=light](0,0)(1,0.8)(0,2)
 \psline[linestyle=dashed](-1,0.8)(1,0.8)
 \psline(-1,0.8)(-1.5,1.6)
 \psline(1,0.8)(1.5,1.6)
 \rput(0,0){\vertex}
 \rput(0,2){\vertex}
 \rput(-1,0.8){\vertex}
 \rput(1,0.8){\vertex}
 \rput(-1.5,1.6){\vertex}
 \rput(1.5,1.6){\vertex}
 \psline[style=fatline]{->}(0,2)(0,0)
 \psline[style=fatline]{->}(-1,0.8)(0,2)
 \psline[style=fatline]{->}(1,0.8)(0,2)
 \rput(-1.6,0.5){$\Pi_1$}
\end{pspicture}
\qquad\qquad
\begin{pspicture}(-2,0)(2,2)
 \pspolygon[style=fyp, fillcolor=light](0,0)(-1,0.8)(0,2)
 \pspolygon[style=fyp, fillcolor=medium](0,0)(1,0.8)(0,2)
 \psline[linestyle=dashed](-1,0.8)(1,0.8)
 \pspolygon[style=fyp](-1,0.8)(-1.5,1.6)(0,2)
 \pspolygon[style=fyp, fillcolor=lightgray](1,0.8)(1.5,1.6)(0,2)
 \psline(1,0.8)(1.5,1.6)
 \rput(0,0){\vertex}
 \rput(0,2){\vertex}
 \rput(-1,0.8){\vertex}
 \rput(1,0.8){\vertex}
 \rput(-1.5,1.6){\vertex}
 \rput(1.5,1.6){\vertex}
 \psline[style=fatline]{->}(0,2)(-1,0.8)
 \psline[style=fatline]{->}(0,2)(1,0.8)
 \psline[style=fatline]{->}(-1,0.8)(0,0)
 \psline[style=fatline]{->}(1,0.8)(0,0)
 \rput(-1.6,0.5){$\Pi_2$}
\end{pspicture}
\end{center}
\caption{The violation of (\cd1) and (\cd2)} \label{CDIndep}
\end{figure}
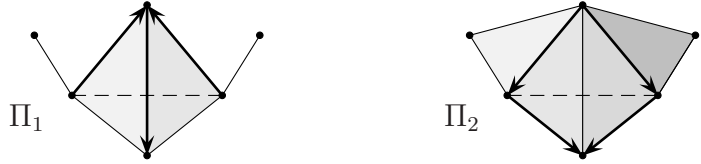
Both complexes contain a $3$-dimensional tetrahedron, and the
additional triangles of $\Pi_2$ are $2$-dimensional cells.
(Perhaps there are similar examples in the class of single
balanced polytopes.) On the other hand, it is clear that every
simplicial complex admits a $K$-theoretic morphism $\iota$ to a
unit simplex as in Proposition \ref{CDimplCE}.
\end{remark}

The next lemma contains the crucial combinatorial properties of
$\Col$-divisible polytopes. A column vector $v$ is called
\emph{terminal} if there exists no base facet $F$ with $\la
F,v\ra>0$.

\begin{lemma}\label{crucdiv}
Let $P$ be a $\Col$-divisible polytope and $u,v,w\in \Col(P)$.
\begin{itemize}
\item[(a)]
Suppose that $\la P_u,w\ra=0$ and that the product $vw$ exists.
Then $\la P_u,v\ra\le 0$.
\item[(b)] If $v$ is not terminal, then there exists
exactly one base facet $F$ with $\la F,v\ra=1$.
\item[(c)] Suppose that $v,w \in \Col(P)$ have the property that there
exist base facets $F,G$, $F\neq G$ with
$$
\la F,v\ra =\la F,w\ra=-1\quad\text{and}\quad\la G,v\ra =\la
G,w\ra =1.
$$
Then $v=w$.
\end{itemize}
\end{lemma}

\begin{proof}
(a) We have to exclude $\la P_u,v\ra>0$. If this were the case,
then the products $vu$ and $(vw)u$ would exist by Proposition
\ref{maincrit}(a). In fact, if $v+u=0$ or $v+w+u=0$, then
$-v\in\Col(P)$ (see Proposition \ref{maincrit}(h) for the second
case), and $P_{-v}$ is the only facet over which $v$ has positive
height. Since $\la P_w,v\ra >0$, one concludes that $P_w=P_u$, a
contradiction to $\la P_u,w\ra=0$. So we can assume that $vu$ and
$(vw)u$ exist.

By (\cd1) there exists $x\in\Col(P)$ such that either $xv=vw$ or
$x(vw)=v$. In the first case $x=w$, and in the second case $x=-w$.
Now the first equality is excluded by Proposition
\ref{maincrit}(d) and, in view of Proposition \ref{maincrit}(f),
the second equality implies the existence of the product $(-w)v$.
In particular, $\la P_v,-w\ra>0$ and $\la P_v, w\ra<0$, in
contradiction with the existence of $vw$.

(b) Suppose first that $-v\in\Col(P)$. Then $F=P_{-v}$ is the only
base facet with $\la F,v\ra=1.$

If $-v\notin\Col(P)$, $u\neq -v$ for each column vector $u$ such
that $F=P_u$, and the product $vu$ exists. Clearly, $\la P_u,v\ra
=1$. Consider another base facet $P_w\neq P_u$. We can also assume
$P_w\neq P_v$. If $\la P_w,u\ra =1$, then $\la P_w,v\ra =0$ since
$P$ is . So assume that $\la P_w,u\ra =0$. Then $\la P_w,v\ra =0$
by (a).

(c) Again we consider the case that $-v\in\Col(P)$ first. Then
$-v=-w$, and so $v=w$, follows from Proposition \ref{maincrit}(g).

By  symmetry we can assume that $-v\notin\Col(P)$ and
$-w\notin\Col(P)$. Then $vu$ exists for $u$ with $G=P_u$, and $wu$
also exists. By (\cd1) one of $v,w$ is divisible by the other. But
this is a contradiction because if, say, $v=tw$ then $v$ and $w$
have different base facets according to Proposition
\ref{maincrit}(c).
\end{proof}

Lemma \ref{crucdiv} makes it easy to identify column vectors and
to control the partial product structure. We extend our notation
as follows (similarly as in connection with Lemma \ref{linear}):
if $v\in\Col(P)$ is not terminal, then $P^v$ denotes the unique
base facet of $P$ with $\la F,v\ra=1$. When we write $G\neq P^v$,
this should not be interpreted as including the existence of
$P^v$. One has $\la F,v\ra=0$ for all facets $F\neq P_v,P^v$. Note
the following \emph{product rule}:
\begin{itemize}
\item The existence of the product $uv$ is equivalent to $P^u=P_v$ and
$P_u\neq P^v$, as follows immediately from Proposition
\ref{maincrit}. Moreover $P_{uv}=P_u\neq P^u=P_v$ and
$P^{uv}=P^v$, provided the latter exists.
\end{itemize}
This simple rule will save us many references to Proposition
\ref{maincrit}.

Next we want to show that for a $\Col$-divisible  polytope $P$
there exist a unit simplex $\Delta_n$ and an embedding
$\iota:\Col(P)\to\Col(\Delta_n)$ which induces an isomorphism of
partial product structures between $\Col(P)$ and the image of
$\iota$.

The critical vectors for the construction are the terminal ones.
Therefore we need some preparation. Terminal elements $u,v$ of
$\Col(P)$ are called \emph{neighbors} if there exists
$t\in\Col(P)$ with $u=tv$ or $v=tu$, and the classes of the finest
equivalence relation on the set of terminal column vectors that
respects neighbors are called \emph{cliques}. We claim:
\begin{itemize}
\item[(i)] If $u,v$, $u\neq v$, belong to the same clique, then they are
neighbors, or there exist $a,b\in\Col(P)$ and a terminal $w$ such
that $u=aw$ and $v=bw$.
\item[(ii)] Two different members of the same clique have
different base facets.
\end{itemize}
Observe that (ii) follows readily from (i): if $u=tv$, then $u$
and $v$ have different base facets  by the rule above. In the
second case, $u=aw$ and $v=bw$, we must use  (\cd1): we have
$P_u=P_a$ and $P_v=P_b$ (again by our rule) and the divisibility
of one of $a,b$ by the other forces these base facets to be
different.

For (i) it is enough to show that every chain $v_1,\dots,v_n$ of
successive, pairwise different neighbors $v_i\neq v_{i+1}$ can be
shortened if one of the following conditions is satisfied for some
$i$:
$$
\text{(1)}\ v_{i-1}=t_1v_i,\ v_i=t_2v_{i+1},\ \text{(2)}\
v_{i+1}=t_1v_i,\ v_i=t_2v_{i-1},\
 \text{(3)}\ v_i=av_{i-1},\ v_i=bv_{i+1}.
$$
In case (1) $v_{i-1}=t_1v_i=t_1(t_2v_{i+1})=(t_1t_2)v_{i+1}$: the
last product exists  by Proposition \ref{maincrit}(f) since
$t_1+t_2\neq0$. Case (2) is symmetric to (1). In case (3) we have
to invoke (\cd2): there exists $t$ such that $v_{i-1}=tv_{i+1}$ or
$v_{i+1}=tv_{i-1}$.

Now we choose a unit simplex $\Delta$ with enough facets, namely
such that it has a facet $E_F$ for each base facet $F$ of $P$, and
a facet $E_C$ for each clique $C$ of terminal column vectors $v$
of $P$. (We assume that all these facets are pairwise different.)

Define $\iota:\Col(P)\to\Col(\Delta)$ as follows:
\begin{itemize}
\item
If $v$  is not terminal then we choose $\iota(v)$ to be $e\in
\Col(\Delta)$ with $\Delta_e=E_{P_v}$ and $\Delta^e=E_{P^v}$.
\item
If $v$ is terminal,we choose $\iota(v)$ as the column vector $e$
of $\Delta$ such that $\Delta_e=E_{P_v}$ and $\Delta^e=E_C$ where
$C$ is the clique of $v$.
\end{itemize}
By Lemma \ref{crucdiv}(c) and claim (ii) above the map $\iota$ is
injective. Moreover, $\iota(v)=-\iota(-v)$ if $v,-v\in\Col(P)$,
and if $\iota(v)=-\iota(w)$, then $v=-w$. Next, if $\la
P_w,v\ra=1$, then $\la P_{\iota(w)},\iota(v)\ra=1$, and
conversely. This shows already that $vw$ exists if and only if
$\iota(v)\iota(w)$ exists. The reader may check that indeed
$\iota(vw)=\iota(v)\iota(w)$ if $vw$ exists: it is enough that
$\la E,\iota(vw)\ra=\la E, \iota(v)\ra+\la E,\iota(w)\ra$ for all
facets $E$ of $\Delta_n$. In view of Subsection \ref{FUNCT} we get

\begin{proposition}\label{CDimplCE}
Let $P$ be a  $\Col$-divisible polytope. Then there exist a unit
simplex $\Delta_n$ and an embedding
$\iota:\Col(P)\to\Col(\Delta_n)$ which defines a $K$-theoretic
morphism $\iota:P\to\Delta_n$. In particular, a system
$V\subset\Col(P)$ is rigid (\Y-rigid) if and only if $\iota(V)$ is
rigid (\Y-rigid).
\end{proposition}

We want to show that $\Col$-divisibility persists under doubling.

\begin{proposition}\label{divper}
Assume $P$ is a $\Col$-divisible polytope and $v\in\Col(P)$. Then
$P^{\sq_v}$ is $\Col$-divisible as well.
\end{proposition}

\begin{proof}
It is important for the proof that $P^{\sq_v}$ has the properties
stated in Lemma \ref{crucdiv}(b) and (c). This follows easily from
the equations (2), (3), (4) in Subsection \ref{DOUBLING}. We are
justified to use the notation $Q^v$, and the product rule
formulated above holds accordingly.

First observe that the conditions (\cd1) and (\cd2) are satisfied
if all the relevant vectors belong to either $\Col(P^-)$ or
$\Col(P^\vt)$. By symmetry between $P^-$ and $P^\vt$ the cases
listed below cover all the other possible situations. To simplify
the notation we set $Q=P^{\sq_v}$.

\noindent\emph{Condition} (\cd1). We have to show that $a$ is
divisible by $b$ or vice versa if $a\neq b$ and the products $ac$
and $bc$ exists. By Remark \ref{CDrem} we can stop checking (\cd1)
if $a$ turns out to be invertible in the case under
consideration.\smallskip

\noindent\emph{Case} (i): $c=\delta^+$. Then
$Q^a=Q^b=Q_{\delta^+}=P^\vt$. Since $Q_a, Q_b\neq
Q_{\delta^+}=P^-$, both $a$ and $b$ are column vectors of $P^-=P$,
and $P^a=P^b=P^\vt\cap P^-=P_v$. So $av$ and $bv$ exist if
$a,b\neq -v$, and the property (\cd1) of $P$ can be applied. We
may assume that $a=-v$ and can stop.\smallskip

\noindent\emph{Case} (ii): $c\in\Col(P^-)$ and $P^c=P_v$. Then
$Q_a, Q_b, Q^a, Q^b\neq Q^c=P^\vt$, and so $a,b\in\Col(P^\vt)$.

We can further assume $a\in\Col(P^\vt)\setminus\Col(P^-)$. Then
Lemma \ref{double}(c) implies $ac=\delta^-$. In particular, $a$ is
an invertible column vector, and we are done. (Use Proposition
\ref{maincrit}(h).)\smallskip

\noindent\emph{Case} (iii): $c\in\Col(P_v)$. We can assume
$a\in\Col(P^-)$. If $a\in\Col(P_v)$, then all three vectors
$a,b,c$ belong to $P^-$ or $P^\vt$, and (\cd1) for $P$ applies.
One has $P_v\neq P^a$ since $P^a=P_c\neq P_v$.

Thus we are left with the case $P_a=P_v$. If $b=a^\vt$ then
$\delta^+$ does the job: $a=\delta^+b$. So we can further assume
$b\neq a^\vt$. We have $b^-\in\Col(P)$, $b^-\neq a$ and the
products $ac$, $b^-c$ both exist. By (\cd1) there exists
$x\in\Col(P)$ such that either $a=xb^-$ or $b^-=xa$. By symmetry
we can assume that the first equation holds. Then $P_a\neq
P_{b^-},P^{b^-}$, i.~e.\ $b^-\in\Col(P_v)$. But this is equivalent
to $b\in\Col(P_v)$ and, thus, all three vectors $a,b,c$ belong to
$\Col(P)$. Done.\smallskip

\noindent\emph{Case} (iv): $c\in\Col(P^-)$ and $P_c=P_v$. In this
situation $Q^a=Q^b=P^\vt$. If neither of $a$ and $b$ is
$\delta^-$, then $a,b,c\in\Col(P^-)$ by Lemma \ref{double}(c). So
we can assume that $a=\delta^-$, and are done since $a$ is
invertible.
\medskip

\noindent\emph{Condition} (\cd2). We have to show that $ab=cd$,
$a\neq c$, implies the existence of $t$ such that $a=ct$, $d=tb$,
or $c=at$, $b=td$. By Remark \ref{CDrem} we can stop the
discussion of a case if one of $a,b,c,d$ turns out to be
invertible.\smallskip

\noindent\emph{Case} (i): $ab=\delta^+$. Then all vectors
$a,b,c,d$ are invertible, and we are done. (See Proposition
\ref{maincrit}(h).)\smallskip

\noindent\emph{Case} (ii): $ab\in\Col(P)$ and $P^\vt=Q^{ab}$.
Without loss of generality we can assume that $a\in\Col(P^\vt)$,
$P^-=Q^a$ and $b=\delta^-$.  All the other possible situations are
either symmetric to this one or reduce to the case in which all
the four vectors belong to $\Col(P)$. Since $b$ is invertible, we
are again done.\smallskip

\noindent\emph{Case} (iii): $ab\in\Col(P)$ and $P^\vt=Q_{ab}$.
Similarly to the previous case we can assume that $a=\delta^+$,
and are done.\smallskip

\noindent\emph{Case} (iv):
$ab\in\Col(P_v)=\Col(P^-)\cap\Col(P^\vt)$. Without loss of
generality we can assume $a,b\in\Col(P)$ such that $P^a=P_v$,
$P_b=P_v$, and $c,d\in\Col(P^\vt)$ such that
$(P^\vt)^c=(P^\vt)_{v^\vt}$ and $(P^\vt)_d=(P^\vt)_{v^\vt}$. (If,
say, $a\in \Col(P_v)$, then also $b\in(P_v)$, and all the four
vectors belong to $\Col(P^-)$ or $\Col(P^\vt$.)

It follows that $(P^\vt)_c=P_a$ and $(P^\vt)^c=(P_a)^\vt$. The
only vector $c$ satisfying these two conditions is
$c=a^\vt=a\delta^+$, and we are done by Remark \ref{CDrem}.
\end{proof}

Now we are ready to present the polyhedral version of Suslin's
argument for the desired simplicity.

\begin{proposition}\label{simple}
Suppose $P$ is a $\Col$-divisible polytope. Then the space
$\X(R,P)$ is simple in dimension $\geq2$.
\end{proposition}

\begin{proof} By Proposition \ref{X}(a,b) we have to show that
$\pi_1(\X'(R,P))=\St(R,P)$ acts trivially on the higher homotopy
groups of the universal cover $\vv(\St(R,P))\to\X'(R,P)$. To this
end we show that for any $z\in\St(R,P)$ the mapping
$\mu_z:\vv(\St(R,P))\to\vv(\St(R,P))$, determined by right
multiplication with $z\in\St(R,P)$, is homotopic to the identity
mapping.

In view of Proposition \ref{divper} (and the equations
$x_v^{\lambda}=[x_{\delta^+}^1,x_{v^\vt}^{\lambda}]$) it is enough
to show that for a $\Col$-divisible polytope $Q$, a vector
$v\in\Col(Q)$ and an element $\lambda\in R$ the simplicial
mappings $\vv(\St(R,Q))\to\vv(\St(R,Q^{\sq_v}))$, induced by
\begin{itemize}
\item[(i)]
the natural group homomorphism $f_v:\St(R,Q)\to\St(R,Q^{\sq_v})$,
\item[(ii)]
$f_v\cdot x_{\delta^+}^1:\St(R,Q)\to\St(R,Q^{\sq_v}),\quad
x\mapsto xx_{\delta^+}^1$,
\item[(iii)]
$f_v\cdot x_{\delta^+}^{-1}:\St(R,Q)\to\St(R,Q^{\sq_v}),\quad
x\mapsto xx_{\delta^+}^{-1}$,
\item[(iv)]
$f_v\cdot x_{v^\vt}^{\lambda}:\St(R,Q)\to\St(R,Q^{\sq_v}),\quad
x\mapsto xx_{v^\vt}^{\lambda}$,
\item[(v)]
$f_v\cdot x_{v^\vt}^{-\lambda}:\St(R,Q)\to\St(R,Q^{\sq_v}),\quad
x\mapsto xx_{v^\vt}^{-\lambda}$
\end{itemize}
are homotopy equivalent.

We denote the induced mappings between the Volodin simplicial sets
by the same names $f_v$, $f_v\cdot x^1_{\delta^+}$, $f_v\cdot
x^{-1}_{\delta^+}$, $f_v\cdot x^{\lambda}_{v^\vt}$ and $f_v\cdot
x^{-\lambda}_{v^\vt}$. Then the desired homotopies $f_v\sim
f_v\cdot x^1_{\delta^+}$, $f_v\sim f_v\cdot x^{-1}_{\delta^+}$,
$f_v\sim f_v\cdot x^{\lambda}_{v^\vt}$, and $f_v\sim f_v\cdot
x^{-\lambda}_{v^\vt}$ are given by
$$
\bigl((\underbrace{0,\dots,0}_s,\underbrace{1,\dots,1}_t)
\times(z_1,\dots,z_{s+t})\bigr)\mapsto
(z_1,\dots,z_s,z_{s+1}x,\dots,z_{s+t}x),
$$
where $x$ is correspondingly $x^1_{\delta^+}$,
$x^{-1}_{\delta^-}$, $x^{\lambda}_{v^\vt}$, or
$x^{-\lambda}_{v^\vt}$.

We only need to make sure that these homotopies do actually exist.
But this holds since the systems $V\cup\{\delta^+\}$ and
$V\cup\{v^\vt\}$ in $\Col(Q^{\sq_v})$ are rigid for every rigid
system $V\subset\Col(Q)$. In fact, by Proposition \ref{CDimplCE}
these sets with partial products are isomorphic to subsets of
$\Col(\Delta_{n+1})$ respectively of the type $W\cup\{\delta^+\}$
and $W\cup\{v^\vt\}$ for some rigid system
$W\subset\Col(\Delta_n)$. (Here $n\in\NN$ and
$\delta^+\notin\Col(\Delta_n)$). But the latter sets are obviously
rigid (Example \ref{basic}(b)).
\end{proof}

\begin{remark}\label{NecCond}
(a) We have used only the following consequence of
$\Col$-divisibility for all members $Q$ of a doubling spectrum
starting from $P$: $V\cup\{\delta^+\}$ and $V\cup\{v^\vt\}$ in
$\Col(Q^{\sq_v})$ are rigid for every rigid system
$V\subset\Col(Q)$. It may hold in a larger class of polytopes, but
we have not found a more general natural sufficient condition for
it than $\Col$-divisibility.

It is not hard to give an example for which $V\cup\{\delta^+\}$ is
not a rigid system. In fact, suppose there exist column vectors
$u,v,w$ such that $\la P_v,u\ra=\la P_w,u\ra=1$, but $P_v\neq
P_w$. (This is the case for the unit pyramid over the unit square;
see Example \ref{rigsyst}(b).) Take $V=\{u,v\}$. Then $uv$ exists,
and both $(uv)\delta_w^+$ and $u\delta_w^+$ exist. This is
impossible if $V\cup\{\delta_w^+\}$ is rigid.

(b) It  is also worth noticing that the proof above does not work
for the \Y-theory: the enlarged systems $V\cup\{\delta^+\}$ and
$V\cup\{v^\vt\}$ may no longer be \Y-rigid even if the original
system $V$ was. This and the difficulty in extending Lemma
\ref{OldStep1} to arbitrary rigid systems make it necessary to use
both versions of Volodin's theory in order to establish the
coincidence with Quillen's theory.
\end{remark}

Next we introduce the notion of \emph{$\merc$-complexity} of a
graph ${\bF}$. (Recall that the graphs we deal with are without
multiple edges, loops and they are assumed to be oriented.) The
\emph{height} $\het x$ of an element $x\in\vert{\bF}$ is by
definition the maximal possible length of an element
$l\in\path{\bF}$ having $x$ as its terminal point. We assume $\het
x=0$ if such a path does not exist.

Consider the two element subsets $\{l,l'\}\subset\path{\bF}$ such
that $l$ and $l'$ meet only at their initial and terminal points.
These couples will be called \emph{regular cycles} and we will use
the notation $l\regcyc l'$ for them. Put
$$
\het(l\regcyc l')=\het(\text{the terminal point of}\ l).
$$
Now consider the triples:
$$
\{\alpha_1,\alpha_2,\beta\in E({\bF})\ |\ \alpha_1\neq\alpha_2,\
[\alpha_1,\beta],[\alpha_2,\beta]\in \path{\bF}\}.
$$
The point $\lambda(\alpha_1,\beta,\alpha_2)$ where such edges meet
will be called a \emph{meeting point} of $\alpha_1$, $\alpha_2$
and $\beta$.

\begin{definition}\label{compl}
The $\merc$-complexity $\comp_{\merc}{\bF}$ of the graph ${\bF}$
is defined by
$$
\comp_{\merc}{\bF}=(A,B)
$$
where $A=\max\bigl(\het(l\regcyc
l'),\het(\lambda(\alpha_1,\beta,\alpha_2)\bigr)$ for
$l,l',\alpha_1,\alpha_2,\beta$ as above and $B$ is the number of
those vertices of ${\bF}$ where this maximum is achieved.

If there are no regular cycles and no meeting points the
$\merc$-complexity is $(0,0)$.
\end{definition}

\begin{lemma}\label{0compl}
If a rigid system $V\subset\Col(P)$ is supported by a graph whose
$\merc$-complexity is $(0,0)$ then $V$ is a \Y-rigid system.
\end{lemma}

\begin{proof}
If ${\bF}$ is a graph supporting $V$ and satisfying the condition
$\comp_{\merc}{\bF}=(0,0)$ then we can form a new graph
${\bF}_{Y}$ by disconnecting edges of ${\bF}$ whenever they meet
at their terminal points -- we split the terminal points of
meeting edges. Using the fact that there are no edges starting
form such meeting points, it is easily seen that ${\bF}_{Y}$ is a
\Y-graph supporting $V$.
\end{proof}

\begin{proposition}\label{resolv}
Assume $P$ is a $\Col$-divisible polytope. Then for any rigid
system $V\subset\Col(P)$ there is a \Y-rigid system
$W\subset\Col(P)$ such that $[V]\subset[W]$.
\end{proposition}

\begin{proof}
Let ${\bF}$ be a graph supporting the system $V$. Assume $l\regcyc
l'$ is a regular cycle of the maximal possible height for some
paths $l=[e_1,\dots,e_n]$ and $l'=[e'_1,\dots,e'_m]$ ($e_i,e'_j\in
E({\bF})$).

If $m=1$ or $n=1$ then we delete the corresponding edge of ${\bf
F}$.  The obtained smaller graph still supports $V$. Therefore,
without loss of generality we can assume that $m,n\geq2$. Consider
the column vectors
$$
a=v_{e_1}\cdots v_{e_{n-1}},\ b=v_{e_n},\ c=v_{e'_1}\cdots
v_{e'_{m-1}},\ d=v_{e'_m}.
$$
($v_e$ is the column vector that corresponds to the edge $e\in
E({\bF})$.) By condition (\cd2) we can assume $at=c$ and $td=b$
for some $t\in\Col(P)$.

If $t\in[V]$ then the subgraph ${\bf G}\subset{\bF}$, obtained
from ${\bF}$ by deleting the edge $e_n$, supports $V$. Clearly,
the maximal possible height of a regular cycle in ${\bf G}$ is at
most $\het(l\regcyc l')$ and the number of regular cycles in ${\bf
G}$ of this height is strictly less than the corresponding number
for the graph ${\bF}$.

Now assume $t\notin[V]$. There are two cases: either
$W=V\cup\{t\}$ is a rigid system or it is not such.

First consider the case when $W=V\cup\{t\}$ is a rigid system. In
this situation $W$ is supported by the graph ${\bf H}$ which is
obtained by deleting the edge $e_n$ and adding a new oriented edge
that connects the  terminal point of $e_{n-1}$ with the initial
point of $e'_m$. Again, any regular cycle in ${\bf H}$ has height
at most $\het(l\regcyc l')$ and the number of such cycles has
strictly decreased.

In  the remaining case when $W$ is not a rigid system, Proposition
\ref{CDimplCE} (together with the description of rigid systems in
a unit simplex, Example \ref{basic}(b),(c)) implies that
$\{w,-w\}\subset[W]$ for some $w\in\Col(P)$. Since $V$ is rigid,
the latter condition is equivalent to the condition $-t\in[V]$.
But then $d=(-t)b$. (We use the properties of a unit simplex).
Therefore, $V$ is supported by the subgraph ${\bf E}\subset{\bF}$
obtained by deleting the edge $e'_m$. Again, we can decrease the
number of the highest regular cycles.

By continuation of this process we finally reach a rigid system
$V'\subset\Col(P)$ supported by a graph without regular cycles and
$[V]\subset[V']$.

Thereafter we carry out a similar procedure to eliminate the
meeting points. This is possible due to condition (\cd1). It is
essential that we do not create regular cycles during the process.
The final result will be a rigid system $W$ supported by a graph
${\bf K}$ such that $[V]\subset[W]$ and $\comp_{\merc}{\bf
K}=(0,0)$. By Lemma \ref{0compl} we are done.
\end{proof}

\begin{theorem}\label{q=v}
Suppose $P$ is a $\Col$-divisible polytope. Then
$$
K_i^{\Qu}(R,P)=K_i^{\vv}(R,P)^\Y=K_i^{\vv}(R,P),\qquad i\geq2.
$$
\end{theorem}

\begin{proof}
By Propositions \ref{divper} and \ref{resolv} we have
$\X(R,P)=\X(R,P)^\Y$. By Proposition \ref{X}(c) we obtain the
equality of Volodin's two theories (for the case of $K_2$ one has
Proposition \ref{v=q2}), and by Lemma \ref{inorder}, Theorem
\ref{acycl} and Proposition \ref{simple} we obtain their
coincidence with Quillen's theory.
\end{proof}

A balanced polytope $P$ could be called \emph{$K$-theoretic
relative to $R$} if the equations in Theorem \ref{q=v} are
satisfied for it.

\begin{question}\label{pyramid}
Is the property of being $K$-theoretic absolute, i.~e.\
independent of the ring? Is the polytope from Remark
\ref{rigsyst}(b) $K$-theoretic? What are the corresponding groups?
Recall, that there is no nice matrix theoretical representation
available for the corresponding stable group of elementary
automorphisms (see \cite[Example 10.3]{BrG5}). It is exactly such
representations that in the polygonal case allow us to perform the
computations in the next section.
\end{question}

\section{Polygonal $K$-theories}\label{POLYG}

The class of $\Col$-divisible polytopes may at first glance seem
rather restricted. However, it follows immediately from Theorem
\ref{DIM2} that \emph{all} balanced polygons are $\Col$-divisible.
In particular, balanced polygons are $K$-theoretic. Classification
of the $\Col$-divisible polytopes is an interesting problem
already in dimension 3.

In Theorem \ref{DIM2} we have grouped all balanced polygons in six
infinite series which, according to \cite[Theorem 10.2]{BrG5},
give rise to the following isomorphism classes of stable
elementary automorphism groups:
\begin{equation}
\tag{a}\EE_a=\E(R),
\end{equation}
\begin{equation}
\tag{b}\EE_b=
\begin{pmatrix}
\E(R)&\End_R(\oplus_{\NN}R)\\
&\\
0&\E(R)
\end{pmatrix},
\end{equation}
\begin{equation}
\tag{c}\EE_c=
\begin{pmatrix}
\E(R)&\End_R(\oplus_{\NN}R)&
\Hom_R(\oplus_{\NN}R,R)\\
&\\
0&\E(R)&\Hom_R(\oplus_{\NN}R,R)\\
&\\
0&0&1
\end{pmatrix},
\end{equation}
\begin{equation}
\tag{d} \EE_{d,t}=
\begin{pmatrix}
\E(R)&\Hom_R(\oplus_{\NN}R,R^t)\\
0&\text{\bf Id}_t
\end{pmatrix},\quad t\in\NN,
\end{equation}
\begin{equation}
\tag{e}\EE_e=\E(R)\times\E(R),
\end{equation}
\begin{equation}
\tag{f}\EE_f=
\begin{pmatrix}
\E(R)&\Hom_R(\oplus_{\NN}R,R)\\
&\\
0&1
\end{pmatrix}
\times
\begin{pmatrix}
\E(R)&\Hom_R(\oplus_{\NN}R,R)\\
&\\
0&1
\end{pmatrix}.
\end{equation}
In view of the remarks above the following theorem identifies all
possible polygonal $K$-groups (under a technical restriction on
rings).

\begin{definition}[\cite{NSu}]\label{muni}
A not necessarily commutative ring $A$ is an \emph{$S(n)$-ring} if
there are $a_1,\dots,a_n\in A^*$ such that the sum of each
nonempty subfamily is a unit. If $A$ is an $S(n)$-ring for all
$n\in\NN$, then $A$ has \emph{many units}.
\end{definition}

The class of rings with many units includes local rings with
infinite residue fields and algebras over rings with many units.

\begin{theorem}\label{comput}
For every (commutative) ring $R$ and every index $i\geq2$ we have:
\begin{itemize}
\item[(a)]
$\pi_i(\B\EE_a^+)=K_i(R)$,
\item[(b)]
$\pi_i(\B\EE_b^+)=K_i(R)\oplus K_i(R)$,
\item[(c)]
$\pi_i(\B\EE_c^+)=K_i(R)\oplus K_i(R)$ if $R$ has many units,
\item[(d)]
$\pi_i(\B\EE_{d,t}^+)=K_i(R)$ if $R$ has many units,
\item[(e)]
$\pi_i(\B\EE_e^+)=K_i(R)\oplus K_i(R)$,
\item[(f)]
$\pi_i(\B\EE_f^+)=K_i(R)\oplus K_i(R)$ if $R$ has many units.
\end{itemize}
\end{theorem}

\begin{proof}
Let $\GG_a,\GG_b,\dots,\GG_f$ denote the groups, obtained by the
corresponding substitution of $\GL(R)$ for $\E(R)$ in the groups
$\EE_a,\EE_b,\dots,\EE_f$. Then we have the equations
$\EE_a=[\GG_a,\GG_a]$,
$\EE_b=[\GG_b,\GG_b]$,$\dots$,$\EE_f=[\GG_f,\GG_f]$.
 We have $\pi_i(\B\GG_a^+)=\pi_i(\EE_a^+)$,
$\pi_i(\B\GG_b^+)=\pi_i(\EE_b^+)$,$\dots$,$\pi_i(\B\GG_f^+)=\pi_i(\EE_f^+)$,
where the + constructions $\B\GG_a^+$ etc. are considered relative
to the normal subgroups $\EE_a\subset\GG_a=\pi_1(\B\GG_a)$ etc.
which are also perfect. (See, for instance, \cite[Theorem
5.2.7]{Ro}.)

Since we have $\B\GG_f^+\approx\B\GG_{d,1}^+\times\B\GG_{d,1}^+$
(essentially due to the uniqueness of the + construction) it is
enough to show that
$$
\B\GL(R)^+\times\B\GL(R)^+\approx\B\GG_b^+\approx\B\GG_c^+.
$$
and
$$
\B\GL(R)^+\approx\B\GG_{d,t}^+,\quad t\in\NN.
$$
Because of the equations
\begin{align*}
 H_1(\GG_b,\ZZ)&=(\GG_b)_{\text{ab}}=\pi_1(\GG_b^+),\\
 H_1(\GG_c,\ZZ)&=(\GG_c)_{\text{ab}}=\pi_1(\GG_c^+),\\
 H_1(\GG_{d,1},\ZZ)&=(\GG_{d,1})_{\text{ab}}=\pi_1(\GG_{d,1}^+)
\end{align*}
in conjunction with Whitehead's theorem it is sufficient to
establish that
\begin{align*}
\tag{i}H_i(\GG_b,\ZZ)&=H_i(\GL(R)\times\GL(R),\ZZ),\quad
 i\in\NN,\\
\tag{ii} H_i(\GG_c,\ZZ)&=H_i(\GL(R)\times\GL(R),\ZZ),\quad
i\in\NN,\\
\tag{iii} H_i(\GG_{d,t},\ZZ)&=H_i(\GL(R),\ZZ),\quad i,t\in\NN.
\end{align*}
Now, (i) is proved in \cite{Qu2} (for not necessarily commutative
rings), and the stronger unstable version of (iii) (for not
necessarily commutative) rings with many units is proved in
\cite{NSu}. It  only remains to notice that the validity of (i)
and (iii) for not necessarily commutative rings implies (ii) as
follows. Put
$$
T=
\begin{pmatrix}
R&R\\
0&R
\end{pmatrix}.
$$
Then $\GL(T)\approx\GG_b$ and the following natural {\it split}
epimorphisms give the result
\begin{multline*}
\begin{pmatrix}
\GL(T)&\Hom(\oplus_{\NN}T,T)\\
&\\
0&1
\end{pmatrix}\approx\\
\displaybreak[3]\\
\approx
\begin{pmatrix}
\GL(R)&\Hom(\oplus_{\NN}R,\oplus_{\NN}R)&\Hom(\oplus_{\NN}R,R)&
\Hom(\oplus_{\NN}R,R)\\
\\
0&\GL(R)&0&\Hom(\oplus_{\NN}R,R)\\
\\
0&0&1&0\\
\\
0&0&0&1
\end{pmatrix}\to
\GG_c\to\\
\\
\to\GL(R)\times\GL(R)
\end{multline*}
where the units refer to the unit elements respectively in $T$ and
in $R$, the second homomorphism is obtained by erasing the third
column and third row, and the third homomorphism is obtained by
picking the first two diagonal entries. (The condition of the
existence of many units is inherited by $T$.)
\end{proof}


\begin{thebibliography}{[ConW]}
\bibitem[BrG1]{BrG1}
W. Bruns and J. Gubeladze, {\em  Polytopal linear groups}, J.
Algebra {\bf 218} (1999), 715--737.

\bibitem[BrG2]{BrG2}
W. Bruns and J. Gubeladze, {\em Polyhedral algebras, arrangements
of toric varieties, and their groups}, Advanced Studies in Pure
Mathematics  {\bf 33} (2002), Computational Commutative Algebra
and Combinatorics, pp.\ 1--51.

\bibitem[BrG3]{BrG3}
W. Bruns and J. Gubeladze, {\em Polytopal linear retractions},
Trans. Amer. Math. Soc. {\bf 354} (2002), 179--203.

\bibitem[BrG4]{BrG4}
W. Bruns and J. Gubeladze, {\em Polytopal linear algebra}, Beitr.
Algebra Geom. {\bf 43} (2002), 479--500.

\bibitem[BrG5]{BrG5}
W. Bruns and J. Gubeladze, {\em Polyhedral $K_2$}, Manuscr. Math.
{\bf 109} (2002), 367--404.

\bibitem[BrG6]{BrG6}
W. Bruns and J. Gubeladze, {\em Normality and covering properties
for affine semigroups}, J. Reine Angew. Math. {\bf 510} (1999),
161--178.

\bibitem[BrGTr]{BrGTr}
W. Bruns, J. Gubeladze and N. V. Trung, {\em Normal polytopes,
triangulations, and Koszul algebras}, J. Reine Angew. Math. {\bf
405} (1997), 123--160.

\bibitem[Ge]{Ge}
 S. M. Gersten, {\em $K_3$ of a ring is $H_3$ of the
Steinberg group}, Proc. Amer. Math. Soc. {\bf 37} (1973),
366--368.

\bibitem[J]{J}
J. P. Jouanolou, {\em Une suite exacte de Mayer-Vietoris en
K-theorie algebrique}, Lecture Notes Math. {\bf 341}, Springer,
(1973), 293--316.

\bibitem[NSu]{NSu}
Y. Nesterenko and A. Suslin, {\em Homology of the general linear
group over a local ring, and Milnor's $K$-theory}, (Russian) Izv.
Akad. Nauk SSSR Ser. Mat. {\bf 53} (1989), no. 1, 121--146.

\bibitem[Qu1]{Qu1}
D. Quillen, {\em On the cohomology and $K$-theory of the general
linear groups over a finite field}, Ann. of Math. (2) {\bf 96}
(1972), 552--586.

\bibitem[Qu2]{Qu2}
D. Quillen, {\em Characteristic classes of representations.
Algebraic $K$-theory}, Lecture Notes in Math. {\bf 551}, Springer,
Berlin, (1976), 189--216.

\bibitem[Ro]{Ro}
J. Rosenberg, {\em Algebraic $K$-theory and its applications},
Graduate Texts in Mathematics {\bf 147}, Springer-Verlag, 1994.

\bibitem[Su1]{Su1}
A. Suslin, {\em On the equivalence of $K$-theories}, Comm. Algebra
{\bf 9} (1981), 1559--1566.

\bibitem[Su2]{Su2}
A. Suslin, {\em Stability in algebraic $K$-theory}, Algebraic
$K$-theory, Part I, Lecture Notes in Math. {\bf 966} (1982),
Springer, 304--333.

\bibitem[Vo]{Vo}
L. Volodin, {\em Algebraic $K$-theory as an extraordinary homology
theory on the category of associative rings with a unit},
(Russian) Izv. Akad. Nauk SSSR Ser. Mat. {\bf 35} (1971),
844--873.

\bibitem[Z]{Z}
G. Ziegler, {\em Lectures on Polytopes}, Graduate Texts in
Mathematics {\bf 152}, Springer-Verlag, 1995, Revised edition
1998.

\end{thebibliography}
\end{document}